\documentclass[11pt,english]{amsart}
\usepackage[T1]{fontenc}
\usepackage{bm}
\usepackage{amstext}
\usepackage{amsthm}
\usepackage{amssymb}
\usepackage{setspace}
\usepackage{tikz}
\usetikzlibrary{arrows}
\usepackage[a4paper,margin=0.9in]{geometry}
\usepackage{hyperref}
\usepackage[nameinlink,noabbrev]{cleveref}
\hypersetup{colorlinks=true,linkcolor=blue,citecolor=blue,urlcolor=blue}
\usepackage{enumitem}
\usepackage{babel}

\setlist{leftmargin=2.2em,labelsep=0.5em,itemsep=0em,parsep=0em}
\numberwithin{equation}{section}
\numberwithin{figure}{section}

\allowdisplaybreaks[4]

\theoremstyle{plain}
\newtheorem{theorem}{Theorem}[section]
\newtheorem{proposition}[theorem]{Proposition}
\newtheorem{corollary}[theorem]{Corollary}
\newtheorem{lemma}[theorem]{Lemma}
\theoremstyle{remark}
\newtheorem{remark}[theorem]{Remark}
\theoremstyle{definition}

\newtheorem*{model}{The model}

\newcommand{\md}{\mathop{}\mathrm{d}}
\newcommand{\R}{\mathbb{R}}
\newcommand{\E}{\mathbb{E}}
\newcommand{\N}{\mathbb{N}}

\newcommand{\Law}{\mathrm{Law}}
\newcommand{\loc}{\mathrm{loc}}
\newcommand{\supp}{\mathrm{supp}}
\newcommand{\Pp}{\mathcal{P}}
\newcommand{\ee}{\mathrm{e}}
\newcommand{\odotimes}{\mathbin{\odot}}
\renewcommand{\S}{\mathbb{S}}

\begin{document}
\title[Propagation of chaos for the Landau equation]{Propagation of chaos for the Landau equation via microcanonical binary collisions}
\author{Kai Du}
\date{\today}
\address{Shanghai Center for Mathematical Sciences, Fudan University, Shanghai 200433, China.}
\email{kdu@fudan.edu.cn}

\begin{abstract}
We develop a fully constructive, conservative, and collision-level realization of Kac's program for the spatially homogeneous Landau equation across the full interaction range, including the Coulomb case.
Our model is the microcanonical binary-collision (MBC) process: a reversible pure-jump $N$-particle Markov process, realizing the grazing-collision mechanism via small conservative rotations of relative velocities.
Under suitable assumptions on the initial data, we prove propagation of chaos in the joint mean-field and grazing-collision limit ($N\to\infty, h\to0$), identifying any limit point with the unique solution in the relevant weak-solution class.
The empirical measures converge in quadratic Wasserstein distance, while fixed-order marginals converge strongly in $L^1$ at every fixed time.
Furthermore, we rigorously derive the Landau master equation as the fixed-$N$ grazing-collision limit of the MBC process.
The proof combines Fisher-information dissipation, compactness, and uniqueness; for $\gamma\in[-3,-2)$, second-order information along collision directions yields the time-integrable $L^\infty$ bound needed for uniqueness.
\end{abstract}

\maketitle
\tableofcontents

\section{Introduction\label{sec:Introduction}}
\label{sec:introduction}

The spatially homogeneous Landau equation is a fundamental kinetic model for weakly coupled plasmas in the grazing-collision regime \cite{Landau1936,Villani2002}. It describes the evolution of a velocity distribution $f_t : \R^{3} \to [0,\infty)$ through
\begin{equation}
\partial_{t}f_{t}(v)=\frac{1}{2}\,\nabla_{v}\cdot\int_{\R^{3}}a(v-w)\big(f_{t}(w)\nabla_{v}f_{t}(v)-f_{t}(v)\nabla_{w}f_{t}(w)\big)\md w,\label{eq:Landau}
\end{equation}
where the collision kernel $a(\cdot)$ is given by
\begin{equation}
a(z)=|z|^{\gamma+2}\,(I-\hat{z}\otimes\hat{z}),\quad\hat{z}=z/|z|,\quad\gamma\in[-3,1].\label{eq:A-kernel}
\end{equation}
The physically most relevant $\gamma=-3$ corresponds to the original Landau equation, 
obtained formally from the Boltzmann equation with Coulomb potential. 
The global well-posedness theory for \eqref{eq:Landau} has been the subject of a decades-long effort, facing particular challenges in the singular Coulomb case; see for instance \cite{Villani1998ARMA,DesvillettesVillani2000I,DesvillettesVillani2000II,FournierGuerin2009JFA,Fournier2010CMP} and references therein. 
A recent breakthrough by Guill\'en--Silvestre \cite{GuillenSilvestre2025}  has transformed this picture by identifying the Fisher information as a Lyapunov functional which is non-increasing along sufficiently decaying classical solutions of \eqref{eq:Landau}; 
this monotonicity rules out blow-up for initial data with a Maxwellian upper bound and opens a new path to global regularity (see~\cite{Villani2025FIKT,ImbertSilvestreVillani2024} for broader perspectives and further developments). 
As we will demonstrate, this functional, and the geometric structure that underpins its dissipation, constitutes the analytic core of our microscopic derivation.

A central theme, originating with Boltzmann~\cite{Boltzmann1872}, is the rigorous derivation of kinetic equations from underlying Newtonian many-particle dynamics. This line of research, highlighted by Lanford's classical theorem \cite{Lanford1975_TimeEvolution} establishing the Boltzmann--Grad limit from hard-sphere dynamics for short times, has recently seen a breakthrough by Deng--Hani--Ma \cite{DengHaniMa2025_longtimeBoltzmann}, who extended this derivation to the entire time interval of existence of the corresponding regular Boltzmann solution. 
These results underscore the central role of microscopic derivations and provide a demanding benchmark for particle-based models in kinetic theory.

Alternatively, Kac~\cite{Kac1956Foundations} proposed his seminal program to derive nonlinear kinetic equations from stochastic many-particle systems rather than from Newtonian dynamics. 
He introduced a \emph{Poisson-like} binary-collision
process on the energy sphere: at random times, a uniformly chosen pair of
scalar velocities is rotated, preserving total energy and the microcanonical measure.
This reversible dynamics provides a microscopic model for a one-dimensional caricature of the spatially
homogeneous Boltzmann equation and yields propagation of chaos.
A decisive
advance was achieved by Mischler--Mouhot~\cite{MischlerMouhot2013Invent},
who developed a robust ``top-down'' framework and obtained quantitative,
uniform-in-time propagation of chaos estimates for realistic unbounded kernels,
including the hard-sphere model (a canonical strong-collision kernel). 
Their analysis, which exploits sharp stability properties of the Boltzmann flow,
sets a natural benchmark for stochastic microscopic derivations of kinetic equations.

In contrast, implementing Kac's program for the Landau equation, and in
particular for the singular Coulomb case $\gamma=-3$, has remained a delicate
problem. 
The main difficulty is the near-diagonal singularity of the Landau
collision operator, which reflects the cumulative effect of numerous weak,
grazing binary interactions. 
Recent important works~\cite{FengWang2025_KacLandau,Tabary2025_PoC_Landau_Coulomb} have established
propagation of chaos for the Landau equation by taming the $N$-particle
Landau master equation or mean-field diffusion models. 
These approaches, however, start from diffusion-type generators and are conceptually distinct from Kac's original blueprint based on discrete binary collisions.

Constructing such a process and proving its chaotic limit is governed by \eqref{eq:Landau} is the central aim of this work.

\begin{model}
For $N\ge2$, the $N$-particle \emph{microcanonical binary-collision} (MBC) process
\[
\bm{v}^{N}(t)=(v_{1}^{N}(t),\dots,v_{N}^{N}(t)),\quad t\ge0
\]
is a \emph{pure-jump} continuous-time Markov process on $\R^{3N}$
defined as follows. Fix a step $h>0$ and let jumps occur at Poisson
times with rate $\lambda=N/h$. At each jump, choose an unordered
pair $\{i,j\}\subset\{1,\dots,N\}$ uniformly from all $\binom{N}{2}$
pairs. If $v_{i}\neq v_{j}$, write the relative direction $e_{ij}=(v_{i}-v_{j})/|v_{i}-v_{j}|\in\mathbb{S}^{2}$,
evolve $e_{ij}$ by spherical Brownian motion with generator $\Delta_{\mathbb{S}^{2}}$
for (virtual) time 
\[
\tau_{ij}=h|v_{i}-v_{j}|^{\gamma},
\]
obtaining $e_{ij}'$, and then update 
\[
v_{i}'=\frac{v_{i}+v_{j}}{2}+\frac{|v_{i}-v_{j}|}{2}\,e_{ij}',\qquad v_{j}'=\frac{v_{i}+v_{j}}{2}-\frac{|v_{i}-v_{j}|}{2}\,e_{ij}',
\]
leaving all other coordinates unchanged (do nothing if $v_{i}=v_{j}$).
Each jump is a rigid rotation in the two-body subspace, hence the
total momentum $\sum_{i}v_{i}$ and energy $\frac{1}{2}\sum_{i}|v_{i}|^{2}$
are conserved along the entire trajectory. 
Equivalently, every sample path remains on the corresponding microcanonical shell (Boltzmann sphere),
and the microcanonical measure is invariant and reversible.
\end{model}

Our main result is the following.
\begin{theorem}
	\label{thm:main}
	Let $\gamma\in[-3,1]$ and $h=h_{N}\to0$ as $N\to\infty$.
	Assume the initial law $f_{0}\in\mathcal{P}_2(\R^{3})$
	has finite Fisher information
	\[
	I(f_{0}):=\int_{\R^{3}}\big|\nabla\log f_{0}(v)\big|^{2}f_{0}(v)\md v<\infty,
	\]
	and a finite $m$-th moment with 
	\[
	m>3(2-\gamma)\ \text{if }\gamma\in[-3,-2),\quad
	m=2\ \text{if }\gamma\in[-2,0),\quad
	m=4\ \text{if }\gamma\in[0,1].
	\]
	Let $(F_{t}^{N})_{t\ge0}$ be the law of the MBC process started from
	$f_{0}^{\otimes N}$, and set $\mu_{t}^{N}=\frac{1}{N}\sum_{i=1}^{N}\delta_{v_{i}^{N}(t)}$.
	
	Then there exists a global weak solution $(f_t)_{t\ge 0}$ of the Landau equation~\textnormal{(\ref{eq:Landau})} with initial data $f_0$, enjoying uniform-in-time bounds on its second moment and Fisher
information, 
such that

1. For every $T>0$,
\[
  \lim_{N\to\infty} \E \Big[\sup_{0\le t\le T} W_2(\mu_t^N,f_t)\Big] = 0,
\]
where $W_2$ is the quadratic Wasserstein distance on $\mathcal{P}_2(\R^3)$.

2. For every $k\ge1$ and $t\ge0$, the $k$-particle marginal $F_t^{k|N}$ admits a Lebesgue density and
  \[
    \lim_{N\to\infty}
    \big\|F_t^{k|N} - f_t^{\otimes k}\big\|_{L^1(\R^{3k})} = 0.
  \]
  In particular, the system is $f_t$-chaotic at each time in the strong
  $L^1$ (total-variation) sense.
\end{theorem}

As described above, the MBC process is a reversible Kac-type binary-collision system with exact conservation of momentum and energy.
Its design realizes the grazing-collision mechanism: at each jump a pair interacts via a small rotation over a virtual time $\tau_{ij}=h|v_i-v_j|^{\gamma}$.
For a fixed pair with relative speed $r=|v_i-v_j|>0$, the conditional mean-square displacement satisfies
\[
 \E\big[|\Delta v_i|^2\mid v_i,v_j\big]
 =\frac{r^2}{2}\big(1-\ee^{-2hr^\gamma}\big)
 =hr^{\gamma+2}+o(h)
 \quad\text{as }h\downarrow0.
\]
In the regime $N\to\infty$, $h\to0$ with jump rate $\lambda=N/h$, these tiny deflections accumulate to produce the effective quadratic variation $|v_i-v_j|^{\gamma+2}$ of the Landau operator \eqref{eq:A-kernel} across the full interaction range $\gamma\in[-3,1]$.

Theorem~\ref{thm:main} makes this heuristic precise: under natural Fisher information and moment assumptions, the MBC dynamics propagates chaos and its empirical measure converges to a solution of the Landau equation in the standard weak sense. Combined with existing uniqueness results (e.g.~\cite{FournierGuerin2009JFA,Fournier2010CMP,FournierHeydecker2021}), this identifies the Landau equation as the kinetic description of this reversible microcanonical binary-collision system in the sense of Kac's program, providing a collision-level realization of Landau dynamics from conservative binary interactions.

The proof is completed in Section~\ref{sec:uniqueness-completion}, using Proposition~\ref{vs:thm:closure}, whose proof is given in Section~\ref{sec:very-soft-closure}.
Finite entropy follows from the finite energy and Fisher-information: Sobolev's inequality gives $f_0\in L^3$ and controls the positive part of $f_0\log f_0$, while the energy controls its negative part.

\begin{remark}[Solution classes]\label{rem:uniqueness-classes}
On each finite interval $[0,T]$, we use uniqueness of weakly continuous, energy-conserving weak solutions with initial datum $f_0$ in the following classes:
\begin{itemize}
\item \emph{Hard potentials} $\gamma\in(0,1]$: $L^\infty(0,T;\Pp_4(\R^3))$; see \cite{FournierHeydecker2021}.
\item \emph{Maxwell molecules} $\gamma=0$: finite energy suffices; see \cite{Villani1998MMMAS}.
\item \emph{Moderately soft potentials} $\gamma\in[-2,0)$: $\int_0^T I(f_t)\md t<\infty$ suffices, by Hardy's inequality and \cite[Theorem~1.3]{FournierGuerin2009JFA}.
\item \emph{Very soft and Coulomb potentials} $\gamma\in[-3,-2)$: $L^\infty(0,T;\Pp_2(\R^3))\cap L^1(0,T;L^\infty(\R^3))$;
see \cite{FournierGuerin2009JFA} for $\gamma>-3$ and \cite{Fournier2010CMP} for $\gamma=-3$.
\end{itemize}
\end{remark}

\subsection{Relation to existing works}\label{sec:related-work}

The rigorous derivation of the Landau equation from many-particle models has motivated a broad and active literature; we do not aim at completeness and refer to the surveys and introductions in, e.g.,~\cite{Villani2002,FournierHauray2016,FournierGuillin2017,FengWang2025_KacLandau,Tabary2025_PoC_Landau_Coulomb} for a more extensive account.
To clarify the positioning of our work, we organize the discussion around the nature of the underlying $N$-particle system:
\begin{enumerate}
	\item PDE-level grazing-collision limits from Boltzmann to Landau;
	\item Propagation of chaos for particle systems whose generators are \emph{a priori} diffusion operators (Landau-type master equations or interacting SDEs);
	\item Derivations from conservative, pure-jump processes built from discrete binary collisions, in the spirit of Kac's original program.
\end{enumerate}
Our work falls into the third category and, as we explain below, also provides a constructive microscopic link to the diffusive models in (2).

\vspace{0.8ex}
\textbf{PDE-level grazing limits.}
At the level of kinetic PDEs, the Landau equation \eqref{eq:Landau} is classically understood as the grazing-collision limit of the Boltzmann equation; see~\cite{Villani2002} and references therein.
These results justify the Landau operator as a macroscopic description of a regime dominated by weak, small-angle deflections, but they do not by themselves provide a fully microscopic, particle-level derivation.

\vspace{0.8ex}
\textbf{Models based on $N$-particle diffusion generators.}
A major line of work starts from an $N$-particle system whose generator already has a diffusive Landau-type structure.

$\bullet$ \emph{Landau master equation.}
A widely used microscopic object in weak-coupling and grazing-collision programs is the Landau master equation with formal generator
\begin{equation}\label{eq:LM-generator}
	\mathcal{L}^{N}F(\bm v)
	= \frac{1}{4(N-1)}\sum_{i\neq j}
	(\nabla_{v_i}-\nabla_{v_j})\cdot\Big(
	a(v_i-v_j)\,(\nabla_{v_i}-\nabla_{v_j})F(\bm v)\Big),
\end{equation}
for $\bm v=(v_1,\dots,v_N)\in\R^{3N}$, with $a$ given by \eqref{eq:A-kernel} and the normalization of \eqref{eq:LN-def}.
The operator \eqref{eq:LM-generator} conserves total momentum and energy and is the natural $N$-particle counterpart of the Landau diffusion.
It underlies several hierarchy and compactness approaches:
Miot--Pulvirenti--Saffirio~\cite{MiotPulvirentiSaffirio2011} introduce a regularized Kac-type model and derive, along subsequences, a Coulomb Landau hierarchy;
Bobylev--Pulvirenti--Saffirio~\cite{BobylevPulvirentiSaffirio2013} obtain a first-order consistency result in a weak-coupling regime of Newtonian dynamics.
More recently, Carrillo--Feng--Guo--Jabin--Wang~\cite{CarrilloFengGuoJabinWang2024} develop relative-entropy methods for Maxwellian Landau equations in this diffusive setting;
Carrillo--Guo~\cite{CarrilloGuo2025_FI_Kac_to_Landau} prove Fisher-information decay for a Kac model and construct Landau--Coulomb hierarchies by compactness.
In much of this literature, the Landau master equation itself is obtained by \emph{formally} passing the Boltzmann master equation or a Kac jump process to the grazing limit.

$\bullet$ \emph{Interacting diffusions.}
A parallel influential line constructs the Landau equation as the nonlinear limit of McKean--Vlasov interacting diffusions.
A prototypical formal model is
\begin{equation}\label{eq:mf-SDE-general}
	\md V^{i,N}_t
	= \frac{1}{N}\sum_{j\ne i} b(V^{i,N}_t - V^{j,N}_t)\md t
	+ \frac{1}{\sqrt{N}}\sum_{j\ne i} \sigma(V^{i,N}_t - V^{j,N}_t)\md B^{ij}_t,
\end{equation}
where $(B^{ij})_{i\ne j}$ are independent Brownian motions, $\sigma\sigma^{\mathsf{T}} = a$, and the limiting law is expected to solve \eqref{eq:Landau}.
Such systems preserve mass, while their total momentum and energy are conserved in expectation under the required integrability assumptions.
In many classical constructions energy is conserved only in expectation, and the microscopic randomness is provided by external Brownian noise rather than by discrete conservative collisions; conservative variants use antisymmetric pair noises and an even square root $\sigma$ to recover exact conservation.
This framework supports several landmark results:
Fontbona--Gu\'erin--M\'el\'eard~\cite{FontbonaGuerinMeleardPTRF09} obtain convergence with explicit rates for Landau-type diffusive particle systems (Maxwellian-type interactions);
Carrapatoso~\cite{Carrapatoso2016KRM} (via a grazing limit of Kac's walk) proves quantitative propagation of chaos and entropic chaos, uniform in time, for Maxwell molecules, together with $N$-uniform convergence to equilibrium;
Fournier--Hauray~\cite{FournierHauray2016} establish propagation of chaos for moderately soft potentials $\gamma\in(-2,0)$ with explicit Wasserstein rates on a subrange;
and Fournier--Guillin~\cite{FournierGuillin2017} derive quantitative convergence from a conservative stochastic particle system to the Landau equation for hard potentials and Maxwell molecules.

Two important works settle propagation of chaos for very soft and Coulomb potentials within this diffusive paradigm.
Feng--Wang~\cite{FengWang2025_KacLandau} start from a conservative $N$-particle system obtained (formally) by passing Kac's walk to the grazing limit, whose generator has the Landau structure \eqref{eq:LM-generator}, and prove propagation of chaos for the full range $\gamma\in[-3,1]$ in weak, Wasserstein, entropic, and strong $L^1$ senses, relying on a duality and cluster-expansion framework (cf.~\cite{BreschDuerinckxJabin2024}) and new weighted and second-order Fisher-information estimates.
Their assumptions include initial boundedness, weighted Fisher information, and high moments.
Independently, Tabary~\cite{Tabary2025_PoC_Landau_Coulomb} considers a suitably regularized Landau-type interacting diffusion and establishes propagation of chaos toward the Landau equation for very soft and Coulomb potentials by a tightness-uniqueness argument on path space, crucially using Fisher-information dissipation and infinite-dimensional affinity.
He treats $\gamma\in[-3,-2]$ under finite energy, entropy, Fisher information, and a moment of order $m>3(2-\gamma)$.
These results provide conservative diffusive realizations of the Landau dynamics in the most singular regime.

\vspace{0.8ex}
\textbf{Our perspective: Kac's binary-collision derivation.}
The present work is of a different, complementary nature.
Building on the MBC dynamics introduced above, we adhere strictly to Kac's binary-collision philosophy: particles interact only through discrete two-body rotations on $\mathbb{S}^2$, with exact conservation of momentum and energy and invariance of the microcanonical measure.
The crucial point is that the generator is a genuine Kac-type jump operator, with no \emph{a priori} Landau diffusion term and outside the diffusive McKean--Vlasov class \eqref{eq:mf-SDE-general}.

In this way we obtain an \emph{end-to-end} derivation of the Landau equation from a conservative binary-collision model, in the precise sense that the only microscopic input is the MBC jump dynamics and the only macroscopic output is the Landau flow.

Moreover, our model connects back to the diffusive approaches.
In Section~\ref{sec:grazing-discussion}, we show that for each fixed $N$, the time-marginal laws of the MBC process converge, as $h\to0$, to those of the $N$-particle Landau diffusion generated by~\eqref{eq:LN-def}.
Thus the MBC dynamics provides a purely jump-based microscopic realization whose grazing limit recovers the Landau master equation used in many existing works.

The Fisher-dissipation geometry used here comes from Guill\'en--Silvestre \cite{GuillenSilvestre2025}.
We use the convex directional second-order functionals introduced by Tabary~\cite{Tabary2025_PoC_Landau_Coulomb} and his framework for passing information bounds to random limits, especially for very soft potentials.
The variational and infinite-dimensional arguments also draw on \cite{HaurayMischler2014,Rougerie2021Fisher}.
We give an independent martingale proof of exact affinity for directional second-order information, sharpening the two-sided comparison in \cite[Theorem~1.6(3)]{Tabary2025_PoC_Landau_Coulomb} to an equality.
Desvillettes \cite{Desvillettes2015EntropyLandau} proves weighted coercivity of entropy dissipation.
Ji \cite[Theorems~1.2 and~1.5]{Ji2024FisherDissipationLandau} establishes a second-order coercivity estimate for Fisher dissipation and constructs global solutions for very soft potentials under finite-entropy and moment assumptions.
A key coercivity step is the explicit reconstruction of the full log-gradient from its tangential differences by the moment identity \eqref{stat:eq:moment-reconstruction}.
Adapting the integration-by-parts argument of \cite[Section~2.3]{Desvillettes2015EntropyLandau}, we use this identity to obtain weighted fourth-power control of the full log-gradient under a positive lower bound on the covariance matrix (Proposition~\ref{vs:prop:static}).
The Fisher-information bound gives uniformly positive covariance for a comparison solution, and we use the weak evolution of the covariance to propagate this condition.
The continuation intervals are chosen using the time integral of the Fisher information.

\vspace{0.8ex}
We remark that the MBC process continues our earlier collision-based designs with synchronous, equal-time grouped pair updates \cite{DuLi2024arXiv,DuLiXY2025JCP}. 
In \cite{DuLi2024arXiv} we introduced a grazing-collision scheme formally derived from the Boltzmann operator and established quantitative propagation of chaos for a regularized Landau-type kernel with smooth bounded coefficients and uniformly elliptic diffusion, under suitable initial assumptions.
In \cite{DuLiXY2025JCP} we developed a structure-preserving collisional particle method for the full Landau operator, including Coulomb.
The present work replaces synchronous grouping by Kac-style random pair updates on the microcanonical shell, retaining the same conservative geometric structure while aligning with genuine binary-collision physics and, crucially, enabling a fully rigorous implementation of Kac's program for the Landau equation.

\subsection{Proof strategy}\label{sec:proof-strategy}

The proof follows a compactness--uniqueness argument: subsequential empirical limits solve the Landau equation, and uniqueness gives convergence of the full sequence to a deterministic solution.
The main difficulties are to retain information about each possible limiting density and to recover enough regularity from information restricted to collision directions.

The first step is to exploit the specific design of the MBC update.
At each jump, a uniformly chosen pair $\{i,j\}$ undergoes a rigid rotation of its
relative velocity on $\mathbb{S}^2$, driven by spherical Brownian motion over a
time $\tau_{ij}=h|v_i-v_j|^\gamma$.
At the level of densities, this two-body move coincides with the ``lifted''
Landau diffusion on the relative direction introduced in
Guill\'en--Silvestre~\cite{GuillenSilvestre2025} to prove Fisher-information
monotonicity for the Landau equation.
As a consequence, the MBC generator is compatible with the same functional structure: in Section~\ref{sec:discrete-information}, we show that if the initial law $F_0^N$ has finite normalized Fisher information, then the Fisher information along the MBC flow is non-increasing and satisfies quantitative dissipation estimates (Proposition~\ref{disc:thm:finite-N}).
The resulting near-diagonal bound on the two-particle marginals controls the singularity of the collision kernel at contact and is a key input for the subsequent compactness and identification estimates.
Together with moment bounds controlling velocity tails, it gives tightness and passage to the weak equation.
The fluctuation and grazing errors vanish independently, so no coupling of $N$ and $h_N$ is needed (Lemmas~\ref{lem:MD} and~\ref{lem:remainder}).

A random empirical limit produces a mixture of product measures, whose mean information may be hidden by averaging.
We prove exact affinity by showing that the $L^2$ representatives of directional distributional derivatives form a martingale across marginal orders (Lemma~\ref{disc:thm:score-affinity}).
We combine this identity with the classical affinity of Fisher information to obtain time-integrable $I(\rho_t)+K_\gamma(\rho_t^{\otimes2})$ for almost every limiting path (Proposition~\ref{disc:thm:level3-diss}).
Here $K_\gamma$ measures second derivatives of the two-particle density along rotations of the relative velocity that preserve pair momentum and energy (see \eqref{disc:eq:Kgamma}).

For $\gamma\in[-2,1]$, existing uniqueness results apply to the available bounds.
Hardy's inequality supplies the singular convolution control for $\gamma\in[-2,0)$, finite energy suffices at $\gamma=0$, and fourth moments give uniqueness for $\gamma\in(0,1]$ (Proposition~\ref{prop:unique}).

In the very soft range $\gamma\in[-3,-2)$, the key is to recover full log-gradient control from information in collision directions.
Fisher information and Hardy's inequality first justify the unsymmetrized weak equation and propagate the required moments directly (Lemma~\ref{vs:lem:moment}).
Second-order information controls fourth-power tangential differences of the log-gradient.
A moment identity then explicitly reconstructs each component, turning this tangential control into full log-gradient control under a positive lower bound on the covariance matrix.
Interpolation with the propagated moment then gives $\|\rho\|_\infty\lesssim1+K_\gamma(\rho^{\otimes2})^{3/4}$; the sublinear exponent makes this bound integrable in time (Proposition~\ref{vs:prop:static}).

The remaining difficulty is to maintain the covariance lower bound using only time-integrated information.
A comparison solution $g$ has time-integrable $L^\infty$ regularity and uniformly bounded Fisher information; the latter gives uniformly positive covariance (Lemma~\ref{vs:prop:comparison}).
The weak equation controls covariance increments by $\int(1+I)$, so on short intervals the static estimate provides $L^1_tL^\infty_v$ regularity and known uniqueness results identify the limit with $g$.
Equality at each endpoint restores the covariance bound.
Since the total Fisher-information integral is finite, we can continue over every finite interval (Proposition~\ref{vs:thm:closure}, proved in Section~\ref{sec:very-soft-closure}).

Combining tightness and the conservation laws, we show in Section~\ref{sec:compactness} that any limit point of $(\mu_t^N)_{N\ge1}$ solves the Landau equation \eqref{eq:Landau} with the regularity and a priori bounds inherited from the information estimates.
Existing uniqueness results for \eqref{eq:Landau} \cite{FournierGuerin2009JFA,Fournier2010CMP,FournierHeydecker2021} then imply that the whole sequence converges to the same deterministic limit, which yields propagation of chaos as stated in Theorem~\ref{thm:main}.
Finally, energy conservation upgrades convergence to uniform-in-time $W_2$ convergence in expectation, while compactness under the Fisher bound gives strong $L^1$ convergence of fixed-order marginals at each fixed time.

\vspace{0.8ex}
\emph{Layout of the paper.}
Sections~\ref{sec:discrete-information}--\ref{sec:very-soft-closure} are devoted to the proof of Theorem~\ref{thm:main}.
Section~\ref{sec:discrete-information} introduces the MBC dynamics and establishes Fisher-information dissipation and the near-diagonal bounds on the $N$-particle laws.
It also passes the information estimates to empirical limits.
Section~\ref{sec:Tightness-of-empirical} uses these estimates to prove tightness of the empirical measures and continuity of subsequential limits.
In Section~\ref{sec:Identification-of-the}, we identify any limit point with a weak solution of the Landau equation \eqref{eq:Landau}.
Section~\ref{sec:uniqueness-completion} combines this with the available uniqueness theory to obtain convergence of the whole sequence and propagation of chaos for all $\gamma\in[-3,1]$.
Section~\ref{sec:very-soft-closure} establishes the regularity estimates for very soft potentials and completes the proof of Proposition~\ref{vs:thm:closure}.
Section~\ref{sec:grazing-discussion} studies the fixed-$N$ grazing-collision limit and shows that, as $h\to0$, the MBC laws converge to the $N$-particle Landau diffusion generated by \eqref{eq:LN-def}.
Appendix~\ref{sec:appendix-ergodicity} collects additional structural properties of the MBC process (including ergodicity), while Appendix~\ref{sec:auxiliary-information} contains auxiliary proofs for the information estimates.

\section{Information estimates for the MBC process}
\label{sec:discrete-information}\label{sec:discrete-to-level3}\label{sec:The-MBC-process}

The particle estimates come from the two-body flow, convexity, and conservation of energy.
The passage to empirical limits then uses the variational structure of information, independently of the particle generator.

\subsection{Pair-flow dissipation and particle estimates}

The natural state space for the MBC process is the microcanonical
shell
\begin{equation}\label{eq:M-shell}
\mathcal{M}_{E,P}^{N}:=\Big\{\bm{v}\in\R^{3N}:\frac{1}{N}\sum_{i}v_{i}=P,\ \frac{1}{N}\sum_{i}|v_{i}|^{2}=2E\Big\}
\end{equation}
with prescribed mean momentum $P\in\R^{3}$ and energy per particle $E>0$, where
$2E-|P|^{2}>0$. It is a compact smooth submanifold of $\R^{3N}$
of dimension $3N-4$; we equip it with the microcanonical measure
$\mu_{\mathrm{mc}}^N$, i.e., the normalized surface measure induced
by the ambient Euclidean structure. The MBC dynamics is a pure-jump
continuous-time Markov process on $\mathcal{M}_{E,P}^{N}$. 

For product initial data, the shell depends on the initial realization, and we use the whole-space formulation below.

To describe its generator, fix $\bm{v}\in\R^{3N}$ and
an unordered pair $\{i,j\}$. Set
\begin{equation}
\bar{v}_{ij}=\frac{v_{i}+v_{j}}{2},\quad r_{ij}=|v_{i}-v_{j}|,\quad e_{ij}=\widehat{v_{i}-v_{j}}\in\mathbb{S}^{2},\label{eq:CoV}
\end{equation}
when $r_{ij}>0$, and let $\bm{v}_{[ij]}$ collect all coordinates except $(v_{i},v_{j})$.
Using the change of variables
\[
\R^{3N}\ni\bm{v}\longleftrightarrow(\bm{v}_{[ij]},\bar{v}_{ij},r_{ij},e_{ij})\in\R^{3(N-2)}\times\R^{3}\times[0,\infty)\times\mathbb{S}^{2},
\]
we may write functions as $\Phi(\bm{v})=\Phi(\bm{v}_{[ij]},\bar{v}_{ij},r_{ij},e_{ij})$. 

We write $P_\tau^{\S^2}=\ee^{\tau\Delta_{\S^2}}$ for the spherical heat semigroup and $P_\tau^{\S^2}(e,\md e')$ for its transition kernel.
Recall the pair-dependent time scale $hr_{ij}^\gamma$, with the ``do nothing'' convention on $\{r_{ij}=0\}$.
The one-step operator for the pair is given by
\begin{equation}
 (\mathcal Q_{ij}^h\Phi)(\bm v)
 =\int_{\S^2}
 \Phi(\bm v_{[ij]},\bar v_{ij},r_{ij},e')
 P_{hr_{ij}^\gamma}^{\S^2}(e_{ij},\md e').
 \label{eq:Q_ij}
\end{equation}

At the density level, the pair transition can be written as
\begin{equation}
 \mathcal Q_{ij}^h=\ee^{hA_{ij}},\qquad
 A_{ij}=\frac14\big(\nabla_{v_i}-\nabla_{v_j}\big)\cdot
 \Big[a(v_i-v_j)\big(\nabla_{v_i}-\nabla_{v_j}\big)\Big].
 \label{eq:MBC-definition}
\end{equation}
In pair coordinates $z=v_i-v_j$ (and $x=v_i+v_j$) and $r=|z|$, we have $\nabla_{v_i}-\nabla_{v_j}=2\nabla_z$ and
\[
 \nabla_z\cdot
 \big[r^{\gamma+2}(I-\hat z\otimes\hat z)\nabla_z\big]
 =r^\gamma\Delta_{\S^2},
\]
so $A_{ij}=r^\gamma\Delta_{\S^2}$ on the relative direction.
Thus $\mathcal Q_{ij}^h=\ee^{hA_{ij}}$ acts by spherical heat flow for time $hr^\gamma$, with the identity update at $r=0$.

With these notations, the embedded one-collision Markov kernel and the generator are
\begin{equation}
 \mathcal Q^{N,h}
 =\frac{2}{N(N-1)}\sum_{i<j}\mathcal Q_{ij}^h,
 \qquad
 \mathcal L^{N,h}=\frac Nh(\mathcal Q^{N,h}-I).
 \label{disc:eq:MBC-generator}
\end{equation}
The same formulas define operators on every bounded Borel function on $\R^{3N}$; paths started on a fixed shell $\mathcal M_{E,P}^N$ remain there for all $t\ge0$.
The corresponding continuous-time semigroup is $P_t^{N,h}=\exp(t\mathcal L^{N,h})$.

\begin{remark}
\label{rem:microcanonical-structure}
Our MBC dynamics is a Kac-type conservative binary jump process on the microcanonical shell.
Invariance, detailed balance and self-adjointness follow directly from the reversibility of the spherical heat semigroup on $\S^2$ and rotational symmetry of the conditional relative direction.
This is fully consistent with the classical Kac framework \cite{MischlerMouhot2013Invent}, but we provide a short self-contained proof of these properties in Appendix~\ref{sec:appendix-ergodicity} tailored to the present kernel.
The same symmetry holds for the extensions to $L^2(\R^{3N})$ with respect to Lebesgue measure, as shown in Corollary~\ref{cor:whole-space}.
\end{remark}

We next record a mild regularity result for the whole-space dynamics: absolute continuity of the initial law is preserved by the MBC process for all positive times.
For any initial probability density, Lebesgue reversibility identifies the adjoint on densities as the same operator, and iteration gives the series representation
\[
 F_t^N=(P_t^{N,h})^*F_0^N=P_t^{N,h}F_0^N
 =\ee^{-Nt/h}\sum_{n=0}^\infty
 \frac{(Nt/h)^n}{n!}(\mathcal Q^{N,h})^nF_0^N.
\]
Each term is a probability density, and the series converges in $L^1$.
Thus, absolute continuity is preserved and the mass is conserved; see also Proposition~\ref{prop:a.c.bound} for the $L^\infty$ contraction.

For a probability measure $P$ on $\R^d$, the (extended) Fisher information is
\begin{equation}
 I(P):=
 \begin{cases}
  4\displaystyle\int_{\R^d}|\nabla\sqrt F|^2\md x,
   &\text{if }P=F\md x\text{ and }\sqrt F\in H^1(\R^d),\\
  +\infty,&\text{otherwise}.
 \end{cases}
 \label{disc:eq:fisher-information}
\end{equation}
The control of the Fisher information (FI) is decisive in our analysis, especially for the very-soft regime $\gamma\in[-3,-2)$, where compactness and identification hinge on uniform FI bounds.
Guill\'en--Silvestre \cite{GuillenSilvestre2025} observed that
the time derivative of the Landau FI can be controlled
via that of a lifted linear diffusion equation. By exploiting the
geometric structure of the lifted Landau operator, they developed
new analytical tools, established monotonicity of FI along the lifted
flow, and then transferred this decay back to the Landau equation.

In our setting, each microcanonical binary update acts on a pair $\{i,j\}$
exactly via that lifted two-body diffusion, up to the constant time normalization in \eqref{eq:MBC-definition}.
Hence the law of the MBC process inherits the same FI monotonicity from the lifted flow; this provides the analytic backbone for the very-soft potentials, where neither entropy nor moments alone suffice.
We also retain the second-order dissipation in the convex directional form developed by Tabary \cite{Tabary2025_PoC_Landau_Coulomb}.
Its contraction under marginalization over the other particle coordinates will allow us to pass the dissipation to empirical limits.
Let $e_k$ be the $k$-th standard basis vector in $\R^3$, set
\begin{equation}
 b_k(z):=e_k\times z,
 \label{disc:eq:rotational-fields}
\end{equation}
and introduce the derivations
\begin{equation}
\partial_{ij,k}
 :=r_{ij}^{\gamma/4-1/2}\,
 b_k(v_i-v_j)\cdot\big(\nabla_{v_i}-\nabla_{v_j}\big).\label{disc:eq:derivation}
\end{equation}
The radial coefficient is constant along these rotation directions, and the vector field in \eqref{disc:eq:derivation} is divergence free.
For a smooth positive function $G$ on $\R^{3N}$, define
\begin{equation}
  K_{ij,\gamma}(G)
 :=
 \sum_{k=1}^3
 \int_{\R^{3N}}
 \frac{|\partial_{ij,k}^2G|^2}{G}
 =
 \sum_{k=1}^3
 \int_{\R^{3N}}
 \frac{r_{ij}^{\gamma-2}
 \big|\big(b_k(v_i-v_j)\cdot(\nabla_{v_i}-\nabla_{v_j})\big)^2G\big|^2}{G}.
\label{disc:eq:Kpair}
\end{equation}
We call $K_{ij,\gamma}$ the second-order Fisher information along the collision directions of the pair $(i,j)$.
The lower semicontinuous extension to nonsmooth measures, used in the integrated estimate below, is constructed in Remark~\ref{disc:rem:singular-affinity}.

\begin{lemma}\label{disc:prop:pair-FI}
For every $\gamma\in[-3,1]$, there is $c_\gamma>0$ with the following property.
Let $G$ be a smooth positive probability density on $\R^{3N}$, symmetric in $v_i,v_j$, with sufficient decay to justify the integrals and integrations by parts below.
Then the lifted pair flow obeys
\begin{equation}
\frac{\md}{\md s}I(\ee^{sA_{ij}}G)
 \le - c_\gamma K_{ij,\gamma}(\ee^{sA_{ij}}G).\label{disc:eq:FI-pair-coercivity}
\end{equation}
For every pair-symmetric probability density $G$ on $\R^{3N}$ with $I(G)<\infty$, the integrated form is
\[
 I(\ee^{hA_{ij}}G)+c_\gamma\int_0^h
 K_{ij,\gamma}(\ee^{sA_{ij}}G)\md s\le I(G)
 \qquad(h>0).
\]
\end{lemma}

\begin{proof}
For a fixed unordered pair $\{i,j\}$, denote by $y=\bm v_{[ij]}\in\R^{3(N-2)}$ all the remaining coordinates and set $G_y(v,w):=G(v,w;y)$.
The Fisher information of $G$ decomposes as
\[
I(G)=\int_{\R^{3(N-2)}}\bigg(I_{ij}(G_{y})+\sum_{k\notin\{i,j\}}\int_{\R^{6}}G_{y}|\nabla_{v_{k}}\log G_{y}|^{2}\md v\md w\bigg)\md y,
\]
where $I_{ij}$ denotes the two-body Fisher information.
Applying $\mathcal Q_{ij}^s=\ee^{sA_{ij}}$ affects only $(v_i,v_j)$.
For any derivative $D$ in a remaining coordinate, the pair kernel is independent of that coordinate, so $D(\ee^{sA_{ij}}G)=\ee^{sA_{ij}}(DG)$.
The Cauchy--Schwarz inequality for this Markov kernel gives
\[
 \big|\ee^{sA_{ij}}(DG)\big|^2
 =\bigg|\ee^{sA_{ij}}\Big(\frac{DG}{\sqrt G}\,\sqrt G\Big)\bigg|^2
 \le\ee^{sA_{ij}}\Big(\frac{|DG|^2}{G}\Big)
       \ee^{sA_{ij}}G.
\]
Dividing by the positive quantity $\ee^{sA_{ij}}G$ yields
\[
 \frac{|D(\ee^{sA_{ij}}G)|^2}{\ee^{sA_{ij}}G}
 \le \ee^{sA_{ij}}\Big(\frac{|DG|^2}{G}\Big).
\]
Integration and preservation of Lebesgue measure show that all non-collisional gradient components are non-increasing.

For the derivatives in $v_i,v_j$, apply Guill\'en--Silvestre's lifted identity and spherical estimate \cite[(8.3), Lemmas~8.4 and~10.1]{GuillenSilvestre2025}, with $\alpha(r)=r^\gamma$, keeping the remaining coordinates fixed.
Pair symmetry permits their spherical estimate.
Integrating over those coordinates and using their Fisher-information contraction proved above gives
\begin{equation}
\begin{aligned}
 \frac{\md}{\md s}I(\ee^{sA_{ij}}G)
 & \le - \frac14\Big(1-\frac{\gamma^2}{19}\Big)
 \sum_{k=1}^3\int_{\R^{3N}}(\ee^{sA_{ij}}G)
 \big|\partial_{ij,k}^{\,2}\log(\ee^{sA_{ij}}G)\big|^2\\
 & \le - \frac{1}{16} \Big(1-\frac{\gamma^2}{19}\Big)
 \sum_{k=1}^3\int_{\R^{3N}}
 \frac{|\partial_{ij,k}^{\,2}(\ee^{sA_{ij}}G)|^2}{\ee^{sA_{ij}}G},
\end{aligned}\label{disc:eq:GS-coercive-chain}
\end{equation}
where the last inequality uses \cite[Proposition~2.1, with $\beta=0$]{Tabary2025_PoC_Landau_Coulomb}.
Therefore, \eqref{disc:eq:FI-pair-coercivity} follows with
\[
 c_\gamma=\frac1{16}\Big(1-\frac{\gamma^2}{19}\Big)
 \ge\frac5{152}\qquad(\gamma\in[-3,1]).
\]

For a pair-symmetric density $G$ with $I(G)<\infty$, standard cutoff and smoothing yield smooth positive pair-symmetric densities $G_n$,
with Gaussian tails and an isotropic Gaussian profile near the pair diagonal, such that $G_n\to G$ in $L^1$ and $I(G_n)\to I(G)$.
Apply the integrated estimate to $G_n$ and let $n\to\infty$.
By the $L^1$ contraction of the pair semigroup, lower semicontinuity, and Fatou's lemma, we obtain
\[
 I(\ee^{hA_{ij}}G)
 +c_\gamma\int_0^h K_{ij,\gamma}(\ee^{sA_{ij}}G)\md s
 \le I(G).
\]
This concludes the proof.
\end{proof}

For $h>0$, define the information averaged along the pair flow by
\begin{equation}
\overline{ K}_{ij,\gamma}^{\,h}(F)
 :=\frac1h\int_0^h
  K_{ij,\gamma}(\ee^{sA_{ij}}F)\md s.
 \label{disc:eq:averaged-information}
\end{equation}

\begin{proposition}
\label{disc:thm:finite-N}
Let $(F_t^N)_{t\ge0}$ be the law of the MBC process started from exchangeable $F_0^N$ with $\sqrt{F_0^N}\in H^1(\R^{3N})$.
Then the Fisher information $I(F_t^N)$ is finite for each $t\ge0$ and non-increasing with respect to $t$.
For every $T>0$, its Fisher information and the averaged information \eqref{disc:eq:averaged-information} satisfy
\begin{equation}
\frac1N I(F_T^N)
 +c_\gamma\int_0^T
 \overline{ K}_{12,\gamma}^{\,h}(F_t^N)\md t
 \le\frac1N I(F_0^N).\label{disc:eq:FI-jump-diss}
\end{equation}
Consequently, if $F_0^N=f_0^{\otimes N}$, then along any sequence $h_N>0$
\begin{equation}
 \frac1N I(F_0^N)=I(f_0),\qquad
 \sup_N\int_0^T
 \overline{ K}_{12,\gamma}^{\,h_N}(F_t^N)\md t
 \le c_\gamma^{-1}I(f_0).
\label{disc:eq:uniform-K}
\end{equation}
\end{proposition}

\begin{proof}
Let $G_n=(\mathcal Q^{N,h})^nF_0^N$.
Apply the integrated estimate of Lemma~\ref{disc:prop:pair-FI} to each pair.
By convexity of Fisher information and exchangeability, we obtain
\[
 I(G_n)-I(G_{n+1})
 \ge c_\gamma h\,\overline K_{12,\gamma}^{\,h}(G_n).
\]
In particular, $I(G_n)\le I(F_0^N)$.
With $p_n(t)=\ee^{-Nt/h}(Nt/h)^n/n!$, the Poisson formula is $F_t^N=\sum_{n\ge0}p_n(t)G_n$.
Therefore
\[
 \begin{split}
 -\frac{\md}{\md t}\sum_{n\ge0}p_n(t)I(G_n)
 &=\frac Nh\sum_{n\ge0}p_n(t)
       \big[I(G_n)-I(G_{n+1})\big]\\
 &\ge N c_\gamma\sum_{n\ge0}p_n(t)
       \overline K_{12,\gamma}^{\,h}(G_n)\\
 &\ge N c_\gamma
       \overline K_{12,\gamma}^{\,h}(F_t^N).
 \end{split}
\]
The series may be differentiated because its coefficients $I(G_n)$ are uniformly bounded; the last inequality uses convexity of the averaged functional.
Integrating and using $I(F_T^N)\le\sum_np_n(T)I(G_n)$ proves \eqref{disc:eq:FI-jump-diss}.
Since the same argument applies starting at any time, the Fisher information is non-increasing.
In Theorem~\ref{thm:main} we start the MBC process from the product law $f_0^{\otimes N}$ with $I(f_0)<\infty$.
A direct computation gives $I(F_0^N)=N I(f_0)$, so the assumption is satisfied by the concrete initial data used in Theorem~\ref{thm:main}, and \eqref{disc:eq:uniform-K} follows.
\end{proof}

The dissipation of Fisher information established above is a structural input to our argument.
In particular, the marginal bound is the mechanism by which Fisher dissipation yields quantitative integrability near the collision set $\{v_i=v_j\}$.
Indeed, the conditional expectation formula for the log-gradient and Jensen's inequality give
\[
 I(F_t^{k|N})\le\frac{k}{N}I(F_t^N)
 \le\frac{k}{N}I(F_0^N),\qquad 1\le k\le N,
\]
where the factor $k/N$ follows from exchangeability.
Combined with Hardy's inequality in the relative variable, this produces the following near-diagonal estimate.

For every probability density $\rho$ with finite Fisher information, applying Hardy's inequality to $\sqrt \rho(\cdot+x)$ gives
\begin{equation}
 \sup_x\int(1+|x-y|^{-2})\rho(y)\md y\le1+I(\rho).
 \label{eq:fisher-hardy}
\end{equation}
For product initial data, we combine Proposition~\ref{disc:thm:finite-N}, marginal contraction, and \eqref{eq:fisher-hardy} to obtain (see \cite[Lemma~5.1]{Tabary2025_PoC_Landau_Coulomb})
\begin{equation}
\sup_{N,t}\E|v_1^N(t)-v_2^N(t)|^{-2}
 \le C I(f_0).\label{eq:NDC}
\end{equation}
Let $f_{t}^{Z}$ be the Lebesgue density of $Z_{t}:=(v_{1}^{N}(t)-v_{2}^{N}(t))/\sqrt{2}$.
By the classical Hardy inequality, applied to $\sqrt{f_{t}^{Z}}$, we get
\[
\E\big[|v_{1}^{N}(t)-v_{2}^{N}(t)|^{-2}\big]=\frac{1}{2}\int_{\R^{3}}\frac{f_{t}^{Z}(z)}{|z|^{2}}\md z\le2\int_{\R^{3}}\Big|\nabla_{z}\sqrt{f_{t}^{Z}(z)}\Big|^{2}\md z=\frac{1}{2}I(f_{t}^{Z}).
\]
With change of variables $w=(v_{1}+v_{2})/\sqrt{2},\ z=(v_{1}-v_{2})/\sqrt{2}$,
we have $f_{t}^{Z}(z)=\int F_{t}^{2|N}(w,z)\md w$.
The same conditional expectation formula and Jensen's inequality yield
\begin{align*}
I(f_{t}^{Z}) & =\int_{\R^{3}}|\nabla_{z}\log f_{t}^{Z}(z)|^{2}f_{t}^{Z}(z)\md z\\
 & \le\int_{\R^{6}}|\nabla_{z}\log F_{t}^{2|N}(w,z)|^{2}F_{t}^{2|N}(w,z)\md w\md z\le I(F_{t}^{2|N}).
\end{align*}
The change of variables is orthogonal, so it preserves Fisher information.
The preceding marginal bound and Proposition~\ref{disc:thm:finite-N} therefore prove \eqref{eq:NDC}.

For the energy bounds and localization, define
\begin{equation}
 E_0:=\int|v|^2f_0(v)\md v,\quad
 M_2^N(t):=\big\langle\mu_t^N,|v|^2\big\rangle,\quad
 \Omega_R^N :=\{M_2^N(0)\le R\}.
 \label{eq:energy-localization}
\end{equation}

\begin{lemma}
\label{prop:M4}
Let $\gamma\in[0,1]$ and assume $f_0\in\Pp_4(\R^3)$.
Then, for every $0<R,T<\infty$,
\[
\begin{gathered}
 \sup_N\sup_{0<h\le1}\sup_{t\le T}
 \E\bigg[
 \frac1N\sum_{i=1}^N\langle v_i^N(t)\rangle^4
 \mathbf1_{\Omega_R^N}\bigg]
 \le C_{T,R,\gamma,f_0},\\
 \sup_N\mathbb P((\Omega_R^N)^c)\le\frac{E_0}{R}.
\end{gathered}
\]
where $\langle v\rangle:=\sqrt{1+|v|^{2}}$.
\end{lemma}

\begin{proof}
It is enough to work with
\[
 M_4^N(\boldsymbol v)=\frac1N\sum_{i=1}^N|v_i|^4,
 \qquad
 M_2^N(\boldsymbol v)=\frac1N\sum_{i=1}^N|v_i|^2,
\]
because the bracketed fourth moment differs by lower moments.
Fix a pair; the case $v_i=v_j$ contributes zero because the update is the identity.
Otherwise, write
\[
 u=\frac{v_i+v_j}{2},\qquad v_i-v_j=re.
\]
The sum of its fourth moments is
\[
 \Big|u+\frac r2e\Big|^4
 +\Big|u-\frac r2e\Big|^4
 =
 2\Big(|u|^2+\frac{r^2}{4}\Big)^2
 +2r^2(u\cdot e)^2.
\]
The pair update preserves $u$ and $r$, hence
\[
 (\mathcal Q_{ij}^h-I)\Big(|u|^2+\frac{r^2}{4}\Big)^2=0.
\]
Only the angular term $2r^2(u\cdot e)^2$ contributes to the change in the pair fourth moment.
The degree-two spherical harmonic identity is exact:
\[
P_\tau^{\S^2}(u\cdot e)^2
 =\ee^{-6\tau}(u\cdot e)^2
 +(1-\ee^{-6\tau})\frac{|u|^2}{3}.
\]
Consequently, the conditional expected change of the pair fourth moment in one update is
\[
2r^2(1-\ee^{-6hr^\gamma})
 \Big(\frac{|u|^2}{3}-(u\cdot e)^2\Big).
\]

A direct calculation gives
\[
r^2\Big(\frac{|u|^2}{3}-(u\cdot e)^2\Big)
 =-\frac16\big(|v_i|^4+|v_j|^4\big)
 +\frac23|v_i|^2|v_j|^2-\frac13(v_i\cdot v_j)^2\\
 \le\frac23|v_i|^2|v_j|^2.
\]
Since a pair update leaves all other coordinates unchanged, the MBC generator \eqref{disc:eq:MBC-generator} acts on $M_4^N$ as
\[
 \mathcal L^{N,h}M_4^N
 =\frac{2}{hN(N-1)}\sum_{i<j}
 (\mathcal Q_{ij}^h-I)\big(|v_i|^4+|v_j|^4\big).
\]
Using the preceding pair estimate and $1-\ee^{-6hr^\gamma}\le6hr^\gamma$, we obtain the pointwise upper bound
\[
 \mathcal{L}^{N,h}M_4^N
 \le C_\gamma
 \frac{2}{N(N-1)}\sum_{i<j}
 |v_i-v_j|^\gamma |v_i|^2|v_j|^2.
 \]
Using $|v_i-v_j|^\gamma\le C_\gamma (|v_i|^\gamma+|v_j|^\gamma)$, the right-hand side is bounded, uniformly in $N,h$, by
\[
 C_\gamma M_2^N\Big(\frac1N\sum_i|v_i|^{\gamma+2}\Big)
 \le C_\gamma M_2^N(1+M_4^N),
\]
where the inequality uses $|v_i|^{\gamma+2}\le1+|v_i|^4$.
On the conserved event $\Omega_R^N$, the second moment is at most $R$, so
\[
 \mathcal L^{N,h}(M_4^N\mathbf1_{\Omega_R^N})
 \le C_\gamma R(1+M_4^N)\mathbf1_{\Omega_R^N}.
\]
For fixed $N$, the function on the left is bounded because $M_4^N\le N(M_2^N)^2\le NR^2$ on $\Omega_R^N$.
Dynkin's formula and Gr\"onwall's inequality thus give directly
\[
 \sup_{t\le T}\E\big[M_4^N(t)\mathbf1_{\Omega_R^N}\big]
 \le C_{T,R,\gamma}\big[1+\E M_4^N(0)\big].
\]
The initial fourth moment is finite.
Finally, $\mathbb P((\Omega_R^N)^c)\le E_0/R$ follows from Markov's inequality and conservation of the expected energy.
\end{proof}

\subsection{Information in empirical limits}

The particle information estimates must pass to each random limiting path.
The following proposition provides the bounds used in the soft-potential uniqueness argument.

Fix $T>0$.
The space $D([0,T],\Pp(\R^3))$ is equipped with the $J_1$ topology induced by weak convergence in $\Pp(\R^3)$, and its continuous-path subspace is denoted by $C_w([0,T],\Pp(\R^3))$.
For a path law $\pi$, write $\rho=(\rho_t)$ for the canonical random path and $\E_\pi$ for expectation under $\pi$.

For an exchangeable $N$-particle density $G$, relabeling particles shows that $K_{ij,\gamma}(G)$ in \eqref{disc:eq:Kpair} is independent of the chosen pair $(i,j)$.
We write $K_\gamma$ for the extension of its two-particle version $K_{12,\gamma}$ to probability measures on $\R^6$; its precise construction is given in Remark~\ref{disc:rem:singular-affinity} below.

\begin{proposition}
\label{disc:thm:level3-diss}
Let $\gamma\in[-3,1]$, let $f_0\in\Pp_2(\R^3)$ be a density with $I(f_0)<\infty$, and let $h_N>0$ tend to zero.
Consider the MBC processes with collision steps $h_N$ and initial laws $F_0^N=f_0^{\otimes N}$.
Assume that the laws of their empirical paths $(\mu_t^N)_{t\in[0,T]}$ are $C$-tight in $D([0,T],\Pp(\R^3))$.

Then, every subsequential limit law $\pi$ is supported on $C_w([0,T],\Pp(\R^3))$.
For each $t\in[0,T]$, the measure $\rho_t$ is absolutely continuous with respect to Lebesgue measure for $\pi$-almost every path, and
\begin{equation}
 \E_\pi I(\rho_t)\le I(f_0),\qquad
 \E_\pi\int_0^T I(\rho_t)\md t\le T I(f_0).
 \label{eq:level3-I}
\end{equation}
It also satisfies
\begin{equation}
\E_\pi\int_0^T
  K_\gamma(\rho_t^{\otimes2})\md t
 \le C_\gamma I(f_0).\label{disc:eq:level3-K}
\end{equation}
In particular, for $\pi$-almost every path,
\begin{equation}
\int_0^T\big[I(\rho_t)+K_\gamma(\rho_t^{\otimes2})\big]\md t<\infty.
\label{disc:eq:pathwise-K}
\end{equation}
\end{proposition}

The proof, given at the end of this subsection, combines lower semicontinuity and affinity: the former passes the particle estimates to mixed marginals, and the latter recovers the mean information of their components.
Lemma~\ref{disc:lem:invisibility} ensures that the short pair flow appearing in the particle estimate has the same marginal limit.

For $\gamma\in[-2,0)$, the Fisher bound suffices for Proposition~\ref{prop:unique};
for $\gamma\in[-3,-2)$, both bounds are used in Proposition~\ref{vs:thm:closure} to obtain the time-integrable $L^\infty$ regularity needed for uniqueness.

For a smooth positive density $p$ and $m=1,2$, the directional information has the form
\[
\mathcal{K}_{X,m}(p) := \int \frac{|X^m p|^2}{p},
\]
where $X=b\cdot\nabla$ and $b$ is a smooth bounded divergence-free vector field.
The variational formula below extends this quantity to arbitrary probability measures, so no density or derivative is assumed before the information bounds are proved;
compare \cite[Lemmas~3.5 and~3.7]{HaurayMischler2014} and \cite[Section~2 and Theorem~1.6(1)--(2)]{Tabary2025_PoC_Landau_Coulomb}.

\begin{lemma}\label{disc:prop:active-calculus}
Let $k_0\ge1$, $m\in\{1,2\}$, and $X=b\cdot\nabla_{1:k_0}$, where $b=b(x_1,\ldots,x_{k_0})\in\R^{3k_0}$ is smooth, bounded and divergence free.
For $j\ge k_0$ and $P\in\Pp(\R^{3j})$, the functional
\begin{equation}
\mathcal{K}_{X,m}(P)
 :=
 \sup_{\psi\in C_c^\infty(\R^{3j})}
 \Big\{
  2(-1)^m\int X^m\psi\md P-\int\psi^2\md P
 \Big\}
 \label{disc:eq:FXm}
\end{equation}
is convex and weakly lower semicontinuous.
It also contracts under marginalization: if $Q$ is obtained from $P$ by integrating out any collection of coordinates beyond the first $k_0$, then
$\mathcal{K}_{X,m}(Q)\le\mathcal{K}_{X,m}(P)$.

Moreover, it is finite if and only if there exists $S\in L^2(P)$, unique up to $P$-null sets, such that
\begin{equation}
X^mP=SP\quad\text{in }\mathcal D',
 \qquad
 \mathcal{K}_{X,m}(P)
 =\int|S|^2\md P.\label{disc:eq:Riesz-score}
\end{equation}
Here $X^mP=SP$ means that $(-1)^m\int X^m\psi\md P=\int S\psi\md P$ for every $\psi\in C_c^\infty(\R^{3j})$.

For a smooth positive density $p$, $\mathcal{K}_{X,m}(p)=\int|X^mp|^2/p$; when this quantity is finite, its $L^2(p)$ representative is $S=X^mp/p$.
\end{lemma}

The proof is given in Appendix~\ref{sec:auxiliary-information}: cutoff tests give marginal contraction, and the Riesz representation theorem gives the distributional representation.

Affinity of the usual level-3 Fisher information is established in \cite[Theorem~5.7]{HaurayMischler2014}; see also \cite[Theorem~2.4]{Rougerie2021Fisher}.
For directional functionals, Tabary's preprint \cite[Theorem~1.6(3)]{Tabary2025_PoC_Landau_Coulomb} proves exact first-order affinity and a two-sided second-order comparison with factor $16$, under finite mean Fisher information and second moment.
Tabary's comparison already suffices for Proposition~\ref{disc:thm:level3-diss} when applied to the bounded annular fields below.
The limiting laws satisfy its hypotheses by the Fisher bound and conserved expected particle energy, and the factor $16$ only changes the constant.
The following lemma gives an independent martingale proof of exact affinity, replacing the second-order comparison by equality.

\begin{lemma}
\label{disc:thm:score-affinity}
Let $\Pi\in\Pp(\Pp(\R^3))$, and write its de Finetti marginals as
\[
 \Pi^j=\int_{\Pp(\R^3)}\rho^{\otimes j}\,\Pi(\md\rho).
\]
For $X$ and $m$ as in Lemma~\ref{disc:prop:active-calculus}, and for $j\ge k_0$,
\begin{equation}
\lim_{j\to\infty}\mathcal{K}_{X,m}(\Pi^j)
 =
 \int_{\Pp(\R^3)}
 \mathcal{K}_{X,m}(\rho^{\otimes k_0})\,\Pi(\md\rho).\label{disc:eq:score-affinity}
\end{equation}
The sequence on the left is increasing in $j$.
\end{lemma}

The proof, given in Appendix~\ref{sec:auxiliary-information}, constructs a martingale of the $L^2$ representatives and recovers the common random measure from the conditional law of large numbers.

\begin{remark}[Singular Landau directions]\label{disc:rem:singular-affinity}
Recall from \eqref{disc:eq:rotational-fields} that $b_k(z)=e_k\times z$.
The Landau directions \eqref{disc:eq:derivation} need not be smooth and bounded, as required in Lemmas~\ref{disc:prop:active-calculus} and~\ref{disc:thm:score-affinity}.
Radial cutoffs extend these results to the collision directions and yield lower semicontinuity, marginal contraction, and exact affinity for $K_{12,\gamma}$, which are used in Proposition~\ref{disc:thm:level3-diss}.
Put $r=r_{12}=|v_1-v_2|$ and
\[
 X_k:= r^{\gamma/4-1/2} b_k(v_1-v_2)\cdot(\nabla_{v_1}-\nabla_{v_2}),\quad k=1,2,3.
\]
The key observation is that $X_kr=0$, so radial cutoffs introduce no derivative terms.

For $X=X_k$, choose nested smooth radial cutoffs such that
\[
 0\le\chi_n\le\chi_{n+1}\le1,\qquad
 \chi_n(r)\uparrow1\quad(r>0),\qquad
 \chi_{n+1}=1\ \text{on }\operatorname{supp}\chi_n,
\]
with each $\chi_n$ vanishing near zero and near infinity.
Then $X_n=\chi_n(r)X$ is smooth, bounded and divergence free, and
\[
 X_n^2=\chi_n^2X^2.
\]
Define the extension away from the diagonal by
\begin{equation}
\mathcal{K}_{X,2}^{\rm off}(P)
 :=\sup_n\mathcal{K}_{X_n,2}(P).\label{disc:eq:off-relaxation}
\end{equation}
Suppose $\mathcal{K}_{X_{n+1},2}(P)<\infty$.
Multiply the representing function for $X_{n+1}^2P$ by $\chi_n^2$.
By the nested-support condition and Lemma~\ref{disc:prop:active-calculus}, this product represents $X_n^2P$.
Since $0\le\chi_n\le1$, this proves $\mathcal{K}_{X_n,2}(P)\le\mathcal{K}_{X_{n+1},2}(P)$; the inequality is immediate when the latter is infinite.
Thus the terms in \eqref{disc:eq:off-relaxation} are increasing for every probability measure.
For a general probability measure $P \in\Pp(\R^6)$, define
\begin{equation}
 K_\gamma(P):=K_{12,\gamma}(P)
 =\sum_{k=1}^3\mathcal{K}_{X_k,2}^{\rm off}(P).
 \label{disc:eq:Kgamma}
\end{equation}
For $P\in\Pp(\R^{3j})$ with $j\ge2$, the same sum defines $K_{12,\gamma}(P)$, with $X_k$ acting only on the first two coordinates.
For a smooth density off the diagonal this definition reduces to \eqref{disc:eq:Kpair}.
As a finite sum of suprema of the functionals in Lemma~\ref{disc:prop:active-calculus}, $ K_{12,\gamma}$ is convex, narrowly lower semicontinuous, extended-real Borel, and contracts under marginalization over the remaining coordinates.

For affinity, apply Lemma~\ref{disc:thm:score-affinity} at each $n$ and commute the two countable suprema to obtain
\[
\begin{split}
 \lim_{j\to\infty}\mathcal{K}_{X,2}^{\rm off}(\Pi^j)
 &=\sup_j\sup_n\mathcal{K}_{X_n,2}(\Pi^j)
 =\sup_n\sup_j\mathcal{K}_{X_n,2}(\Pi^j)\\
 &=\sup_n\int
 \mathcal{K}_{X_n,2}(\rho^{\otimes2})\,\Pi(\md\rho)
 =\int\mathcal{K}_{X,2}^{\rm off}
 (\rho^{\otimes2})\,\Pi(\md\rho),
 \end{split}
\]
where only the final equality uses the monotone convergence theorem.
Summing over $k$ yields exact affinity for $K_{12,\gamma}$.

For application to weak solutions, suppose $I(\rho)<\infty$, so that $\rho$ has a density and its product puts no mass on $v_1=v_2$.
If in addition the extended functional is finite, the nested representing functions are compatible.
We can therefore choose $S_k\in L^2(\rho^{\otimes2})$ representing $X_k^2\rho^{\otimes2}$ away from the diagonal, with
\[
 \mathcal{K}_{X_k,2}^{\rm off}(\rho^{\otimes2})
 =\int|S_k|^2\md\rho^{\otimes2}.
\]
For a smooth density this is the quotient integral in \eqref{disc:eq:Kpair}.
\end{remark}

The following estimate shows that a short pair flow does not affect the limiting marginals.

\begin{lemma}\label{disc:lem:invisibility}
Let $\gamma\in[-3,1]$, $h_N\to0$, and let $F^N$ be exchangeable probability densities satisfying
\begin{equation}
\sup_N\int_{\R^{3N}}
 \frac1N\sum_{i=1}^N|v_i|^2F^N(\boldsymbol v)\md\boldsymbol v<\infty.\label{disc:eq:energy}
\end{equation}
For every fixed $j\ge2$ and every bounded Lipschitz $\psi:\R^{3j}\to\R$,
\[
\sup_{0\le s\le h_N}\Big|
 \int\psi\md(\ee^{sA_{12}}F^N)^{j|N}
 -\int\psi\md F^{j|N}
\Big|\to 0
\quad\text{as } h_N \to 0.
\]
The convergence is uniform over families satisfying \eqref{disc:eq:energy}.
\end{lemma}

\begin{proof}
Only $v_1,v_2$ move.
Conditional on $z=v_1-v_2=re$, the new direction is spherical Brownian motion at time $sr^\gamma$, and
\[
 \E_e|e_s-e|^2\le C\min\{1,sr^\gamma\}.
 \]
Thus
\begin{equation}
\E\big(|\Delta v_1|^2+|\Delta v_2|^2\mid r\big)
 \le Cr^2\min\{1,sr^\gamma\}.\label{disc:eq:pair-displacement}
\end{equation}
The same three-region argument works for every $\gamma$.
Fix $0<\varepsilon<R$.
On $\{r\le\varepsilon\}$, the deterministic bound $|\Delta v_1|+|\Delta v_2|\le2r$ gives an error at most $2\operatorname{Lip}(\psi)\varepsilon$.
On $\{\varepsilon<r\le R\}$, \eqref{disc:eq:pair-displacement} and the Cauchy--Schwarz inequality give an error at most $C_{\varepsilon,R,\gamma}\operatorname{Lip}(\psi)\sqrt{s}$.
On $\{r>R\}$, boundedness of $\psi$ and the energy estimate $\mathbb P(r>R)\le C/R^2$ give an error at most $C\|\psi\|_\infty/R^2$.
Consequently, for every $0\le s\le h_N$,
\[
 \Big|
  \int\psi\md(\ee^{sA_{12}}F^N)^{j|N}
  -\int\psi\md F^{j|N}
 \Big|
 \le2\operatorname{Lip}(\psi)\varepsilon
 +C_{\varepsilon,R,\gamma}\operatorname{Lip}(\psi)\sqrt{s}
 +\frac{C\|\psi\|_\infty}{R^2}.
\]
First choose $\varepsilon$ small and $R$ large, then use $s\le h_N\to0$ and let $N\to\infty$.
All constants use only the uniform energy bound, proving the asserted uniformity.
\end{proof}

\begin{proof}[Proof of Proposition~\ref{disc:thm:level3-diss}]
Let $\pi$ be a subsequential limit law, $\pi_t$ the law of $\rho_t$ under $\pi$, and define the de Finetti marginals
\begin{equation}
\pi_t^j:=\int_{\Pp(\R^3)}\eta^{\otimes j}\,\pi_t(\md\eta).\label{disc:eq:definetti}
\end{equation}

By exchangeability, sampling distinct particle labels gives the marginal $F_t^{j|N}$, whereas sampling from $\mu_t^N$ with replacement can repeat a label.
The probability of a repetition is at most $j(j-1)/(2N)$, so for every bounded continuous $\psi$ on $\R^{3j}$,
\begin{equation}
 \Big|
 \int\psi\md F_t^{j|N}
 -\E\int\psi\md(\mu_t^N)^{\otimes j}
 \Big|
 \le\frac{j(j-1)}N\|\psi\|_\infty.
 \label{disc:eq:empirical-marginal-comparison}
\end{equation}
Since the limiting paths are continuous, fixed-time evaluation preserves the path-law convergence, and \eqref{disc:eq:empirical-marginal-comparison} gives
\begin{equation}
F_t^{j|N}\rightharpoonup\pi_t^j\label{disc:eq:marginal-limit}
\end{equation}
for every fixed $j,t$.
Apply \cite[Theorem~5.7(2)]{HaurayMischler2014} to the marginal convergence \eqref{disc:eq:marginal-limit}.
By Proposition~\ref{disc:thm:finite-N}, we obtain, at every fixed time,
\[
 \E_\pi I(\rho_t)
 \le\liminf_{N\to\infty}\frac1N I(F_t^N)
 \le I(f_0).
\]
By the definition \eqref{disc:eq:fisher-information}, $\rho_t$ has a probability density for $\pi$-almost every path at each fixed $t$.
Tonelli's theorem gives the time-integrated inequality in \eqref{eq:level3-I}.

We now prove the directional second-order information bound.
Fix $j\ge2$ and denote the pair-flow averaged density by
\[
 G_t^N:=\frac1{h_N}\int_0^{h_N}\ee^{sA_{12}}F_t^N\md s.
\]
Convexity and contraction under marginalization give
\begin{equation}
 K_{12,\gamma}(G_t^{j|N})
 \le\frac1{h_N}\int_0^{h_N}
  K_{12,\gamma}((\ee^{sA_{12}}F_t^N)^{j|N})\md s
 \le\overline{ K}_{12,\gamma}^{\,h_N}(F_t^N).
\label{disc:eq:K-marginal}
\end{equation}

Averaging Lemma~\ref{disc:lem:invisibility} over $s\in[0,h_N]$ and using \eqref{disc:eq:marginal-limit} give $G_t^{j|N}\rightharpoonup\pi_t^j$ for every $t$.
Lower semicontinuity gives $K_{12,\gamma}(\pi_t^j)\le\liminf_{N\to\infty}K_{12,\gamma}(G_t^{j|N})$ for every $t$.
Since the functional is Borel measurable, Fatou's lemma yields
\[
 \int_0^T K_{12,\gamma}(\pi_t^j)\md t
 \le\liminf_{N\to\infty}\int_0^T
  K_{12,\gamma}(G_t^{j|N})\md t
  \le c_\gamma^{-1}I(f_0)
 \qquad(j\ge2),
\]
where the second inequality is due to \eqref{disc:eq:K-marginal} and Proposition~\ref{disc:thm:finite-N}.

Finally let $j\to\infty$.
By Lemma~\ref{disc:thm:score-affinity} and its extension in Remark~\ref{disc:rem:singular-affinity}, the limit is exactly the mean second-order functional.
The sequence is increasing by marginal contraction, so the monotone convergence theorem gives
\[
 \lim_{j\to\infty}\int_0^T K_{12,\gamma}(\pi_t^j)\md t
 =
 \int_0^T\lim_{j\to\infty} K_{12,\gamma}(\pi_t^j)\md t
 =
 \E_\pi\int_0^T
  K_\gamma(\rho_t^{\otimes2})\md t.
\]
This proves \eqref{disc:eq:level3-K}.
Together with \eqref{eq:level3-I}, Tonelli's theorem gives \eqref{disc:eq:pathwise-K}, hence finite $I(\rho_t)+K_\gamma(\rho_t^{\otimes2})$ for almost every time along $\pi$-almost every path.
At these times, $\rho_t$ has a density, and Remark~\ref{disc:rem:singular-affinity} gives the distributional $L^2$ representation.
\end{proof}

\section{Tightness of empirical measures\label{sec:Tightness-of-empirical}}
\label{sec:tightness}

To obtain the tightness of empirical measures $(\mu_{\cdot}^{N})_{N\ge2}$ in $D([0,T],\Pp(\R^3))$, we invoke Jakubowski's tightness criterion (cf. \cite[Theorem~3.1]{Jakubowski1986}).
In the present setting, it is enough to verify the following two conditions:
\begin{enumerate}
\item For any $\varepsilon>0$ there is a compact $K\subset\R^3$ such that
\[
 \sup_N\mathbb P\Big(\sup_{t\le T}\mu_t^N(K^c)>\varepsilon\Big)<\varepsilon.
\]
\item There exists a countable family $\mathcal G\subset C_c^\infty(\R^3)$, dense in $C_0(\R^3)$ for the uniform norm and closed under addition, such that for every $\varphi\in\mathcal G$, the real-valued processes $\langle\mu_\cdot^N,\varphi\rangle$ are tight in $D([0,T],\R)$.
\end{enumerate}
The maps $\nu\mapsto\langle\nu,\varphi\rangle$, $\varphi\in\mathcal G$, separate probability measures, as required by the criterion.
Throughout this section, we impose the assumptions of Theorem~\ref{thm:main} and fix $T>0$.
Since $h_N\to0$, discarding finitely many indices allows us to assume $0<h_N\le1$.

This first condition holds in our model by Chebyshev's inequality
and pathwise energy conservation:
\begin{align*}
\sup_{t\le T}\mu_{t}^{N}(\{|v|>R\}) & \le\frac{1}{R^{2}}\sup_{t\le T}\int|v|^{2}\mu_{t}^{N}(\md v)=\frac{1}{R^{2}}\frac{1}{N}\sum_{i}|v_{i}^{N}(0)|^{2}\quad\text{a.s.}
\end{align*}
Taking expectations yields
\[
\E\Big[\sup_{t\le T}\mu_{t}^{N}(\{|v|>R\})\Big]\le\frac{E_0}{R^{2}}\to0\quad\text{as }R\to+\infty.
\]

More explicitly, the sets $\{\nu\in\Pp(\R^3):\int|v|^2\nu(\md v)\le R\}$ are weakly compact.
Energy conservation and Markov's inequality give
\[
 \sup_N\mathbb P\Big(\sup_{t\le T}\langle\mu_t^N,|v|^2\rangle>R\Big)
 \le E_0/R,
\]
which verifies compact containment in $\Pp(\R^3)$.

For the second condition (tightness of the test processes), we use the martingale decomposition and the following estimates.
They rely on energy conservation and the near-diagonal control \eqref{eq:NDC}.
For $\gamma\ge0$, we also use the moment bound of Lemma~\ref{prop:M4} on the conserved energy event $\Omega_R^N$ defined in \eqref{eq:energy-localization}.
The second-order information enters later, when placing the limits for very soft potentials in a uniqueness class.

For $\varphi\in C_c^2(\R^3)$, define the symmetric collision kernel
\begin{equation}
 H_\varphi(v,w)
 :=
 \frac14a(v-w):
 \big(D^2\varphi(v)+D^2\varphi(w)\big)
 +\frac12b(v-w)\cdot
 \big(\nabla\varphi(v)-\nabla\varphi(w)\big),
\label{eq:Hphi}
\end{equation}
where $a(z) = |z|^{\gamma+2}(\operatorname{Id}-\hat{z}\otimes\hat{z})$ and $b(z)=-2|z|^\gamma z$.

\begin{lemma}\label{lem:MD}
For $\varphi\in C_c^2(\R^3)$, define the drift by
\[
 B_s^{N,\varphi}
 :=\big(\mathcal L^{N,h_N}\langle\mu^N,\varphi\rangle\big)(\bm v^N(s-)).
\]
Then 
\begin{equation}
 M_t^{N,\varphi}
 :=\langle\mu_t^N,\varphi\rangle-\langle\mu_0^N,\varphi\rangle
 -\int_0^tB_s^{N,\varphi}\md s
 \label{eq:empirical-martingale}
\end{equation}
is a square-integrable martingale starting from zero, with respect to the full natural filtration of the MBC process.
For every $T>0$, one has
\begin{equation}
\E[M^{N,\varphi}]_T\le\frac{C_T}{N},
 \qquad
 \sup_{s\le T}\E|B_s^{N,\varphi}|^p\le C_{T,p}\label{eq:MD}
\end{equation}
for some $p>1$ when $\gamma<0$.
When $\gamma\ge0$, for every $R<\infty$ one has instead
\begin{equation}
\E\big([M^{N,\varphi}]_T\mathbf1_{\Omega_R^N}\big)
 \le\frac{C_{T,R}}N,\qquad
 \sup_{s\le T}\E\big(
 |B_s^{N,\varphi}|^p\mathbf1_{\Omega_R^N}\big)
 \le C_{T,R,p}.\label{eq:MD-local}
\end{equation}
For $\gamma<0$, the bounds use only energy and \eqref{eq:NDC}.
At $\gamma=0$, the unlocalized quadratic-variation estimate also uses only energy; all other estimates for $\gamma\ge0$ use Lemma~\ref{prop:M4}.
\end{lemma}

\begin{proof}
Write $h=h_N$.
For each pair $(i,j)$, put $u=(v_i+v_j)/2$, $r=r_{ij}$, and $e=e_{ij}$, and set
\[
 U_{ij}^{N,\varphi}(e')
 :=N^{-1}\big[\varphi(u+re'/2)+\varphi(u-re'/2)\big].
\]
In the following two formulas, all pair quantities are evaluated at $\bm v^N(s-)$, and pairs with $r_{ij}=0$ contribute zero.
The generator gives the drift
\begin{equation}
 B_s^{N,\varphi}
 =\frac{2}{h(N-1)}\sum_{i<j}
 \Big((P_{hr_{ij}^\gamma}^{\S^2}-I)U_{ij}^{N,\varphi}\Big)(e_{ij}).
 \label{eq:pair-drift}
\end{equation}
The purely discontinuous martingale $M^{N,\varphi}$ has jumps $U_{ij}^{N,\varphi}(e')-U_{ij}^{N,\varphi}(e_{ij})$, since the compensating integral in \eqref{eq:empirical-martingale} is continuous.
Its predictable quadratic variation is
\begin{equation}
 \langle M^{N,\varphi}\rangle_t
 =\frac{2}{h(N-1)}\sum_{i<j}\int_0^t\int_{\S^2}
 \big|U_{ij}^{N,\varphi}(e')-U_{ij}^{N,\varphi}(e_{ij})\big|^2
 P_{hr_{ij}^\gamma}^{\S^2}(e_{ij},\md e')\md s.
 \label{eq:pair-quadratic-variation}
\end{equation}

We use two basic estimates:
\begin{itemize}
\item \emph{Deterministic jump bound}.
For any state and any $e'\in\S^2$,
\[
 \big|\Delta\langle\mu^N,\varphi\rangle\big|
 \le\frac{\|\nabla\varphi\|_{\infty}}{N}\,r\,|e'-e|
 \le\frac{2\|\nabla\varphi\|_{\infty}}{N}\,r.
\]
\item \emph{Second-moment of the rotation}.
If $e'\sim P_\tau^{\S^2}(e,\cdot)$ on $\S^2$, then
\[
 \E[|e'-e|^{2}\mid e]=2(1-\ee^{-2\tau})\le4\tau.
\]
In particular, with $\tau=hr^\gamma$, this is bounded by $C\min(1,hr^\gamma)$.
\end{itemize}
Plugging into \eqref{eq:pair-quadratic-variation} and using $\lambda=N/h$ give
\begin{align*}
\E[M^{N,\varphi}]_{T}
= \E\langle M^{N,\varphi}\rangle_{T}
\le\frac{4\|\nabla\varphi\|_{\infty}^{2}}{N}\int_{0}^{T}\E\Big[\frac{2}{N(N-1)}\sum_{i<j}r_{ij}^{\gamma+2}(t)\Big]\md t.
\end{align*}
We then consider three cases:
\begin{itemize}
\item For $\gamma\in[-3,-2)$, the elementary bound $r^{\gamma+2}\le r^{-2}+r^{2}$ yields a uniform bound by \eqref{eq:NDC} and the energy bound.
\item For $\gamma\in[-2,0]$, we have $r^{\gamma+2}\le1+r^{2}$ and obtain the same conclusion using only energy.
\item For $\gamma\in(0,1]$, multiply by $\mathbf1_{\Omega_R^N}$ and use Lemma~\ref{prop:M4}, noting that $\gamma+2<4$.
\end{itemize}

For the drift, direct tangential differentiation at $v_i=u+re/2$ and $v_j=u-re/2$ gives
\begin{equation}
\begin{aligned}
 \Delta_{\mathbb S^2}U_{ij}^{N,\varphi}(e)
 &=\frac{r^2}{4N}(I-e\otimes e):
 \big(D^2\varphi(v_i)+D^2\varphi(v_j)\big)\\
 &\quad-\frac rN e\cdot
 \big(\nabla\varphi(v_i)-\nabla\varphi(v_j)\big),\\
 Nr^\gamma\Delta_{\mathbb S^2}U_{ij}^{N,\varphi}(e)
 &=H_\varphi(v_i,v_j).
\end{aligned}
\label{eq:spherical-Landau}
\end{equation}
The gradient difference is bounded by $r\|D^2\varphi\|_\infty$, so both terms in the first identity are of order $r^2/N$.
On the sphere $\mathbb{S}^{2}$ the heat semigroup satisfies
\[
\|(P_{\tau}^{\S^2}-I)G\|_{\infty}\le\tau\,\|\Delta_{\mathbb{S}^{2}}G\|_{\infty},\quad G\in C^{2}(\mathbb{S}^{2}).
\]
Applying to \eqref{eq:pair-drift} with $G=U_{ij}^{N,\varphi}$ and $\tau=hr_{ij}^{\gamma}$
gives
\[
|B_{t}^{N,\varphi}|
\le\frac{2\lambda}{N(N-1)}\sum_{i<j}h r_{ij}^{\gamma}\,\|\Delta_{\mathbb{S}^{2}}U_{ij}^{N,\varphi}\|_{\infty}
\le\frac{C\lambda h}{N^2(N-1)}\sum_{i<j} r_{ij}^{\gamma+2}\,\|D^{2}\varphi\|_{\infty}.
\]
Using $\lambda h=N$, we obtain
\[
|B_{t}^{N,\varphi}|\le C\,\|D^{2}\varphi\|_{\infty}\cdot\frac{2}{N(N-1)}\sum_{i<j}r_{ij}^{\gamma+2}(t).
\]
By Jensen's inequality and exchangeability,
\[
 \E|B_t^{N,\varphi}|^p
 \le C\E r_{12}^{p(\gamma+2)}.
\]
If $\gamma+2<0$, choose $p>1$ with $-p(\gamma+2)\le2$ and use \eqref{eq:NDC}.
If $0<\gamma+2<2$, choose $p>1$ with $p(\gamma+2)\le2$ and use energy.
At $\gamma=-2$, any fixed $p>1$ works.
If $\gamma\ge0$, choose $p>1$ sufficiently close to one that $p(\gamma+2)<4$, multiply by $\mathbf1_{\Omega_R^N}$, and use Lemma~\ref{prop:M4}.
This proves \eqref{eq:MD}--\eqref{eq:MD-local}.
\end{proof}

We now prove tightness of the empirical measures and show that every subsequential limit has weakly continuous trajectories.
Continuity is precisely what upgrades $J_1$-convergence (natural for jump processes) to uniform-in-time convergence needed to pass to the Landau limit.

\begin{proposition}
\label{prop:tight}
Under the assumptions of Theorem~\ref{thm:main}, the empirical path laws are $C$-tight in $D([0,T],\Pp(\R^3))$.
Consequently,
there exist a subsequence $(N_{\ell})$, a probability space $(\tilde{\Omega},\tilde{\mathcal{F}},\tilde{\mathbb{P}})$,
and processes $(\tilde{\mu}_{\cdot}^{N_{\ell}})_{\ell\ge1}$ and $\tilde{\mu}_{\cdot}$
on $(\tilde{\Omega},\tilde{\mathcal{F}},\tilde{\mathbb{P}})$, such
that
\begin{itemize}
\item $\Law(\tilde{\mu}_{\cdot}^{N_{\ell}})=\Law(\mu_{\cdot}^{N_{\ell}})$
for all $\ell$;
\item $\tilde{\mu}_{\cdot}^{N_{\ell}}\to\tilde{\mu}_{\cdot}$ almost surely
in $D([0,T],\Pp(\R^{3}))$ with the $J_{1}$-topology.
\end{itemize}
Moreover, every such subsequential limit $\tilde{\mu}_{\cdot}$ has
continuous trajectories in $\Pp(\R^{3})$ for the weak topology.
In particular, for every $\varphi\in C_{c}^{2}(\R^{3})$,
\[
\sup_{t\in[0,T]}\big|\langle\tilde{\mu}_{t}^{N_{\ell}}-\tilde{\mu}_{t},\varphi\rangle\big|\to0\quad\text{almost surely as }\ell\to\infty.
\]
\end{proposition}

\begin{proof}
We first verify $C$-tightness of the real-valued test processes.
Fix $\varphi\in\mathcal G$.
The initial values are bounded by $\|\varphi\|_\infty$.
For $\gamma<0$, H\"older's inequality gives
\[
 \bigg|\int_s^t B_r^{N,\varphi}\md r\bigg|
 \le |t-s|^{1-1/p}
 \bigg(\int_0^T|B_r^{N,\varphi}|^p\md r\bigg)^{1/p}.
\]
Lemma~\ref{lem:MD} bounds the $L^p$ norm on the right in probability, uniformly in $N$.
Together with the bounded initial values, the Arzel\`a--Ascoli theorem gives tightness of $\langle\mu_0^N,\varphi\rangle+\int_0^\cdot B_r^{N,\varphi}\md r$ in $C([0,T],\R)$.
Doob's maximal inequality and Lemma~\ref{lem:MD} give
\[
 \E\sup_{t\le T}|M_t^{N,\varphi}|^2
 \le4\E[M^{N,\varphi}]_T\le C_T/N.
\]
Thus the martingale tends to zero in probability uniformly on $[0,T]$, and the test processes are $C$-tight.

For $\gamma\ge0$, apply the same argument to $t\mapsto\mathbf1_{\Omega_R^N}\langle\mu_t^N,\varphi\rangle$, using the localized estimates of Lemma~\ref{lem:MD};
the event depends only on the initial state, so $\mathbf1_{\Omega_R^N}M^{N,\varphi}$ remains a martingale.
These localized observables are $C$-tight for each $R$, and they agree with the original observables on $\Omega_R^N$.
Since $\mathbb P((\Omega_R^N)^c)\le E_0/R$, letting $R\to\infty$ proves scalar $C$-tightness for $\gamma\ge0$ as well.

Choosing a countable
dense family $\mathcal{G}\subset C_{c}^{\infty}(\R^{3})$ and invoking
Jakubowski's criterion yields tightness of $(\mu_{\cdot}^{N})_{N\ge2}$
in $D([0,T],\Pp(\R^{3}))$. Prokhorov's theorem and Skorokhod's
theorem (cf. \cite{Billingsley1999}) then produce a subsequence $(N_{\ell})$
and a coupling $(\tilde{\mu}_{\cdot}^{N_{\ell}},\tilde{\mu}_{\cdot})$
with the same laws and almost sure convergence in $J_{1}$-topology
as desired.

Next we verify the continuity of subsequential limits via $C$-tightness (cf. \cite[Chapter~VI]{JacodShiryaev2003}).
For each $\varphi\in\mathcal G$, composition with $\nu\mapsto\langle\nu,\varphi\rangle$ preserves $J_1$ convergence.
Scalar $C$-tightness therefore makes $t\mapsto\langle\tilde\mu_t,\varphi\rangle$ continuous almost surely.
Since $\mathcal G$ is countable, this holds for all its members on one event of probability one.
At any $t>0$, continuity gives $\langle\tilde\mu_t,\varphi\rangle=\langle\tilde\mu_{t-},\varphi\rangle$ for every $\varphi\in\mathcal G$.
The separating property then gives $\tilde\mu_t=\tilde\mu_{t-}$, so $\tilde\mu$ is weakly continuous.
Since the limit is continuous, $J_1$-convergence strengthens to uniform convergence on $[0,T]$ (cf. \cite[Section~12]{Billingsley1999}).
Finally, for an arbitrary $\psi\in C_c^2(\R^3)$, choose $\varphi\in\mathcal G$ with $\|\psi-\varphi\|_\infty<\varepsilon$; then
\[
 \sup_{t\le T}\big|\langle\tilde\mu_t^{N_\ell}-\tilde\mu_t,\psi\rangle\big|
 \le2\varepsilon+
 \sup_{t\le T}\big|\langle\tilde\mu_t^{N_\ell}-\tilde\mu_t,\varphi\rangle\big|.
\]
Letting $\ell\to\infty$, and $\varepsilon\downarrow0$, proves the claim.
\end{proof}

\section{Identification of the Landau limit\label{sec:Identification-of-the}}
\label{sec:compactness}

Fix $T>0$.
Let $\gamma\in[-3,1]$ and let the stepsize satisfy $h=h_N\to0$ as $N\to\infty$.
Consider the MBC process with initial law $F_0^N=f_0^{\otimes N}$ under the assumptions of Theorem~\ref{thm:main}.
Let $\pi$ be an arbitrary limit law provided by Proposition~\ref{prop:tight}, and let $\rho=(\rho_t)_{t\le T}$ be its canonical path.
The de Finetti marginals $\pi_t^j$ satisfy \eqref{disc:eq:definetti}--\eqref{disc:eq:marginal-limit}.

\begin{lemma}\label{lem:remainder}
For $\varphi\in C_c^4(\R^3)$, write
\begin{equation}
 \mathcal L^{N,h}\langle\mu^N,\varphi\rangle
 =\int H_\varphi\md(\mu^N)^{\odotimes2}
   +\mathcal R^{N,h}_\varphi,
 \qquad
 (\mu^N)^{\odotimes2}
 :=\frac1{N(N-1)}\sum_{i\ne j}\delta_{(v_i,v_j)}.
 \label{eq:decomp}
\end{equation}
Here $(\mu^N)^{\odotimes2}$ is the empirical measure obtained by sampling two distinct particle labels in order.
Then, for every $T>0$,
\begin{equation}
 \int_0^T\big|\mathcal R^{N,h_N}_\varphi(\bm v^N(s))\big|\md s
 \longrightarrow0\quad\text{in probability}.
 \label{eq:integrated-Rto0}
\end{equation}
\end{lemma}

\begin{proof}
Write $h=h_N$ and take $N$ large enough that $h\le1$.
For the spherical heat semigroup $P_{\tau}^{\S^2}=\exp(\tau\Delta_{\S^2})$, the Taylor expansion with remainder gives
\[
 R_{\tau}G:=P_{\tau}^{\S^2}G-G-\tau\Delta_{\S^2}G
\]
and
\[
 \|R_{\tau}G\|_{\infty}
 \le \min\Big\{\frac{1}{2}\tau^{2}\|\Delta_{\S^2}^{2}G\|_{\infty},
 2\tau\|\Delta_{\S^2}G\|_{\infty}\Big\}
\]
for all $\tau\ge0$ and $G\in C^4(\S^2)$.
We apply this with $G=U_{ij}^{N,\varphi}$ and $\tau_{ij}:=h r_{ij}^{\gamma}$, where $r_{ij}=|v_i-v_j|$.
Apply the spherical Laplacian once more to \eqref{eq:spherical-Landau}.
The even-order derivative terms have factors $r^2$ and $r^4$.
The odd-order terms contain $r(\nabla\varphi(v_i)-\nabla\varphi(v_j))$ or $r^3(D^3\varphi(v_i)-D^3\varphi(v_j))$.
The mean value theorem supplies one additional factor of $r$ in each difference, so
\[
\begin{aligned}
 \|\Delta_{\S^2}U_{ij}^{N,\varphi}\|_\infty
 &\le\frac{C}{N}r^2\|D^2\varphi\|_\infty,\\
 \|\Delta_{\S^2}^2U_{ij}^{N,\varphi}\|_\infty
 &\le\frac{C}{N}\big(r^2\|D^2\varphi\|_\infty+r^4\|D^4\varphi\|_\infty\big).
\end{aligned}
\]
Hence, for every $e\in\S^2$,
\begin{align*}
 \big|(R_{\tau_{ij}}U_{ij}^{N,\varphi})(e)\big|
 &\le \min\Big\{\tfrac{1}{2}\tau_{ij}^{2}
 \|\Delta_{\S^2}^{2}U_{ij}^{N,\varphi}\|_{\infty},
 2\tau_{ij}\|\Delta_{\S^2}U_{ij}^{N,\varphi}\|_{\infty}\Big\}\\
 &\le \frac{C\|\varphi\|_{C^{4}}}{N}
 \min\Big\{h^{2}r_{ij}^{2\gamma}\big(r_{ij}^{2}+r_{ij}^{4}\big),
 h r_{ij}^{\gamma+2}\Big\}.
\end{align*}

By \eqref{eq:spherical-Landau}, multiplying the first-order term by the jump rate and averaging over pairs gives the Landau term in \eqref{eq:decomp}.
The remaining term is therefore
\begin{equation}
 \mathcal R^{N,h}_\varphi(\bm v)
 =\frac{2}{h(N-1)}\sum_{i<j}
 \big(R_{hr_{ij}^\gamma}U_{ij}^{N,\varphi}\big)(e_{ij}),
 \label{eq:pair-remainder}
\end{equation}
with zero contribution from pairs with $r_{ij}=0$.
From this definition we obtain
\[
 |\mathcal R^{N,h}_\varphi(\bm v)|
 \le \frac{C\|\varphi\|_{C^{4}}}{h}\cdot
 \frac{2}{N(N-1)}\sum_{i<j}
 \min\Big\{h^{2}r_{ij}^{2\gamma}\big(r_{ij}^{2}+r_{ij}^{4}\big),
 h r_{ij}^{\gamma+2}\Big\}.
\]
By exchangeability,
\[
 \E|\mathcal R^{N,h}_\varphi(\bm v^N(t))|
 \le C\E\min\big\{h(r^{2\gamma+2}+r^{2\gamma+4}),r^{\gamma+2}\big\}.
\]
Here $r=|v_1^N(t)-v_2^N(t)|$.

We estimate this expectation by splitting the relative speed into three regions.
Fix $0<\varepsilon<1<L$.
On $r<\varepsilon$, use the second bound in the minimum and $r^{\gamma+2}\le\varepsilon^{\gamma+4}r^{-2}$.
The near-diagonal estimate \eqref{eq:NDC} bounds this contribution by $C\varepsilon^{\gamma+4}$.
On $\varepsilon\le r\le L$, the first bound is at most $C_{\varepsilon,L}h$.
On $r>L$, use the second bound again: $r^{\gamma+2}\le L^\gamma r^2$ if $\gamma<0$, and $r^{\gamma+2}\le L^{\gamma-2}r^4$ otherwise.
Energy conservation controls the second moment; for $\gamma\ge0$, Lemma~\ref{prop:M4} controls the fourth moment on $\Omega_A^N$.
Thus, for $\gamma<0$,
\[
 \sup_{t\le T}\E|\mathcal R^{N,h}_\varphi(\bm v^N(t))|
 \le C\varepsilon^{\gamma+4}
       +C_{\varepsilon,L}h+CL^\gamma,
\]
whereas for $\gamma\ge0$ and each fixed $A$,
\[
 \sup_{t\le T}\E\big[
 |\mathcal R^{N,h}_\varphi(\bm v^N(t))|\mathbf1_{\Omega_A^N}\big]
 \le C\varepsilon^{\gamma+4}
       +C_{\varepsilon,L}h+C_{T,A}L^{\gamma-2}.
\]

Integrate in time and let $N\to\infty$, then $L\to\infty$ and $\varepsilon\downarrow0$.
The time integral tends to zero in $L^1$ when $\gamma<0$.
When $\gamma\ge0$, the same holds after multiplication by $\mathbf1_{\Omega_A^N}$ for each fixed $A$.
Markov's inequality proves \eqref{eq:integrated-Rto0} directly when $\gamma<0$.
For $\gamma\ge0$, the same inequality and $\mathbb P((\Omega_A^N)^c)\le E_0/A$ give, for every $\delta>0$,
\[
\begin{aligned}
 &\mathbb P\Big(\int_0^T|\mathcal R^{N,h_N}_\varphi(\bm v^N(s))|\md s>\delta\Big)\\
 &\qquad\le\frac{E_0}{A}
 +\frac1\delta\E\Big[\mathbf1_{\Omega_A^N}
 \int_0^T|\mathcal R^{N,h_N}_\varphi(\bm v^N(s))|\md s\Big].
\end{aligned}
\]
First let $N\to\infty$ for fixed $A$, then let $A\to\infty$ to conclude.
\end{proof}

\begin{proposition}\label{prop:identify}
Every limit law is supported on weak solutions of \eqref{eq:Landau} starting from $f_0$.
Explicitly, for $\pi$-almost every $\rho$,
\begin{equation}
\begin{aligned}
 \langle\rho_t,\varphi\rangle
 &=\langle f_0,\varphi\rangle
 +\frac14\int_0^t\iint
 a(v-w):\big(\nabla^2\varphi(v)+\nabla^2\varphi(w)\big)
 \rho_s(\md v)\rho_s(\md w)\md s\\
 &\quad+\frac12\int_0^t\iint
 b(v-w)\cdot\big(\nabla\varphi(v)-\nabla\varphi(w)\big)
 \rho_s(\md v)\rho_s(\md w)\md s.
\end{aligned}\label{eq:weak}
\end{equation}
for all $\varphi\in C_c^4(\R^3)$ and $t\le T$.
Moreover, almost every such path satisfies
\begin{equation}
\int |v|^2\rho_t(\md v)=
 \int |v|^2f_0(v)\md v\qquad(t\le T),\label{eq:path-energy}
\end{equation}
and, if $\gamma\ge0$,
\begin{equation}
\int_0^T\int\langle v\rangle^4\rho_t(\md v)\md t<\infty.\label{eq:path-M4}
\end{equation}
For every $\gamma\in[-3,1]$, \eqref{eq:weak} consequently extends to the standard symmetric $C_b^2$ weak-solution class.
\end{proposition}

\begin{remark}\label{rem:limit-regularity}
The information bounds of Proposition~\ref{disc:thm:level3-diss} apply to every limit law obtained above.
Tonelli's theorem then gives, for $\pi$-almost every path,
\[
\int_0^T\big[
 I(\rho_t)+ K_\gamma(\rho_t^{\otimes2})
 \big]\md t<\infty.
\]
For $\gamma\in[-2,0)$, the time-integrated Fisher bound controls the singular convolution through \eqref{eq:fisher-hardy} and supplies the estimate needed for uniqueness in Proposition~\ref{prop:unique}.
For $\gamma\in[-3,-2)$, the second-order term is a key input in Proposition~\ref{vs:thm:closure} which uses these information bounds to obtain $L^1(0,T;L^\infty)$ regularity and then uniqueness.
\end{remark}

\begin{proof}[Proof of Proposition~\ref{prop:identify}]
Let $(N_\ell)$ be a subsequence whose empirical path laws converge to $\pi$, as provided by Proposition~\ref{prop:tight}.
Apply the Skorokhod representation theorem \cite{Billingsley1999} jointly to these paths and the initial empirical energies $M_2^N(0)$, which converge to $E_0$ by the law of large numbers.
We retain the same notation on this representation and write $N$ for $N_\ell$.
Then $\mu^N\to\rho$ in the path space and $M_2^N(0)\to E_0$ almost surely.
Since $\rho$ is weakly continuous, the path convergence is uniform in time for the weak topology.

Fix $\varphi\in C_c^4(\R^3)$ and $t>0$.
For each $N$ we have the semimartingale decomposition from Lemmas~\ref{lem:MD} and~\ref{lem:remainder},
\[
 \langle\mu_t^N,\varphi\rangle
 =\langle\mu_0^N,\varphi\rangle+M_t^{N,\varphi}
 +\int_0^t\iint H_\varphi(v,w)\md(\mu_s^N)^{\odotimes2}\md s
 +\int_0^t\mathcal R^{N,h_N}_\varphi(\bm v^N(s))\md s.
\]
The drift and remainder are symmetric functions of the particle configuration, hence functions of the empirical measure.
We pass to the limit along $(N_\ell)$, term by term.

\emph{Step 1. Initial data, martingale, and remainder.}

By independence at $t=0$,
\[
 \langle\mu_0^N,\varphi\rangle\to\langle f_0,\varphi\rangle
 \quad\text{in }L^2\text{ and in probability}.
\]
By Doob's maximal inequality and Lemma~\ref{lem:MD}, the martingale tends to zero in probability uniformly on $[0,T]$.
For $\gamma\ge0$, use the estimate localized on $\Omega_A^N$ with any fixed $A>E_0$; its indicator tends to one almost surely by the joint energy convergence.
Lemma~\ref{lem:remainder} gives
\[
 \int_0^T|\mathcal R^{N,h_N}_\varphi(\bm v^N(s))|\md s
 \longrightarrow0\quad\text{in probability}.
\]

To identify the limit of the Landau term, we first record the moment bounds needed below.
The particle energy is conserved.
By lower semicontinuity of truncated quadratic moments, $\int |v|^2\rho_t(\md v)\le E_0$ simultaneously for all $t$.
When $\gamma\ge0$, keep this energy threshold $A>E_0$ fixed.
Lemma~\ref{prop:M4} and Fatou's lemma bound the expected time integral of the fourth moment of $\rho$.
Tonelli's theorem therefore gives \eqref{eq:path-M4}.

\emph{Step 2. Passing to the limit away from the diagonal.}

We first pass to the limit away from the diagonal and at bounded velocities.
For fixed diagonal and velocity cutoffs, denote the resulting bounded continuous kernel by $H_\varphi^{\varepsilon,R}$.
Almost sure convergence of $\mu_s^N$ and the dominated convergence theorem give
\[
 \int_0^t\iint H_\varphi^{\varepsilon,R}
 \md(\mu_s^N)^{\otimes2}\md s
 \longrightarrow
 \int_0^t\iint H_\varphi^{\varepsilon,R}
 \md\rho_s^{\otimes2}\md s
\]
uniformly for $t\in[0,T]$.
The off-diagonal statistic differs from the tensor square here by at most $2\|H_\varphi^{\varepsilon,R}\|_\infty/N$.

\emph{Step 3. Control near the diagonal.}

	On the set $\{|v-w|<\varepsilon\}$, we use the structure of $a$ and $b$.
	Since $\varphi$ is smooth and compactly supported,
	\[
	|H_{\varphi}(v,w)|
	\lesssim |v-w|^{\gamma+2}
	\]
	for $|v-w|$ small.
	Moreover, for $\gamma\ge-3$,
	\[
	|v-w|^{\gamma+2}
	\le \varepsilon^{\gamma+4}\,|v-w|^{-2}
	\quad\text{whenever }|v-w|\le\varepsilon.
	\]

The particle contribution is controlled by the near-diagonal estimate \eqref{eq:NDC}.
For the limit, use the Fisher-information bound of Proposition~\ref{disc:thm:level3-diss} and \eqref{eq:fisher-hardy}.
These bounds give
\begin{align*}
 \sup_N\E\int_0^T\iint
 |H_\varphi|\mathbf1_{\{|v-w|\le\varepsilon\}}
 \md(\mu_s^N)^{\odotimes2}\md s
 &\le C_{\varphi,T}\varepsilon^{\gamma+4}I(f_0),\\
 \E_\pi\int_0^T\iint
 |H_\varphi|\mathbf1_{\{|v-w|\le\varepsilon\}}
 \md\rho_s^{\otimes2}\md s
 &\le C_{\varphi,T}\varepsilon^{\gamma+4}I(f_0).
\end{align*}
Letting $\varepsilon\to0$ yields that the near-diagonal contributions vanish in $L^1$, uniformly in $N$.

\emph{Step 4. Removing the velocity cutoff.}

It remains to remove the velocity cutoff.
Since $\varphi$ is compactly supported, $H_\varphi$ vanishes unless at least one variable lies in a fixed ball, and
\[
 |H_\varphi(v,w)|\le C_\varphi(1+|v-w|^{\gamma+2}).
\]
For $\gamma<0$, its growth is strictly subquadratic and energy gives uniform integrability.
For $\gamma\ge0$, keep the energy threshold $A>E_0$ fixed and let $R$ denote only the velocity-cutoff radius.
Take $R$ larger than twice the support radius of the derivatives of $\varphi$.
For $H_\varphi$ to be nonzero, at least one variable must remain in that fixed support.
Hence
\[
 |H_\varphi(v,w)|
 \mathbf1_{\{\max(|v|,|w|)>R\}}
 \le C_\varphi R^{\gamma-2}
 \big(\langle v\rangle^4+\langle w\rangle^4\big).
\]
Consequently, pair averaging and Lemma~\ref{prop:M4} give
\[
 \sup_N\E\bigg[
 \mathbf1_{\Omega_A^N}\int_0^T\iint
 |H_\varphi(v,w)|
 \mathbf1_{\{\max(|v|,|w|)>R\}}
 \md(\mu_s^N)^{\odotimes2}\md s\bigg]
 \le C_{\varphi,T,A}R^{\gamma-2}.
\]
Because $\gamma-2\le-1$, this tends to zero as $R\to\infty$.
The analogous limit-path tail tends to zero by \eqref{eq:path-M4}.

\emph{Step 5. Identification.}

Combining Steps 2--4, we first let $N\to\infty$, then $R\to\infty$, and finally $\varepsilon\downarrow0$.
For $\gamma\ge0$, this is done on $\Omega_A^N$, whose indicator tends to one for the fixed $A>E_0$.
Consequently, the time-integrated Landau term converges in probability to its counterpart for $\rho$.
Finally, passing to the limit in the semimartingale identity and using the convergences established in Step~1, we obtain
\[
 \langle\rho_t,\varphi\rangle
 =\langle f_0,\varphi\rangle+
 \int_0^t\iint H_\varphi(v,w)\rho_s(\md v)\rho_s(\md w)\md s
\]
almost surely for the fixed $\varphi,t$.
Expanding $H_\varphi$ in terms of $a$ and $b$ gives precisely the claimed weak formulation in Proposition~\ref{prop:identify}.
First take a countable family of tests, dense in the $C^2$ norm on each fixed compact support, and rational times.
Weak continuity of $\rho$ and the preceding integral bounds then extend the identity to all $C_c^4$ tests and all $t\le T$ on one event of probability one.

Finally, for every $\gamma\in[-3,1]$,
\[
 |H_\varphi(v,w)|\le C_\varphi
 \big(1+|v-w|^{\gamma+2}\big),
 \qquad \varphi\in C_b^2(\R^3).
\]
For $\gamma<0$, the near-diagonal part is controlled by $|v-w|^{-2}$ and the Fisher-information bound of Proposition~\ref{disc:thm:level3-diss}, exactly as above, while the large-velocity part is controlled by the energy upper bound proved above.
For $\gamma\ge0$, use \eqref{eq:path-M4}.
These bounds justify the standard cutoff and mollification passage from $C_c^4$ to $C_b^2$.
Hence \eqref{eq:weak} holds simultaneously for all $C_b^2$ tests on the same full-probability event.

To prove energy conservation, take a radial cutoff $\chi\in C_c^\infty(\R^3)$, with $0\le\chi\le1$, which equals one near the origin, and set $\varphi_R(v)=|v|^2\chi(v/R)$.
Its Hessian is bounded uniformly in $R\ge1$, and the mean value theorem gives
\[
 |H_{\varphi_R}(v,w)|\le C|v-w|^{\gamma+2}.
\]
The right-hand side is integrable against $\rho_s\otimes\rho_s\md s$ on $[0,T]$.
For $\gamma<0$, use the energy bound at large relative velocities and the Fisher-information bound of Proposition~\ref{disc:thm:level3-diss}, Tonelli's theorem, and \eqref{eq:fisher-hardy} near the diagonal.
For $\gamma\ge0$, use \eqref{eq:path-M4}.
For $v\ne w$, $H_{\varphi_R}(v,w)\to H_{|v|^2}(v,w)=0$, since
\[
 \operatorname{tr}a(z)=2|z|^{\gamma+2},
 \qquad b(z)\cdot z=-2|z|^{\gamma+2}.
\]
The dominated convergence theorem applied in \eqref{eq:weak} therefore gives \eqref{eq:path-energy} at every $t$.
The same argument with $v\chi(v/R)$ gives momentum conservation:
\[
 \int v\,\rho_t(\md v)=\int v f_0(v)\md v\qquad(t\le T).
\]
The proof is complete.
\end{proof}

\section{Uniqueness and proof of the main theorem}
\label{sec:identification-convergence}\label{sec:uniqueness-completion}

By Proposition~\ref{prop:identify}, every subsequential limit is supported on weak solutions of the spatially homogeneous Landau equation with initial data $f_0$.
To complete the proof of Theorem~\ref{thm:main}, it remains to show
that the whole sequence $(\mu_{\cdot}^{N})_{N\ge2}$
converges, which reduces to proving uniqueness of the limiting dynamics.

Global convergence follows once uniqueness holds in a class containing all such limits.
The available results are summarized in Remark~\ref{rem:uniqueness-classes}.
For soft potentials, we use a comparison solution; Proposition~\ref{vs:thm:closure} handles the range $\gamma\in[-3,-2)$ and is proved in Section~\ref{sec:very-soft-closure}.
We use the weighted velocity moments
\begin{equation}
 M_\ell(\rho):=\int_{\R^3}\langle v\rangle^\ell \rho(\md v),
 \qquad \rho\in\Pp(\R^3).
 \label{eq:moment-definition}
\end{equation}

\begin{lemma}
\label{vs:prop:comparison}
Let $\gamma\in[-3,0)$, and let $f_0$ satisfy the assumptions of Theorem~\ref{thm:main}.
There exists $g\in C_w([0,\infty),\Pp(\R^3))$ starting from $f_0$, satisfying \eqref{eq:weak} for every $C_b^2$ test, conserving energy and momentum, and obeying
\[
 I(g_t)\le I(f_0),\qquad M_2(g_t)=M_2(f_0)
 \quad(t\ge0).
\]
In particular, $g_t$ has a density at every time.
If $\gamma\in[-3,-2)$, then, for every $T>0$,
\begin{equation}
 \sup_{t\le T}M_m(g_t)<\infty,\qquad
 \|g_t\|_\infty\le C_T(1+t^{-3/4})
 \quad(0<t\le T),
 \label{vs:eq:g-Linf}
\end{equation}
where $m$ is the moment order in Theorem~\ref{thm:main}.
Consequently, $g\in L^1(0,T;L^\infty(\R^3))$ in this range.
\end{lemma}

\begin{proof}
Cut off and mollify $\sqrt{f_0}$, using nonnegative functions, and normalize in $L^2$.
This gives nonnegative $u_n\in C_c^\infty(\R^3)$, with $\|u_n\|_2=1$, converging to $\sqrt{f_0}$ in $H^1$ and in $L^2(\langle v\rangle^m\md v)$,
where $m=2$ for $\gamma\in[-2,0)$ and is the moment order in Theorem~\ref{thm:main} otherwise.
Set $f_0^n=u_n^2$.
Using the Cauchy--Schwarz inequality and $I(u_n^2)=4\|\nabla u_n\|_2^2$, we obtain
\[
 f_0^n\to f_0\ \hbox{in }L^1,\qquad
 M_m(f_0^n)\to M_m(f_0),\qquad I(f_0^n)\to I(f_0).
\]
Each $f_0^n$ is bounded by a Maxwellian, so \cite[Theorems~1.1 and~1.2]{GuillenSilvestre2025} give a global classical solution with nonincreasing Fisher information.
After rescaling time by $t\mapsto t/2$ to match the normalization in \eqref{eq:Landau}, denote this solution by $g^n$.
These Maxwellian bounds supply the classical construction, while the estimates below are uniform in $n$ by the common initial Fisher information and moment bounds.
Thus $I(g_t^n)\le I(f_0^n)$ and $M_2(g_t^n)=M_2(f_0^n)$ uniformly in time.

In the very soft range, the energy and Fisher bounds give uniform entropy bounds.
We apply the moment estimate of Carrapatoso--Desvillettes--He \cite[Lemma~8]{CarrapatosoDesvillettesHe2017} to the classical solutions and rescale the resulting estimate in time.
This gives $\sup_{n,t\le T}M_m(g_t^n)<\infty$.
Its constants depend only on $\gamma,m,T$ and the common initial mass, energy, entropy, and $m$-moment bounds.
Sobolev's inequality gives a uniform $L^3$ bound, and H\"older's inequality yields
\[
 \int\langle v\rangle^{m/2}(g_t^n)^2\md v
 \le M_m(g_t^n)^{1/2}\|g_t^n\|_{L^3}^{3/2}.
\]
Interpolating $L^\infty_tL^3_v$ with the propagated $m$-moment thus yields a uniform weighted $L^2$ bound.
By the $L_k^p\to L^\infty$ smoothing criterion of \cite[Theorem~2.11]{Ji2024FisherDissipationLandau}, the classical solutions satisfy
\[
 p\in\Big(\frac32,\frac3{3+\gamma}\Big),\quad
 k>\frac1p\max\Big\{3p-3,-\frac32\gamma-1\Big\}
 \quad\Longrightarrow\quad
 \|g_t^n\|_\infty\lesssim1+t^{-\frac3{2p}},
\]
where $\|f\|_{L_k^p}=\|\langle v\rangle^k f\|_{L^p}$ and the upper endpoint for $p$ is understood as $+\infty$ at $\gamma=-3$.
This improves the weighted criterion in Silvestre's estimate \cite[Theorem~3.8]{Silvestre2017UpperBounds}.
We apply it with $p=2$ and $k=m/4$.
Here $p=2$ lies in the admissible range, and $m>3(2-\gamma)$ implies the required condition $m>\max\{6,-3\gamma-2\}$.
The estimate gives
\[
 \|g_t^n\|_\infty\le C_T(1+t^{-3/4})
 \qquad(0<t\le T).
\]
The constant depends only on the common mass, energy, entropy, and weighted $L^2$ bounds above, so $C_T$ is uniform in $n$; the time rescaling changes only this constant.

For every $\gamma\in[-3,0)$, the energy bound gives tightness in velocity.
Use the weak equation with the uniform energy and Fisher bounds and \eqref{eq:fisher-hardy} to obtain a uniform time modulus against smooth tests.
More precisely, for every $\varphi\in C_c^\infty(\R^3)$, these bounds give
\[
 |\langle g_t^n-g_s^n,\varphi\rangle|
 \le C_\varphi|t-s|,
\]
where $C_\varphi$ depends only on the common energy and Fisher bounds.
The Arzel\`a--Ascoli theorem and a diagonal subsequence therefore give convergence to a path $g$ in $C_w([0,T],\Pp(\R^3))$ on every finite interval.
The estimate
\[
 \iint_{|v-w|<\varepsilon}
 |H_\varphi(v,w)|g_t^n(v)g_t^n(w)\md v\md w
 \le C_\varphi\varepsilon^{\gamma+4}I(g_t^n)
\]
controls the diagonal cutoff, while the subquadratic growth of $H_\varphi$ and the energy bound control its tails.
The weak identity passes to $g$, for every $C_b^2$ test, as in Proposition~\ref{prop:identify}.
The convergence of the initial data gives $g_0=f_0$, and lower semicontinuity gives $I(g_t)\le I(f_0)$ at every time.
In the very soft range, lower semicontinuity gives the moment bound in \eqref{vs:eq:g-Linf}.
The uniform $L^\infty$ estimate also passes to the limit: for every fixed $t>0$ and every nonnegative $\psi\in C_c(\R^3)$, weak convergence gives
\[
 \langle g_t,\psi\rangle
 \le C_T(1+t^{-3/4})\int\psi(v)\md v.
\]
This proves the $L^\infty$ bound at every positive time and completes \eqref{vs:eq:g-Linf}.
Finally, the cutoff argument in Proposition~\ref{prop:identify} proves conservation of energy and momentum.
\end{proof}

For $\gamma\in[-3,0)$, fix a comparison solution $g$ from Lemma~\ref{vs:prop:comparison}.

Let $\gamma\in[-3,0)$ and let $h\in C_w([0,T],\Pp(\R^3))$ have bounded $M_2(h_t)$ and $\int_0^T I(h_t)\md t<\infty$.
If $h$ satisfies \eqref{eq:weak} for every $C_b^2$ test, it also satisfies the unsymmetrized weak formulation.
Indeed, the coefficients obey
\[
 |a(z)|+|b(z)|
 \le C_\gamma\big(1+|z|^2+|z|^{-2}\big)
 \qquad(z\ne0),
\]
so \eqref{eq:fisher-hardy} gives, for almost every $t$,
\begin{equation}
 \iint\big(|a(v-w)|+|b(v-w)|\big)
 h_t(\md v)h_t(\md w)
 \le C_\gamma\big(1+M_2(h_t)+I(h_t)\big).
 \label{vs:eq:ordinary-coefficient-bound}
\end{equation}
Hence the unsymmetrized collision kernel
\begin{equation}
 L_\varphi(v,w)
 :=\tfrac12 a(v-w):D^2\varphi(v)
     +b(v-w)\cdot\nabla\varphi(v)
 \label{vs:eq:unsymmetrized-kernel}
\end{equation}
is absolutely integrable in space and time for every $\varphi\in C_b^2(\R^3)$.
Exchanging $v,w$ therefore gives $\iint L_\varphi\,h_t^{\otimes2}=\iint H_\varphi\,h_t^{\otimes2}$, and converts \eqref{eq:weak} into the unsymmetrized weak identity.

\begin{proposition}
\label{vs:thm:closure}
Let $\gamma\in[-3,-2)$ and $T>0$, and assume that $f_0$ satisfies the initial-data assumptions of Theorem~\ref{thm:main}.
Suppose that $f\in C_w([0,T],\Pp(\R^3))$ starts from $f_0$, satisfies the weak identity \eqref{eq:weak} for every $C_b^2$ test, and obeys
\begin{equation}
 \int_0^T\big[I(f_t)+K_\gamma(f_t^{\otimes2})\big]\md t<\infty.
 \label{vs:eq:IK}
\end{equation}
Then $f$ agrees with the comparison solution $g$:
\begin{equation}
 f_t=g_t\quad(0\le t\le T),
 \qquad f\in L^1(0,T;L^\infty(\R^3)).
 \label{vs:eq:identification}
\end{equation}
\end{proposition}

The proof is deferred to Section~\ref{sec:very-soft-closure}, where moment and covariance estimates yield the regularity needed for uniqueness on successive short time intervals.

Proposition~\ref{vs:thm:closure} identifies the limits for $\gamma\in[-3,-2)$; in the remaining ranges, we apply the known uniqueness theorems.
We then use energy conservation and the Fisher bounds to obtain the two convergence statements of Theorem~\ref{thm:main}.

For the soft-potential uniqueness estimates, define the singular convolution quantity
\begin{equation}
 J_\alpha(h):=\sup_{x\in\R^3}\int_{\R^3}|x-y|^\alpha h(\md y),
 \qquad \alpha\in(-3,0),\quad h\in\Pp(\R^3).
 \label{eq:singular-convolution}
\end{equation}

\begin{proposition}
\label{prop:unique}
Every subsequential limit of the empirical path laws is concentrated on the same Landau solution $f$ with initial datum $f_0$.
\end{proposition}

\begin{proof}
We verify the uniqueness conditions in Remark~\ref{rem:uniqueness-classes} for almost every limiting path.
For $\gamma\in[-2,0)$, \eqref{eq:fisher-hardy} gives, for every probability density $h$,
\[
 J_\gamma(h)\le1+J_{-2}(h)\le1+I(h).
\]
Apply Propositions~\ref{prop:identify} and~\ref{disc:thm:level3-diss} together with this inequality.
Almost every limiting path $\rho$ then belongs to $L^\infty(0,T;\Pp_2)\cap L^1(0,T;J_\gamma)$.
Lemma~\ref{vs:prop:comparison} gives the same bounds for the comparison solution $g$ and ensures that $g_t$ has a density at every time.
Applying \eqref{vs:eq:ordinary-coefficient-bound} to $\rho$ and $g$ converts their symmetric weak identities into the unsymmetrized weak formulation.
The stability theorem \cite[Theorem~1.3]{FournierGuerin2009JFA} gives $\rho=g$.

For $\gamma\in[-3,-2)$, almost every limiting path satisfies the hypotheses of Proposition~\ref{vs:thm:closure}, by Propositions~\ref{prop:identify} and~\ref{disc:thm:level3-diss}.
Hence these paths agree with the comparison solution $g$ and belong to $L^\infty(0,T;\Pp_2)\cap L^1(0,T;L^\infty)$.

At $\gamma=0$, finite-energy weak solutions are unique \cite{Villani1998MMMAS}; see also \cite[Section~1.6]{FournierHeydecker2021}.
For $\gamma\in(0,1]$, the moment bound in Proposition~\ref{prop:identify} makes the kernel $L_\varphi$ in \eqref{vs:eq:unsymmetrized-kernel} integrable and gives $\rho\in L^\infty(0,T;\Pp_2)\cap L^1(0,T;\Pp_{2+\gamma})$.
Together with energy conservation, this is the weak-solution class of \cite[Definition~1]{FournierHeydecker2021}.
Since $f_0\in\Pp_4$, apply the moment estimates of Desvillettes--Villani \cite{DesvillettesVillani2000I}, as recalled in \cite[Theorem~2(a)]{FournierHeydecker2021}.
They give
\[
 \sup_{s\le T}M_4(\rho_s)
 +\int_0^T M_{4+\gamma}(\rho_s)\md s<\infty.
\]
The corresponding existence result is recalled in Theorem~2(b) there, and \cite[Theorem~8]{FournierHeydecker2021}, with $p=4$, gives uniqueness.
This identifies all limiting paths in every interaction range.
\end{proof}

\begin{proof}[Proof of Theorem~\ref{thm:main}]
Fix $T>0$.
Propositions~\ref{prop:identify} and~\ref{prop:unique} show that every limit law equals $\delta_f$.
The solutions obtained on finite time intervals agree on overlaps by uniqueness and therefore define a single global path.
Tightness therefore gives convergence in law of the full sequence to this Dirac mass.
Convergence in law to a deterministic limit implies convergence in probability \cite[Section~3]{Billingsley1999}, hence
\[
 \mu^N\longrightarrow f
 \quad\hbox{in probability in }D([0,T],\Pp(\R^3)).
\]
Let $d_{\mathrm w}$ be a bounded complete metric generating the weak topology on $\Pp(\R^3)$, for instance the Prokhorov metric.
Since $f$ is continuous, the time-change characterization of $J_1$ convergence makes this equivalent to uniform weak convergence in time; see \cite[Section~12]{Billingsley1999} and \cite[Chapter~VI]{JacodShiryaev2003}:
\begin{equation}
 \sup_{t\le T}d_{\mathrm w}(\mu_t^N,f_t)\longrightarrow0
 \quad\hbox{in probability}.
 \label{eq:uniform-path-conv}
\end{equation}
Proposition~\ref{prop:identify} gives energy conservation for $f$.
Proposition~\ref{disc:thm:level3-diss}, applied with $\pi=\delta_f$, gives $I(f_t)\le I(f_0)$ at every time.

Next we use the control of the second moments. Since the MBC dynamics
conserves the energy pathwise, we have
\[
  \int_{\R^3}|v|^2\,\mu_t^{N}(\md v)
  = \int_{\R^3}|v|^2\,\mu_0^{N}(\md v)
  = \frac{1}{N}\sum_{i=1}^{N}|v_i^{N}(0)|^2
  \qquad\text{for all }t\ge0.
\]

By the law of large numbers applied at time $t=0$, this converges in probability to $\int_{\R^3}|v|^2f_0(v)\md v=E_0$.
From any subsequence, extract a further subsequence on which \eqref{eq:uniform-path-conv} and this convergence of second moments hold almost surely.
Combining the weak convergence $\mu_t^{N}\rightharpoonup f_t$
with the (pathwise) convergence of the second moments and using the classical
characterization of $W_2$-convergence in $\Pp_2(\R^3)$
(cf.~\cite[Proposition~7.1.5]{AmbrosioGigliSavare2008}), we deduce that
\[
  W_2(\mu_t^{N},f_t)\longrightarrow0
  \qquad\text{for every fixed }t\in[0,T]\ \text{almost surely}.
\]

Weak continuity and the constant second moment also imply $f\in C([0,T],\Pp_2)$ in the $W_2$ topology.
We now strengthen this to $\sup_{t\le T}W_2(\mu_t^N,f_t)\to0$ almost surely.
Indeed, otherwise one could choose $t_N\to t$ such that this distance stays positive.
But uniform weak convergence gives $\mu_{t_N}^N\rightharpoonup f_t$, and particle energy conservation gives
\[
 \int|v|^2\mu_{t_N}^N(\md v)=M_2^N(0)
 \longrightarrow E_0=\int|v|^2f_t(v)\md v.
\]
The same characterization \cite[Proposition~7.1.5]{AmbrosioGigliSavare2008} then gives $W_2(\mu_{t_N}^N,f_t)\to0$, contradicting the choice of $t_N$ and the $W_2$-continuity of $f$.
The subsequence criterion proves convergence in probability for the full sequence.
Using the standard estimate
\[
  W_2^2(\nu,\eta)
  \le
  2\int_{\R^3}|v|^2\,\nu(\md v)
  + 2\int_{\R^3}|v|^2\,\eta(\md v),
  \quad \forall\,\nu,\eta\in\Pp_2(\R^3),
\]

together with the energy bounds above, we get
\[
 \sup_{t\le T}W_2^2(\mu_t^N,f_t)
 \le2\big(M_2^N(0)+E_0\big),\qquad
 \E M_2^N(0)=E_0.
\]
Thus the family $\{\sup_{0\le t\le T}W_2(\mu_t^N,f_t)\}_{N\ge2}$ is uniformly integrable.
Combining this with convergence in probability proves the first assertion of Theorem~\ref{thm:main}.

Finally, we prove strong convergence of the marginals.
Fix $k\ge1$, $t\ge0$, and $\psi\in C_c^\infty(\R^{3k})$.
By exchangeability,
\[
 \big\langle F_t^{k|N},\psi\big\rangle
 =\E\Bigg[\frac1{(N)_k}
 \sum_{\text{distinct }i_1,\dots,i_k}
 \psi\big(v_{i_1}^N(t),\dots,v_{i_k}^N(t)\big)\Bigg],
\]
where $(N)_k:=N(N-1)\cdots(N-k+1)$.
Comparing sampling without/with replacement and using $N^k-(N)_k\le C_kN^{k-1}$, we obtain the standard replacement bound \eqref{disc:eq:empirical-marginal-comparison}.
Since $\nu\mapsto\langle\nu^{\otimes k},\psi\rangle$ is bounded and continuous for the weak topology, \eqref{eq:uniform-path-conv} gives
\[
 \E\big\langle(\mu_t^N)^{\otimes k},\psi\big\rangle
 \longrightarrow\langle f_t^{\otimes k},\psi\rangle.
\]
Passing to the limit in the replacement bound yields $\langle F_t^{k|N},\psi\rangle\to\langle f_t^{\otimes k},\psi\rangle$.
Finally, using the uniform second-moment bound (energy conservation), a standard truncation-approximation argument upgrades the convergence from $C_c^\infty$ to $C_b(\R^{3k})$.
Consequently,
\begin{equation}
 F_t^{k|N}\rightharpoonup f_t^{\otimes k}.
 \label{eq:chaos}
\end{equation}

Fisher contraction under marginalization \cite[Lemma~3.7(i)]{HaurayMischler2014} and Proposition~\ref{disc:thm:finite-N} give
\[
 I(F_t^{k|N})\le\frac kN I(F_t^N)\le kI(f_0).
\]
This implies $\{(F_t^{k|N})^{1/2}\}_{N\ge k}$ is bounded in $H^1(\R^{3k})$.
Hence, from any subsequence, we can extract a further subsequence along which $(F_t^{k|N})^{1/2}\to u$ weakly in $H^1(\R^{3k})$ and strongly in $L^2_{\loc}(\R^{3k})$, by the Rellich--Kondrachov compactness theorem and a diagonal argument.
For every bounded Borel set $B\subset\R^{3k}$, the Cauchy--Schwarz inequality gives
\[
 \|F_t^{k|N}-u^2\|_{L^1(B)}
 \le\|(F_t^{k|N})^{1/2}-u\|_{L^2(B)}
 \big(\|(F_t^{k|N})^{1/2}\|_{L^2(B)}+\|u\|_{L^2(B)}\big)
 \longrightarrow0.
\]
On the other hand, \eqref{eq:chaos} gives weak convergence to $f_t^{\otimes k}$.
The uniqueness of the limit yields $u^2=f_t^{\otimes k}$ almost everywhere, so $F_t^{k|N}\to f_t^{\otimes k}$ strongly in $L^1_{\loc}(\R^{3k})$.

To upgrade the local $L^1$ convergence to the global one, we use the second-moment bound.
For any $R>0$ and any $N\ge k$,
\[
 \int_{|v|>R}F_t^{k|N}(v)\md v
 \le\frac1{R^2}\int_{\R^{3k}}|v|^2F_t^{k|N}(v)\md v
 =\frac{kE_0}{R^2},
\]
and the same bound holds for $f_t^{\otimes k}$ by energy conservation.
Hence for any $R>0$,
\[
 \int_{|v|>R}\big(F_t^{k|N}(v)+f_t^{\otimes k}(v)\big)\md v
 \le\frac{2kE_0}{R^2}.
\]
Given $\varepsilon>0$, choose $R$ so large that $2kE_0/R^2<\varepsilon/2$.
On the ball $B_R$ we have $F_t^{k|N}\to f_t^{\otimes k}$ in $L^1(B_R)$, so for $N$ large enough along the further subsequence,
\[
 \int_{|v|\le R}|F_t^{k|N}(v)-f_t^{\otimes k}(v)|\md v<\frac\varepsilon2.
\]
Combining the estimates on $B_R$ and its complement, we obtain
\[
 \|F_t^{k|N}-f_t^{\otimes k}\|_{L^1(\R^{3k})}
 \le\int_{|v|\le R}|F_t^{k|N}-f_t^{\otimes k}|\md v
 +\int_{|v|>R}\big(F_t^{k|N}+f_t^{\otimes k}\big)\md v
 <\varepsilon
\]
for all sufficiently large $N$ along this further subsequence.
The subsequence criterion proves the second assertion of Theorem~\ref{thm:main} for the full sequence.
\end{proof}

\section{Regularity for very soft potentials}
\label{sec:very-soft-closure}\label{sec:very-soft-static}
\label{sec:very-soft-continuation}

This section proves Proposition~\ref{vs:thm:closure}.
Throughout, assume $\gamma\in[-3,-2)$, and let $g$ be the comparison solution from Lemma~\ref{vs:prop:comparison}.
We combine moment and covariance bounds with a weighted gradient estimate to obtain the regularity needed for uniqueness on short time intervals.
Agreement with $g$ at each endpoint restores the covariance lower bound and allows the argument to continue.
The proof is completed at the end of this section.

We use the convolution estimate \eqref{eq:fisher-hardy} both to propagate the moments $M_\ell$ defined in \eqref{eq:moment-definition} and to control the variation of the covariance.

\begin{lemma}
\label{vs:lem:moment}
Let $f\in C_w([0,T],\Pp(\R^3))$ satisfy \eqref{eq:weak} for every $C_b^2$ test and $\int_0^T I(f_t)\md t<\infty$.
If $\ell\ge2$ and $M_\ell(f_0)<\infty$, then
\begin{equation}
 \sup_{0\le t\le T}M_\ell(f_t)
 \le
 C_{\gamma,\ell}
 \big(1+M_\ell(f_0)\big)
 \exp\bigg[
 C_{\gamma,\ell}
 \int_0^T\big(1+I(f_s)\big)\md s
 \bigg].
 \label{vs:eq:moment}
\end{equation}
\end{lemma}

\begin{proof}
We adapt the strategy of \cite[proof of Lemma~8]{CarrapatosoDesvillettesHe2017}, replacing the entropy dissipation in the near-diagonal estimate by the Fisher-information bound.
For the cutoff argument, set
\[
 \Phi_R(v)=\langle v\rangle^\ell
 \Big(1+\frac{|v|^2}{R^2}\Big)^{-\ell/2},\qquad R\ge2.
\]
They belong to $C_b^2$ and increase to $\langle v\rangle^\ell$ as $R\to\infty$.
As in that argument, we have $|\nabla\Phi_R|+|D^2\Phi_R|\le C_\ell\Phi_R$, also for $\ell=2$.
Hence $\Phi_R$ is uniformly comparable on unit balls, with constants independent of $R$.
Applying the mean value theorem to the drift term in \eqref{eq:Hphi} when $|v-w|<1$, and using $|v-w|^{\gamma+2},|v-w|^{\gamma+1}\le1$ when $|v-w|\ge1$, we obtain
\[
 |H_{\Phi_R}(v,w)|
 \le C_{\gamma,\ell}(1+|v-w|^{-2})
 [\Phi_R(v)+\Phi_R(w)].
\]
By symmetry and \eqref{eq:fisher-hardy}, the weak equation gives
\[
 \langle f_t,\Phi_R\rangle
 \le\langle f_0,\Phi_R\rangle
 +C_{\gamma,\ell}\int_0^t
 (1+I(f_s))\langle f_s,\Phi_R\rangle\md s.
\]
An application of Gr\"onwall's inequality followed by the monotone convergence theorem yields \eqref{vs:eq:moment}.
\end{proof}

\begin{lemma}
\label{vs:lem:weak-modulus}
Under the hypotheses of Lemma~\ref{vs:lem:moment} with $\ell=2$, the mean of $f_t$ is constant.
Define its mean and covariance by
\begin{equation}
 \bar v:=\int v f_0(\md v),\qquad
 \Sigma_f(t):=\int(v-\bar v)\otimes(v-\bar v)f_t(\md v).
 \label{vs:eq:mean-covariance}
\end{equation}
The covariance satisfies
\begin{equation}
 \|\Sigma_f(t)-\Sigma_f(s)\|_2
 \le C_\gamma\int_s^t(1+I(f_u))\md u
 \qquad(0\le s\le t\le T).
 \label{eq:weak-modulus}
\end{equation}
\end{lemma}

\begin{proof}
Lemma~\ref{vs:lem:moment} gives bounded second moments.
Fix a smooth compactly supported cutoff $\chi$, equal to one near the origin.
Use the bounded tests
\[
 (v_i-\bar v_i)\chi\big((v-\bar v)/R\big),\qquad
 \big(\xi\cdot(v-\bar v)\big)^2\chi\big((v-\bar v)/R\big),
\]
for unit vectors $\xi$.
Their Hessians are bounded by $C/R$ and $C$, respectively, uniformly in $\xi$ and $R\ge1$.
By the mean value theorem and \eqref{eq:Hphi}, the corresponding kernels are bounded by $C_\gamma R^{-1}|v-w|^{\gamma+2}$ and $C_\gamma|v-w|^{\gamma+2}$.
Since $|z|^{\gamma+2}\le1+|z|^{-2}$, we can pass to the linear and quadratic tests by the dominated convergence theorem, using \eqref{eq:fisher-hardy} and the second-moment bound.
The linear kernel vanishes in the limit, proving conservation of the mean.
For the quadratic test $(\xi\cdot(v-\bar v))^2$, formula \eqref{eq:Hphi} gives the identity
\[
 \xi^\top[\Sigma_f(t)-\Sigma_f(s)]\xi
 =\int_s^t\iint
 \big[|v-w|^{\gamma+2}
 -3|v-w|^\gamma(\xi\cdot(v-w))^2\big]
 f_u(\md v)f_u(\md w)\md u.
\]
The bracket has absolute value at most $2|v-w|^{\gamma+2}$ for every unit $\xi$.
Applying \eqref{eq:fisher-hardy} and taking the supremum over these vectors gives \eqref{eq:weak-modulus}.
\end{proof}

We next recover an $L^\infty$ bound from the directional information under a positive lower bound on the covariance matrix.
For the two-particle estimates, we use the fields $X_k$ from Remark~\ref{disc:rem:singular-affinity}, with coordinates $(v,w)$ and $r=|v-w|$.

\begin{proposition}
\label{vs:prop:static}
Let $m>3(2-\gamma)$, and fix $\kappa>0$ and $M<\infty$.
Let $\rho$ be a probability density whose covariance matrix has all eigenvalues greater than or equal to $\kappa$, and suppose that
\[
 M_m(\rho)\le M,\qquad I(\rho)+K_\gamma(\rho^{\otimes2})<\infty.
\]
Then
\begin{equation}
 \|\rho\|_{L^\infty}
 \le C_{\gamma,m,\kappa,M}
 \big[1+K_\gamma(\rho^{\otimes2})^{3/4}\big].
 \label{vs:eq:static}
\end{equation}
The functional $K_\gamma$ is understood in the extended sense of \eqref{disc:eq:Kgamma}.
\end{proposition}

The proof rests on the moment identity \eqref{stat:eq:moment-reconstruction}: it reconstructs each log-gradient component from tangential differences when the covariance matrix has a positive lower bound.
Combining this identity with Lemma~\ref{stat:lem:active-difference} gives weighted fourth-power control of the full log-gradient; the local form of the Gagliardo--Nirenberg inequality then yields the stated $L^\infty$ bound.
The proof is given below, after Lemma~\ref{stat:lem:active-difference}.

For a probability density $\rho$ with finite Fisher information, we write $\nabla\log \rho=\nabla \rho/\rho$ on $\{\rho>0\}$ and set it to zero on $\{\rho=0\}$.
To measure fourth-power tangential differences of the log-gradient, define
\begin{equation}
 \mathcal D_\gamma(\rho)
 :=\iint \rho(v)\rho(w)\langle v-w\rangle^{\gamma-2}
 \big|(v-w)\times\big(\nabla\log \rho(v)-\nabla\log \rho(w)\big)\big|^4\md v\md w.
 \label{stat:eq:tangential-information}
\end{equation}

The following estimate is a version of Tabary's fourth-power directional bound \cite[Proposition~2.1 and Lemma~8.3]{Tabary2025_PoC_Landau_Coulomb}.
The proof uses the same integration-by-parts identity; angular regularization justifies it under the stated finite-information assumptions.

\begin{lemma}
\label{stat:lem:active-difference}
If $I(\rho)+K_\gamma(\rho^{\otimes2})<\infty$, then
\begin{equation}
 \mathcal D_\gamma(\rho)
 \le\frac{27}{4}K_\gamma(\rho^{\otimes2}).
 \label{stat:eq:active-difference}
\end{equation}
\end{lemma}

The angular regularization argument is given in Appendix~\ref{sec:auxiliary-information}.

\begin{proof}[Proof of Proposition~\ref{vs:prop:static}]
Translate by the mean $\bar v=\int v\rho(v)\md v$ and rotate to diagonalize the covariance, retaining the notation $\rho$ in these coordinates.
The quantity $\mathcal D_\gamma(\rho)$ defined in \eqref{stat:eq:tangential-information} is unchanged under this transformation.
Since $|\bar v|\le M^{1/m}$, the transformed moments are bounded in terms of $m,M$, and the weights $\langle v\rangle^{\gamma-2}$ before and after translation are comparable with constants depending only on $\gamma,m,M$.
Write
\[
 \int v_i\rho(v)\md v=0,\qquad
 \int v_iv_j\rho(v)\md v=a_i\delta_{ij},
 \qquad a_1,a_2,a_3\ge\kappa.
\]
Finite Fisher information gives $\rho\in W^{1,1}$ and $\nabla \rho=\rho\nabla\log \rho$, including on the vacuum set.
Moreover, $\int\langle w\rangle^2|\nabla \rho(w)|\md w\le M_4(\rho)^{1/2}I(\rho)^{1/2}<\infty$, so integration by parts against polynomials of degree at most two is justified by cutoffs.

For $i\ne j$, set
\[
\begin{aligned}
 Q_{ij}(v,w)
 &:= (e_i\times e_j)\cdot
 \big[(v-w)\times\big(\nabla\log \rho(v)-\nabla\log \rho(w)\big)\big]\\
 &=(v_i-w_i)\big(\partial_j\log \rho(v)-\partial_j\log \rho(w)\big)
 -(v_j-w_j)\big(\partial_i\log \rho(v)-\partial_i\log \rho(w)\big).
\end{aligned}
\]
The key step is the following moment identity:
\begin{equation}
 \partial_j\log \rho(v) 
 = - \frac{1}{a_i} \int w_i\rho(w)Q_{ij}(v,w)\md w - \frac{v_j}{a_i}.
\label{stat:eq:moment-reconstruction}
\end{equation}
Indeed, by the integration by parts we have
\[
\begin{gathered}
 \int w_i\rho(w)Q_{ij}(v,w)\md w
 =-a_i\partial_j\log \rho(v)-\int w_i(v_i-w_i)\partial_j\rho(w)\md w\\
 +\int w_i(v_j-w_j)\partial_i\rho(w)\md w=-a_i\partial_j\log \rho(v)-v_j.
\end{gathered}
\]
The bound $\langle v-w\rangle\le\sqrt2\langle v\rangle\langle w\rangle$ gives
\[
 \int \rho(w)|w_i|^{4/3}\langle v-w\rangle^{(2-\gamma)/3}\md w
 \le C_\gamma\langle v\rangle^{(2-\gamma)/3}M_{(6-\gamma)/3}(\rho).
\]
Since $(6-\gamma)/3<m$, weighted H\"older's inequality therefore yields
\[
\begin{split}
 \Big|\int w_i\rho(w)Q_{ij}(v,w)\md w\Big|^4
 &\le
 \Big(\int \rho(w)\langle v-w\rangle^{\gamma-2}|Q_{ij}(v,w)|^4\md w\Big)\\
 &\quad\cdot
 \Big(\int \rho(w)|w_i|^{4/3}
 \langle v-w\rangle^{(2-\gamma)/3}\md w\Big)^3\\
 &\le C_{\gamma,m,M}\langle v\rangle^{2-\gamma}
 \int \rho(w)\langle v-w\rangle^{\gamma-2}|Q_{ij}(v,w)|^4\md w.
\end{split}
\]
For each $j$, choose any $i\ne j$.
Use \eqref{stat:eq:moment-reconstruction} with $a_i\ge\kappa$, multiply by $\rho(v)\langle v\rangle^{\gamma-2}$, and integrate in $v$.
Since $|Q_{ij}(v,w)|\le |(v-w)\times (\nabla\log \rho(v)-\nabla\log \rho(w))|$ and $\langle v\rangle^{\gamma-2}|v_j|^4\le1$ for $\gamma\le-2$, summation over $j$ gives
\[
 A:=\int \rho(v)\langle v\rangle^{\gamma-2}|\nabla\log \rho(v)|^4\md v
 \le C_{\gamma,m,\kappa,M}\big(1+\mathcal D_\gamma(\rho)\big).
\]
The weight comparability under translation gives the same bound for $A$ in the original coordinates, to which we now return.
Lemma~\ref{stat:lem:active-difference} therefore yields
\[
 A\le C_{\gamma,m,\kappa,M}\big(1+K_\gamma(\rho^{\otimes2})\big).
\]
Setting $u=\rho^{1/4}$, we have
\[
 |\nabla u|^4 = 4^{-4}\rho|\nabla\log \rho|^4.
\]
Apply the Gagliardo--Nirenberg inequality with a cutoff from $B_2(x)$ to $B_1(x)$ to obtain
\[
 \|u\|_{L^\infty(B_1(x))}^4
 \le C\bigg[
 \Big(\int_{B_2(x)}|\nabla u|^4\Big)^{3/4}
 \Big(\int_{B_2(x)}u^4\Big)^{1/4}
 +\int_{B_2(x)}u^4
 \bigg].
\]
On $B_2(x)$, $\langle v\rangle$ and $\langle x\rangle$ are comparable.
Hence
\[
 \int_{B_2(x)}|\nabla u|^4
 \le C\langle x\rangle^{2-\gamma}A,
 \qquad
 \int_{B_2(x)}u^4
 \le C M\langle x\rangle^{-m}.
\]
It follows that
\[
 \|\rho\|_{L^\infty(B_1(x))}
 \le C\Big[
 A^{3/4}M^{1/4}
 \langle x\rangle^{[3(2-\gamma)-m]/4}
 +M\langle x\rangle^{-m}
 \Big].
\]
The powers of $\langle x\rangle$ are nonpositive, so taking the supremum in $x$ proves \eqref{vs:eq:static}.
\end{proof}

\begin{proof}[Proof of Proposition~\ref{vs:thm:closure}]
By Lemma~\ref{vs:lem:moment}, we have
\begin{equation}
 \sup_{t\le T}M_m(f_t)<\infty.
 \label{vs:eq:f-moment}
\end{equation}
By Lemma~\ref{vs:lem:weak-modulus}, $f$ and $g$ have the same constant mean $\bar v=\int vf_0(v)\md v$.
Set
\[
 \kappa:=\frac1{2(1+I(f_0))}.
\]
Recall the covariance notation in \eqref{vs:eq:mean-covariance}.
For every unit vector $\xi$, integration by parts and the Cauchy--Schwarz inequality give
\[
 1=-\int\xi\cdot(v-\bar v)\,\xi\cdot\nabla g_t(v)\md v
 \le\big[\xi^\top\Sigma_g(t)\xi\,I(g_t)\big]^{1/2}.
\]
Finite energy and Fisher information justify the integration by parts.
By Lemma~\ref{vs:prop:comparison}, $I(g_t)\le I(f_0)$.
The preceding estimate therefore gives
\[\Sigma_g(t)\ge2\kappa\operatorname{Id}\quad \forall\,t\in[0,T].\]
This covariance bound is a consequence of the Cram\'er--Rao inequality for translation families.

Next, by \eqref{vs:eq:IK}, we have $\int_0^T(1+I(f_r))\md r<\infty$.
By absolute continuity of this integral, choose a finite partition $0=t_0<\cdots<t_M=T$ such that
\[
 \int_{t_k}^{t_{k+1}}(1+I(f_r))\md r\le\kappa/C_\gamma
 \quad(0\le k<M),
\]
where $C_\gamma$ is the constant in Lemma~\ref{vs:lem:weak-modulus}.
We prove by induction that $f=g$ on each interval.
At $t_0=0$, the initial data agree.
Suppose that $f_{t_k}=g_{t_k}$, so $\Sigma_f(t_k)=\Sigma_g(t_k)$.
By Lemma~\ref{vs:lem:weak-modulus}, for every $t\in[t_k,t_{k+1}]$,
\[
 \Sigma_f(t)\ge\Sigma_g(t_k)
 -C_\gamma\int_{t_k}^t(1+I(f_r))\md r\,\operatorname{Id}
 \ge\kappa\operatorname{Id}.
\]
Thus the covariance eigenvalues of $f_t$ are bounded below by $\kappa$ throughout this interval.
The moment bound \eqref{vs:eq:f-moment} and the assumption \eqref{vs:eq:IK} verify the remaining hypotheses of Proposition~\ref{vs:prop:static} at almost every time.
Applying this proposition and H\"older's inequality, we obtain
\[
 \int_{t_k}^{t_{k+1}}\|f_t\|_\infty\md t
 \le C_T\bigg[(t_{k+1}-t_k)
 +(t_{k+1}-t_k)^{1/4}
 \Big(\int_{t_k}^{t_{k+1}}
 K_\gamma(f_t^{\otimes2})\md t\Big)^{3/4}\bigg]<\infty.
\]

For $f$, Lemma~\ref{vs:lem:moment} gives bounded energy, and \eqref{vs:eq:IK} gives finite time-integrated Fisher information.
For $g$, the same bounds follow from Lemma~\ref{vs:prop:comparison}.
Hence \eqref{vs:eq:ordinary-coefficient-bound} applies to both paths, and their weak identities also hold in the unsymmetrized form.
Subtract the identity at $t_k$ from the identity at $t$.
On $[t_k,t_{k+1}]$, both paths therefore solve the same weak equation with the same initial datum.
Moreover, Lemma~\ref{vs:prop:comparison} gives $g\in L^1(t_k,t_{k+1};L^\infty)$.

If $\gamma\in(-3,-2)$, then, for $h=f,g$, the singular convolutions \eqref{eq:singular-convolution} satisfy
\[
 J_\gamma(h_t)+J_{\gamma+1}(h_t)
 \le2+C_\gamma\|h_t\|_\infty.
\]
Both solutions have bounded energy, and $g_t$ has a density at every time by Lemma~\ref{vs:prop:comparison}.
The Fournier--Gu\'erin stability theorem \cite[Theorem~1.3]{FournierGuerin2009JFA} therefore gives $f=g$ on $[t_k,t_{k+1}]$.
For $\gamma=-3$, apply \cite[Theorem~2]{Fournier2010CMP} to obtain the same conclusion.
In particular, $f_{t_{k+1}}=g_{t_{k+1}}$, so the induction continues.
It follows that $f=g$ on $[0,T]$.
Finally, Lemma~\ref{vs:prop:comparison} gives \eqref{vs:eq:identification}.
\end{proof}

\section{Discussion: Grazing-collision limit of the MBC process}\label{sec:grazing-discussion}

The MBC dynamics depends on two parameters: the particle number $N$
and the grazing-collision parameter $h$.
The Landau equation corresponds to the joint
mean-field and grazing-collision regime, where both $N\to\infty$
and $h\to0$.
Classically, one follows a two-step route: starting from a binary
collision model or the Boltzmann equation, a grazing-collision limit
at fixed $N$ leads to an $N$-particle Landau diffusion (the
grazing-collision master equation), and propagation of chaos as
$N\to\infty$ then yields the Landau equation.
Our main result, Theorem~\ref{thm:main}, shows that the Landau flow
can be obtained directly along the diagonal limit
$N\to\infty$, $h\to0$ for the MBC process, without explicitly
passing through this intermediate level.

It is natural to compare this diagonal limit with the ``upper edge''
of the classical picture, where one first lets $h\to0$ at fixed $N$
and only then sends $N\to\infty$.
In this section we rigorously justify this upper edge:
for each fixed $N$, the time-marginal laws of the MBC process converge, as $h\to0$, to those of the
standard $N$-particle Landau diffusion.
The complementary ``lower edge'' (first $N\to\infty$ at fixed $h$
and then $h\to0$ for the resulting nonlinear dynamics) may involve
some distinct ideas and is left for future work.

\medskip
\begin{center}
	{\small
		\begin{tikzpicture}[node distance=5cm,>=stealth]
			\node (A) at (0,3) {%
				\begin{tabular}{c}
					{MBC Process} \\
					$(N, h)$
				\end{tabular}
			};
			\node (B) at (6.5,3) {%
				\begin{tabular}{c}
					{$N$-particle Landau Diffusion} \\
					$(N, 0)$
				\end{tabular}
			};
			\node (C) at (0,0) {%
				\begin{tabular}{c}
					{MV Jump Process} \\
					$(\infty, h)$
				\end{tabular}
			};
			\node (D) at (6.5,0) {%
				\begin{tabular}{c}
					{Landau Equation} \\
					$(\infty, 0)$
				\end{tabular}
			};
			\draw[->] (A) -- node[above] {$h \to 0$} node[sloped, below] {Proposition~\ref{thm:fixed-N-grazing}} (B);
			\draw[->] (C) -- node[below] {$h \to 0$} (D);
			\draw[->] (A) -- node[left] {$N \to \infty$} (C);
			\draw[->] (B) -- node[right] {$N \to \infty$} (D);
			\draw[->] (A) -- node[sloped, above] {$N \to \infty$, $h\to 0$} node[sloped, below] {Theorem~\ref{thm:main}}(D);
		\end{tikzpicture}
	}
\end{center}
\vspace{0.8ex}

We recall the $N$-particle Landau diffusion operator which describes the
grazing-collision master equation at fixed $N$.
For $\Phi\in C_c^4(\R^{3N})$, the limiting differential expression is
\begin{equation}
\label{eq:LN-def}
\begin{aligned}
\mathcal{L}^{N}\Phi(\bm v)
&:=\frac{2}{N-1}\sum_{i<j}A_{ij}\Phi(\bm v)\\
&=\frac{1}{4(N-1)}
\sum_{i\neq j}
a(v_{i}-v_{j}):\big(\nabla_{v_{i}}^{2}\Phi(\bm v)
+\nabla_{v_{j}}^{2}\Phi(\bm v)
-\nabla_{v_i}\nabla_{v_j}\Phi(\bm v)
-\nabla_{v_j}\nabla_{v_i}\Phi(\bm v)\big)\\
&\quad
+\frac{1}{2(N-1)}
\sum_{i\neq j}
b(v_{i}-v_{j})\cdot
\big(\nabla_{v_{i}}\Phi(\bm v)-\nabla_{v_{j}}\Phi(\bm v)\big),
\end{aligned}
\end{equation}
with $a$ given by \eqref{eq:A-kernel} and $b(z)=-2|z|^{\gamma}z$.
Here $A_{ij}$ is the pair operator in \eqref{eq:MBC-definition}; the self-adjoint realization of $\mathcal L^N$ is specified in the proof below.

\begin{proposition}
	\label{thm:fixed-N-grazing}
	Let $N\ge2$ and $\gamma\in[-3,1]$.
	Assume that the initial law $F_{0}^{N}\in\Pp(\R^{3N})\cap L^{\infty}(\R^{3N})$
	is exchangeable and has finite second moment and finite Fisher information.
	For each $h>0$, let $(F_{t}^{N,h})_{t\ge0}$ be the law of the
	$N$-particle MBC process started from $F_{0}^{N}$.
	
	Then there exists a weakly continuous curve
	$(F_{t}^{N})_{t\ge0}\subset\Pp(\R^{3N})$ such that, for every
	$T>0$ and every $\Phi\in C_{c}^{4}(\R^{3N})$,
	\begin{equation}
		\label{eq:fixedN-conv}
		\lim_{h\to0}\ \sup_{t\in[0,T]}
		\big|\langle F_{t}^{N,h}-F_{t}^{N},\Phi\rangle\big|=0.
	\end{equation}
	
	Moreover, the limit $(F_{t}^{N})_{t\ge0}$ satisfies the $N$-particle Landau master equation
	\[
	\partial_{t}F_{t}^{N}=(\mathcal{L}^{N})^{*}F_{t}^{N} 
	\]
	in the weak sense, and it inherits the same a priori bounds:
	for all $t\ge0$,
	\begin{equation}\label{eq:Fn-bound}
		\|F_{t}^{N}\|_{\infty}\le\|F_{0}^{N}\|_{\infty},\quad
		I(F_{t}^{N})\le I(F_{0}^{N}),\quad
		\int_{\R^{3N}}|\bm v|^{2}F_{t}^{N}(\md\bm v)
		=\int_{\R^{3N}}|\bm v|^{2}F_{0}^{N}(\md\bm v),
	\end{equation}
	so total mass, momentum and energy are conserved and the limiting
	Landau dynamics remains in the same regularity class as the MBC
	scheme.
\end{proposition}

\begin{remark}
	Proposition~\ref{thm:fixed-N-grazing} gives a fixed-$N$ grazing-collision
	limit for the MBC scheme:
	it shows that the associated time-marginal laws converge to those of the classical
	$N$-particle Landau diffusion (grazing-collision master equation),
	while preserving uniform $L^\infty$, Fisher information, and energy--momentum
	bounds.
	This provides a concrete, structure-preserving microscopic realization
	of the Landau master equation, complementing earlier studies
	of the Landau master equation (cf.~\cite{KiesslingLancellotti2004}), constructions based on Boltzmann-to-Landau grazing limits (cf.~\cite{MiotPulvirentiSaffirio2011, Carrapatoso2016KRM}), and diffusive particle systems (cf. \cite{FontbonaGuerinMeleardPTRF09, Fournier2009KRM}).
	In particular, it shows that the intermediate $N$-particle Landau
	structure arises directly from the MBC process, consistently with the
	diagonal limit $N\to\infty$, $h\to0$ in Theorem~\ref{thm:main}.
\end{remark}

\begin{proof}[Proof of Proposition~\ref{thm:fixed-N-grazing}]
	Fix $N\ge2$ and $\gamma\in[-3,1]$.
	For each $h>0$, let $(F_{t}^{N,h})_{t\ge0}$ denote the law of the
	$N$-particle MBC process with generator
	$\mathcal{L}^{N,h}=\frac{N}{h}(\mathcal{Q}^{N,h}-I)$ and initial data
	$F_{0}^{N}$.
	
	\emph{Step 1: Uniform a priori bounds.}
	
	By Propositions~\ref{prop:a.c.bound} and~\ref{disc:thm:finite-N}, and
	conservation of momentum and energy for the MBC dynamics,
	for all $t\ge0$ and all $h>0$,
	\[
	\|F_{t}^{N,h}\|_{\infty}\le\|F_{0}^{N}\|_{\infty},\quad
	I(F_{t}^{N,h})\le I(F_{0}^{N}),\quad
	\int_{\R^{3N}}|\bm v|^{2}F_{t}^{N,h}(\md\bm v)
	=\int_{\R^{3N}}|\bm v|^{2}F_{0}^{N}(\md\bm v).
	\]
	Hence $\{F_{t}^{N,h}:t\in[0,T],\,0<h\le1\}$ is tight in
	$\Pp(\R^{3N})$ for every $T>0$.
	
	The density bounds also give
\[
 \sup_{h>0,\,t\ge0}\|F_t^{N,h}\|_2
 \le\|F_0^N\|_\infty^{1/2}.
\]

	\emph{Step 2: Equicontinuity of observables.}

Let $\Phi\in C_{c}^{4}(\R^{3N})$ and set
	$m_{\Phi}^{h}(t):=\langle F_{t}^{N,h},\Phi\rangle$.
	The laws $F_t^{N,h}$ are exchangeable, and the generators commute with particle permutations.
All pairings below are therefore unchanged when $\Phi$ is replaced by its average over particle permutations, so we may assume that $\Phi$ is symmetric.
Since $F_{t}^{N,h}$ solves the forward Kolmogorov equation associated
	with $\mathcal{L}^{N,h}$, we have
	\[
	\frac{\md}{\md t}m_{\Phi}^{h}(t)
	=\langle F_{t}^{N,h},\mathcal{L}^{N,h}\Phi\rangle.
	\]
	Using the pairwise expansion of the spherical heat semigroup, we write
	\[
	\mathcal{L}^{N,h}\Phi
	=\mathcal{L}^{N}\Phi+\mathcal{R}^{N,h}\Phi,
	\]
	where $\mathcal{L}^{N}$ is given by \eqref{eq:LN-def} and
	$\mathcal{R}^{N,h}\Phi$ is a remainder of the same structure as in
	Lemma~\ref{lem:remainder}.
	
	We record the bounds needed for this expansion on the full particle space.
Choose $R_\Phi$ with $\supp\Phi\subset\{\bm v:|\bm v|<R_\Phi\}$.
Pair rotations preserve $|\bm v|$, so every term in the expansion vanishes outside this ball.
For fixed remaining coordinates and pair centre, symmetry makes the function of the relative velocity even.
The mean value theorem and the spherical differentiation used in Lemma~\ref{lem:remainder} then give
\[
 |\Delta_{\S^2}\Phi|+|\Delta_{\S^2}^2\Phi|
 \le C_\Phi r_{ij}^2.
\]
The second-order remainder estimate for the heat semigroup consequently gives
\[
 |\mathcal R^{N,h}\Phi(\bm v)|
 \le C_{N,\Phi}\sum_{i<j}
 \min\big\{h r_{ij}^{2\gamma+2},r_{ij}^{\gamma+2}\big\}
 \mathbf1_{\{|\bm v|<R_\Phi\}}.
\]
On $r_{ij}<1$, each summand is bounded by $C\min\{h r_{ij}^{-4},r_{ij}^{-1}\}$, uniformly for $\gamma\in[-3,1]$.
On $r_{ij}\ge1$, it is bounded by $C_\Phi h$.
Integrating the square in relative polar coordinates and splitting at $r=h^{1/3}$ gives
\[
 \|\mathcal R^{N,h}\Phi\|_2^2
 \le C_{N,\Phi}
 \bigg[\int_0^{h^{1/3}}\md r
 +h^2\int_{h^{1/3}}^1r^{-6}\md r+h^2\bigg]
 \le C_{N,\Phi}h^{1/3}.
\]
Together with the uniform $L^2$ bound in Step~1, the Cauchy--Schwarz inequality yields
\begin{equation}
 \big|\langle F_t^{N,h},\mathcal R^{N,h}\Phi\rangle\big|
 \le C_\Phi h^{1/6}
 \qquad(t\ge0,\ 0<h\le1).
 \label{eq:fixedN-rem}
\end{equation}
The same differentiation gives
\[
 |\mathcal L^N\Phi(\bm v)|
 \le C_{N,\Phi}\sum_{i<j}(1+r_{ij}^{-1})
 \mathbf1_{\{|\bm v|<R_\Phi\}},
\]
so $\mathcal L^N\Phi\in L^2(\R^{3N})$ and, for each $T>0$,
\begin{equation}
 \sup_{0<h\le1}\sup_{t\in[0,T]}
 \big|\langle F_t^{N,h},\mathcal L^N\Phi\rangle\big|
 \le C_{\Phi,T}<\infty.
 \label{eq:fixedN-LN-bound}
\end{equation}
Combining these,
	\[
	|m_{\Phi}^{h}(t)-m_{\Phi}^{h}(s)|
	\le\int_{s}^{t}\big|\langle F_{r}^{N,h},\mathcal{L}^{N}\Phi\rangle\big|\md r
	+\int_{s}^{t}\big|\langle F_{r}^{N,h},\mathcal{R}^{N,h}\Phi\rangle\big|\md r
	\le(C_{\Phi,T}+C_{\Phi})|t-s|
	\]
	for all $0\le s\le t\le T$ and all $0<h\le1$.
	Thus $\{m_{\Phi}^{h}\}_{0<h\le1}$ is equi-Lipschitz on $[0,T]$.
	
	By exchangeability, the same conclusion holds for every $\Phi\in C_c^4$.

	\emph{Step 3: Extraction of a limit curve.}

Let $\{\Phi_{k}\}_{k\ge1}\subset C_{c}^{4}(\R^{3N})$ be a countable
	dense family in $C_0(\R^{3N})$ separating probability measures.
	By the equicontinuity above and tightness at each fixed time,
	the Arzel\`a--Ascoli theorem plus a diagonal argument yield a subsequence
	$h_{n}\downarrow0$ and a weakly continuous curve
	$(F_{t}^{N})_{t\ge0}\subset\Pp(\R^{3N})$ such that, for every
	$T>0$ and every $\Phi\in C_{c}^{4}(\R^{3N})$,
	\[
	\sup_{t\in[0,T]}
	\big|\langle F_{t}^{N,h_{n}}-F_{t}^{N},\Phi\rangle\big|\to0.
	\]
	The limit remains exchangeable, and lower semicontinuity gives the density and Fisher-information bounds in \eqref{eq:Fn-bound}.
For every bounded continuous function $\psi$ on $[0,\infty)$, pathwise energy conservation gives
\[
 \int\psi(|\bm v|^2)F_t^{N,h}(\md\bm v)
 =\int\psi(|\bm v|^2)F_0^N(\md\bm v).
\]
This identity passes to the weak limit.
Taking $\psi(s)=\min\{s,R\}$ and then letting $R\to\infty$ proves the energy equality in \eqref{eq:Fn-bound}.
The same argument, applied to bounded cutoffs of the conserved total momentum, proves momentum conservation.

	\emph{Step 4: Identification as Landau master equation.}

For each symmetric $\Phi\in C_{c}^{4}(\R^{3N})$ and $t\ge0$ we have
	\[
	\langle F_{t}^{N,h},\Phi\rangle
	=\langle F_{0}^{N},\Phi\rangle
	+\int_{0}^{t}\langle F_{s}^{N,h},\mathcal{L}^{N}\Phi\rangle\md s
	+\int_{0}^{t}\langle F_{s}^{N,h},\mathcal{R}^{N,h}\Phi\rangle\md s.
	\]
	The weak convergence from Step~3 and the uniform $L^2$ bound imply weak $L^2$ convergence of the densities at each fixed time.
Since $\mathcal L^N\Phi\in L^2$, its pairing therefore converges, and \eqref{eq:fixedN-LN-bound} permits the dominated convergence theorem in time.
Passing to the limit along $h_{n}\to0$, using
	\eqref{eq:fixedN-rem}, \eqref{eq:fixedN-LN-bound} and the convergence
	of $F_{t}^{N,h_{n}}$ against $\Phi$ and $\mathcal{L}^{N}\Phi$, we obtain
	for all $t\ge0$,
	\[
	\langle F_{t}^{N},\Phi\rangle
	=\langle F_{0}^{N},\Phi\rangle
	+\int_{0}^{t}\langle F_{s}^{N},\mathcal{L}^{N}\Phi\rangle\md s.
	\]
	Permutation symmetry extends this identity to every $\Phi\in C_c^4$.
For such tests, the collision term is absolutely integrable: without symmetry, its singularity is at most $C_\Phi r_{ij}^{-2}$ on a bounded set, which is locally integrable in the three-dimensional relative variable and can be paired with the bounded density $F_s^N$.
Thus $(F_{t}^{N})_{t\ge0}$ is a weak solution of the
	$N$-particle Landau master equation
	$\partial_{t}F_{t}^{N}=(\mathcal{L}^{N})^{*}F_{t}^{N}$ with
	$F_{0}^{N}$.
	
	\emph{Step 5: Uniqueness of the limit and full convergence.}

We use the self-adjoint pair flows to identify a single limiting semigroup.
For $u\in L^2(\R^{3N})$, set
\[
 \mathcal E_h(u)
 :=-\langle u,\mathcal L^{N,h}u\rangle
 =\frac{2}{N-1}\sum_{i<j}
 \frac{\langle u,(I-\ee^{hA_{ij}})u\rangle}{h}.
\]
Each $-A_{ij}$ is nonnegative and self-adjoint, by the spherical heat representation and the whole-space reversibility in Corollary~\ref{cor:whole-space}.
The spectral theorem shows that $\mathcal E_h$ increases as $h\downarrow0$ to
\[
 \mathcal E(u)
 =\frac{2}{N-1}\sum_{i<j}\|(-A_{ij})^{1/2}u\|_2^2,
 \qquad
 D(\mathcal E)=\bigcap_{i<j}D((-A_{ij})^{1/2}).
\]
This finite sum is a closed form on the indicated intersection.
It is densely defined because $C_c^\infty(\R^{3N})\subset D(\mathcal E)$; the coefficients $a(v_i-v_j)$ are locally integrable for $\gamma\in[-3,1]$.
Let $-\mathcal L^N$ denote its associated nonnegative self-adjoint operator.
For smooth compactly supported functions, integration by parts identifies its form with the differential expression \eqref{eq:LN-def}:
\[
 \mathcal E(u)
 =\frac{1}{2(N-1)}\sum_{i<j}
 \int_{\R^{3N}}
 (\nabla_{v_i}u-\nabla_{v_j}u)\cdot
 a(v_i-v_j)(\nabla_{v_i}u-\nabla_{v_j}u)\md\bm v.
\]
The monotone convergence theorem for closed forms and its functional-calculus consequence \cite[Theorems~3.1 and~4.2]{Simon1978MonotoneForms} give
\[
 P_t^{N,h}u\longrightarrow\ee^{t\mathcal L^N}u
 \quad\text{in }L^2\quad(t\ge0).
\]
The convergence is uniform on bounded time intervals.
Indeed, for $u\in D(\mathcal E)$ the spectral theorem gives
\[
 \|(P_t^{N,h}-P_s^{N,h})u\|_2^2
 \le2|t-s|\mathcal E_h(u)
 \le2|t-s|\mathcal E(u).
\]
A finite time grid therefore makes the pointwise convergence uniform for these $u$, and density of $D(\mathcal E)$ together with $L^2$ contraction gives the result for every $u\in L^2$.
Since $F_0^N\in L^1\cap L^\infty\subset L^2$, we conclude that
\[
 \sup_{t\le T}
 \|F_t^{N,h}-\ee^{t\mathcal L^N}F_0^N\|_2\to0.
\]
In particular, the curve constructed in Step~3 is $F_t^N=\ee^{t\mathcal L^N}F_0^N$.
Hence any limit of $\{F_{\cdot}^{N,h}\}_{h>0}$ must coincide with
	$F_{\cdot}^{N}$.
	Therefore the convergence \eqref{eq:fixedN-conv} holds for the full
	limit $h\to0$, and the stated a priori bounds follow from Step~1 by
	lower semicontinuity and the conservation argument in Step~3.
\end{proof}

\subsection*{Acknowledgements}
The author wishes to thank Lei Li for valuable discussions regarding the model, which is based on their previous collaboration.
This work was financially supported by the National Natural Science Foundation of China (12222103, 42450269), and LMNS at Fudan University.

\appendix

\section{Structural properties of the MBC process}
\label{sec:appendix-ergodicity}

In this appendix we justify the structural properties of the MBC dynamics.
We work throughout with fixed parameters $(N,h,E,P)$ and simplify notation
accordingly.

Recall from Section~\ref{sec:The-MBC-process} the microcanonical (Boltzmann sphere)
$\mathcal{M}:=\mathcal{M}^{N}_{E,P}$ defined in \eqref{eq:M-shell} and the normalized surface measure
$\mu_{\mathrm{mc}}$ on it.
By a Galilean shift we may, without loss of generality, assume $P=0$;
all statements below are invariant under this change of frame.
We also recall from~\eqref{disc:eq:MBC-generator} the operators
\[
\mathcal{Q}_{ij}=\mathcal{Q}^{h}_{ij},\quad
\mathcal{Q}=\mathcal{Q}^{N,h},\quad
\mathcal{L}=\mathcal{L}^{N,h},\quad
P_t=P_t^{N,h}.
\]

\begin{proposition}\label{thm:MBC-properties-app}
For the MBC process on $\mathcal M$ defined above, the following hold.
\begin{enumerate}
		\item[\emph{(i)}] \textup{(Invariance)}
		$\mu_{\mathrm{mc}}$ is invariant for $\mathcal{Q}$ and for $(P_{t})_{t\ge0}$:
		for all bounded measurable $F$ on $\mathcal{M}$ and all $t\ge0$,
		\[
		\mu_{\mathrm{mc}}(\mathcal{Q}F)=\mu_{\mathrm{mc}}(F),
		\qquad
		\mu_{\mathrm{mc}}(P_{t}F)=\mu_{\mathrm{mc}}(F).
		\]
		
		\item[\emph{(ii)}] \textup{(Reversibility)}
		$\mathcal{Q}$ and $P_{t}$ are $\mu_{\mathrm{mc}}$-reversible:
		for all bounded measurable $F,G$ on $\mathcal{M}$,
		\[
		\mu_{\mathrm{mc}}(F\,\mathcal{Q}G)=\mu_{\mathrm{mc}}(G\,\mathcal{Q}F),
		\qquad
		\mu_{\mathrm{mc}}(F\,P_{t}G)=\mu_{\mathrm{mc}}(G\,P_{t}F).
		\]
		Equivalently, $\mathcal{Q}$, $\mathcal{L}$ and $P_{t}$ are
		self-adjoint on $L^{2}(\mathcal{M},\mu_{\mathrm{mc}})$.

\item[\emph{(iii)}] \textup{(Feller regularity)}
$\mathcal Q$ is a Feller Markov operator on $\mathcal M$.
Consequently, $(P_t)_{t\ge0}$ is a Feller semigroup.
\item[\emph{(iv)}] \textup{(Uniqueness of absolutely continuous invariant measures)}
$\mu_{\mathrm{mc}}$ is the unique invariant probability measure of $\mathcal Q$, and hence of $(P_t)_{t\ge0}$, that is absolutely continuous with respect to $\mu_{\mathrm{mc}}$.
\item[\emph{(v)}] \textup{(Ergodicity and mixing)}
For any initial law on $\mathcal M$ that is absolutely continuous with respect to $\mu_{\mathrm{mc}}$, the law of $\bm v^N(t)$ converges in total variation to $\mu_{\mathrm{mc}}$ as $t\to\infty$.
Moreover, there exist constants $c>0$ and $\rho\in(0,1)$, depending on $(N,E,h,\gamma)$, such that
\[
 \|\mathcal Q^nF\|_{L^2(\mu_{\mathrm{mc}})}
 \le\rho^n\|F\|_{L^2(\mu_{\mathrm{mc}})}
 \quad\forall\,n\in\N,
\]
and
\[
 \|P_tF\|_{L^2(\mu_{\mathrm{mc}})}
 \le\ee^{-ct}\|F\|_{L^2(\mu_{\mathrm{mc}})}
 \quad\forall\,t\ge0,
\]
for all $F\in L^2(\mathcal M,\mu_{\mathrm{mc}})$ with $\mu_{\mathrm{mc}}(F)=0$.
\end{enumerate}
\end{proposition}

The proof below uses disintegration for invariance and reversibility, continuity of the pair update for the Feller property, and comparison with the Kac spectral gap for mixing.

\begin{corollary}\label{cor:whole-space}
Extend $\mathcal Q$, $\mathcal L$ and $(P_t)_{t\ge0}$ to bounded Borel functions on $\R^{3N}$ by the same formulas \eqref{disc:eq:MBC-generator}.
Then $\mathcal Q$ and $P_t$ are self-adjoint contractions on $L^2(\R^{3N})$ with respect to Lebesgue measure, and $\mathcal L$ is bounded and self-adjoint.
Every mixture of microcanonical measures is invariant for the whole-space dynamics:
\[
		\nu(\md\bm v)
		=\int \mu_{\mathrm{mc}}^{N}(E,P,\md\bm v)\,\pi(\md E,\md P),
		\]
where $\pi$ is any probability measure on $\{(E,P):2E>|P|^2\}$.
On the degenerate shell $2E=|P|^2$, the point mass at $v_1=\cdots=v_N=P$ is also invariant.
\end{corollary}

\begin{proof}
In the whole-space formulation, the pair change of variables \eqref{eq:CoV} gives
\[
 \md\bm v=r_{ij}^2\md\sigma(e_{ij})\md r_{ij}\md\bar v_{ij}\md\bm v_{[ij]},
\]
where $\sigma$ is surface measure on $\S^2$.
The update preserves all variables except the relative direction.
Fubini's theorem and symmetry of the spherical heat operator therefore give
\[
 \int(\mathcal Q_{ij}F)G\md\bm v
 =\int F(\mathcal Q_{ij}G)\md\bm v.
\]
Averaging and using the Poisson representation proves whole-space self-adjointness.
Jensen's inequality for the Markov kernel, followed by invariance of Lebesgue measure, gives the $L^2$ contractions.
The invariance of each mixture follows by integrating item \emph{(i)} of Proposition~\ref{thm:MBC-properties-app}; on a degenerate shell every update is the identity.
\end{proof}

We next prove a mild regularity result for the whole-space dynamics: absolute continuity of the initial law is preserved by the MBC process for all positive times, and its $L^\infty$-norm does not increase.
This provides the baseline regularity used in Section~\ref{sec:discrete-information}.

\begin{proposition}
\label{prop:a.c.bound}If the initial law $F_{0}^{N}$ on $\R^{3N}$ has a density with respect to the Lebesgue measure, then for every $t>0$ the law $F_{t}^{N}$ has a density (still denoted by $F_{t}^{N}$).
Moreover, if $F_{0}^{N}\in L^{\infty}(\R^{3N})$, then $F_{t}^{N}\in L^{\infty}(\R^{3N})$ and $\|F_{t}^{N}\|_{\infty}\le\|F_{0}^{N}\|_{\infty}$.
\end{proposition}

\begin{proof}
Write $\tau_{ij}=hr_{ij}^{\gamma}$ and $\lambda=N/h$.
For an unordered pair $\{i,j\}$, with respect to (\ref{eq:CoV}) we have
\[
\md\bm{v}=r_{ij}^{2}\md\sigma(e_{ij})\md r_{ij}\md\bar{v}_{ij}\md\bm{v}_{[ij]}.
\]
Let $F\in L^{1}(\R^{3N})$ be a density of a law $\nu$, i.e., $\nu(\md\bm{v})=F(\bm{v})\md\bm{v}$.
For any bounded test function $\Phi$, Fubini's theorem and the symmetry of the heat operator $P_t^{\S^2}$ give
\begin{align*}
\int_{\R^{3N}}(\mathcal{Q}^{N,h}\Phi)(\bm{v})F(\bm{v})\md\bm{v} & =\frac{2}{N(N-1)}\sum_{i<j}\int\md\bm{v}_{[ij]}\int\md\bar{v}_{ij}\int_{0}^{+\infty}r^{2}\md r\\
 & \quad\qquad\cdot\int_{\mathbb{S}^{2}}\big[P_{\tau_{ij}}^{\S^2}\Phi(\bm{v}_{[ij]},\bar{v}_{ij},r_{ij},\cdot)\big](e_{ij})\,F(\bm{v}_{[ij]},\bar{v}_{ij},r_{ij},e_{ij})\md\sigma(e_{ij})\\
 & =\frac{2}{N(N-1)}\sum_{i<j}\int\md\bm{v}_{[ij]}\int\md\bar{v}_{ij}\int_{0}^{+\infty}r^{2}\md r\\
 & \quad\qquad\cdot\int_{\mathbb{S}^{2}}\Phi(\bm{v}_{[ij]},\bar{v}_{ij},r_{ij},e_{ij})\,\big[P_{\tau_{ij}}^{\S^2}F(\bm{v}_{[ij]},\bar{v}_{ij},r_{ij},\cdot)\big](e_{ij})\md\sigma(e_{ij}).
\end{align*}
Since $\Phi$ is arbitrary, this identifies the adjoint on densities:
\begin{align*}
(\mathcal{Q}^{N,h})^{*}F(\bm{v}) & =\frac{2}{N(N-1)}\sum_{i<j}\big[P_{\tau_{ij}}^{\S^2}F(\bm{v}_{[ij]},\bar{v}_{ij},r_{ij},\cdot)\big](e_{ij})\\
 & =\frac{2}{N(N-1)}\sum_{i<j}\int_{\mathbb{S}^{2}}F(\bm{v}_{[ij]},\bar{v}_{ij},r_{ij},e')\,P_{\tau_{ij}}^{\S^2}(e_{ij},\md e')
\end{align*}
for a.e. $\bm{v}\in\R^{3N}$ (the set $\{r_{ij}=0\}$ is null with respect to both Lebesgue measure and the microcanonical surface measure).
This formula shows that $(\mathcal{Q}^{N,h})^{*}F$ is given pointwise by an integral of $F$ against a probability kernel on $\mathbb{S}^{2}$, hence $(\mathcal{Q}^{N,h})^{*}F$ is a Borel function.
Moreover, taking test function $\Phi\equiv\bm{1}$ gives
\[
(\mathcal{Q}^{N,h}\bm{1})(\bm{v})=\frac{2}{N(N-1)}\sum_{i<j}(P_{\tau_{ij}}^{\S^2}\bm{1})(e_{ij})=1,
\]
which yields
\[
\int_{\R^{3N}}(\mathcal{Q}^{N,h})^{*}F(\bm{v})\md\bm{v}=\int_{\R^{3N}}(\mathcal{Q}^{N,h}\bm{1})(\bm{v})F(\bm{v})\md\bm{v}=\int_{\R^{3N}}F(\bm{v})\md\bm{v}.
\]
Thus, the absolute continuity is preserved and the mass is conserved.

Furthermore, if $F\in L^{\infty}(\R^{3N})$, then for a.e. $\bm{v}\in\R^{3N}$, (\ref{disc:eq:MBC-generator}) gives
\[
|(\mathcal{Q}^{N,h})^* F(\bm{v})|=|\mathcal{Q}^{N,h}F(\bm{v})|\le\frac{2}{N(N-1)}\sum_{i<j}\int_{\mathbb{S}^{2}}\|F\|_{\infty}\,P_{\tau_{ij}}^{\S^2}(e_{ij},\md e')=\|F\|_{\infty},
\]
by using the self-adjointness of $\mathcal{Q}^{N,h}$ from Corollary~\ref{cor:whole-space}.
Thus $\|\mathcal{Q}^{N,h}F\|_{\infty}\le\|F\|_{\infty}$.
By iteration, $(\mathcal{Q}^{N,h})^{n}F_{0}^{N}$ is a density in $L^{\infty}(\R^{3N})$ with $\|(\mathcal{Q}^{N,h})^{n}F_{0}^{N}\|_{\infty}\le\|F_{0}^{N}\|_{\infty}$ for each $n\ge0$.
The series representation of the semigroup then implies
\[
F_{t}^{N}=(P_{t}^{N,h})^{*}F_{0}^{N}=P_{t}^{N,h}F_{0}^{N}={\rm e}^{-\lambda t}\sum_{n=0}^{\infty}\frac{(\lambda t)^{n}}{n!}(\mathcal{Q}^{N,h})^{n}F_{0}^{N},
\]
which is a pointwise absolutely convergent sum of $L^{\infty}$ functions.
Hence $F_{t}^{N}$ is a density and 
\[
\|F_{t}^{N}\|_{\infty}\le{\rm e}^{-\lambda t}\sum_{n=0}^{\infty}\frac{(\lambda t)^{n}}{n!}\|F_{0}^{N}\|_{\infty}=\|F_{0}^{N}\|_{\infty}.
\]
This concludes the proof.
\end{proof}

\begin{remark}\label{rem:equilibrium-consistency}
For fixed $(N,E,P)$, Proposition~\ref{thm:MBC-properties-app} shows exponential convergence in $L^2(\mathcal M,\mu_{\mathrm{mc}})$ for densities in this space.
By standard equivalence-of-ensembles arguments, the one-particle marginal of $\mu_{\mathrm{mc}}$ converges, as $N\to\infty$ with $E,P$ fixed, to the Maxwellian with the same mean momentum and energy per particle.
Together with our main results, this indicates that the long-time/large-$N$ behavior of the MBC process is consistent with the Landau dynamics converging to the corresponding Maxwellian.
\end{remark}

The rest of the appendix is devoted to a proof of Proposition~\ref{thm:MBC-properties-app}, based on the geometric structure of the microcanonical shell and the spherical heat kernel.

\subsection{Invariance and reversibility}

We prove items \emph{(i)} and \emph{(ii)} of Proposition~\ref{thm:MBC-properties-app}.
Fix an unordered pair $\{i,j\}$ and recall the notation $\bar v_{ij},r_{ij},e_{ij}$ from \eqref{eq:CoV}.
We use only the standard properties of the spherical heat semigroup $P_t^{\S^2}$: its kernel is mass-preserving and symmetric with respect to the uniform probability measure $\omega_{\S^2}$:
\[
 \int_{\S^2}P_t^{\S^2}(e,\md e')=1,\qquad
 P_t^{\S^2}(e,\md e')\,\omega_{\S^2}(\md e)
 =P_t^{\S^2}(e',\md e)\,\omega_{\S^2}(\md e').
\]
Disintegrate $\mu_{\mathrm{mc}}$ with respect to $y=(\bm v_{[ij]},\bar v_{ij},r_{ij})$.
Rotational symmetry of the shell makes the conditional direction uniform on $\S^2$ whenever $r_{ij}>0$.
Writing $\eta(\md y)$ for the marginal law of $y$, Fubini's theorem and preservation of the uniform angular measure give
\[
\begin{split}
 \mu_{\mathrm{mc}}(\mathcal Q_{ij}F)
 &=\int\eta(\md y)\int_{\S^2}\int_{\S^2}
 F(y,e')P_{hr_{ij}^\gamma}^{\S^2}(e,\md e')\,\omega_{\S^2}(\md e)\\
 &=\mu_{\mathrm{mc}}(F).
\end{split}
\]
The zero-relative-speed update is the identity and causes no change.
Averaging over $i<j$ yields $\mu_{\mathrm{mc}}(\mathcal QF)=\mu_{\mathrm{mc}}(F)$, so $\mu_{\mathrm{mc}}$ is invariant for $\mathcal Q$.
Since $P_t=\exp(t\lambda(\mathcal Q-I))$, with $\lambda=N/h$, $\mu_{\mathrm{mc}}$ is also invariant for $(P_t)_{t\ge0}$.
This proves item \emph{(i)}.

Let $F,G$ be bounded Borel functions on $\mathcal M$.
For each $\{i,j\}$,
\[
 \mu_{\mathrm{mc}}(F\,\mathcal Q_{ij}G)
 =\int\eta(\md y)\int_{\S^2}\int_{\S^2}
 F(y,e)G(y,e')P_{hr_{ij}^\gamma}^{\S^2}(e,\md e')\,\omega_{\S^2}(\md e).
\]
Using symmetry in $(e,e')$ and exchanging $F,G$ gives
\[
\mu_{\mathrm{mc}}(F\,\mathcal{Q}_{ij}G)
=\mu_{\mathrm{mc}}(G\,\mathcal{Q}_{ij}F).
\]
Summing over $i<j$ yields
$\mu_{\mathrm{mc}}(F\,\mathcal{Q}G)
=\mu_{\mathrm{mc}}(G\,\mathcal{Q}F),$
equivalently,
\[
\mathcal{Q}(\bm v,\md\bm v')\,\mu_{\mathrm{mc}}(\md\bm v)
=\mathcal{Q}(\bm v',\md\bm v)\,\mu_{\mathrm{mc}}(\md\bm v').
\]
Thus $\mathcal{Q}$ is $\mu_{\mathrm{mc}}$-reversible.
Since $\mathcal{L}=\lambda(\mathcal{Q}-I)$ and
$P_{t}=\exp(t\mathcal{L})$ are functions of $\mathcal{Q}$,
they are also self-adjoint on $L^{2}(\mathcal{M},\mu_{\mathrm{mc}})$.
This proves item \emph{(ii)} of
Proposition~\ref{thm:MBC-properties-app}.

\subsection{Feller regularity and uniqueness of absolutely continuous invariant measures}

Since $P_t=\exp(t\lambda(\mathcal Q-I))$, a probability measure is invariant for $(P_t)_{t\ge0}$ if and only if it is invariant for the embedded kernel $\mathcal Q$.
It is therefore enough to work with $\mathcal Q$.
By \eqref{eq:Q_ij}, for continuous $F$ the map $\bm v\mapsto\mathcal Q_{ij}F(\bm v)$ is continuous away from $r_{ij}=0$, by continuity of the spherical heat kernel and of the pair coordinates.
At $r_{ij}=0$, every possible pair update moves the configuration by at most $\sqrt2r_{ij}$.
Uniform continuity of $F$ on the compact shell therefore gives $\mathcal Q_{ij}F(\bm v)-F(\bm v)\to0$ as $r_{ij}\to0$.
Hence $\mathcal Q$ is Feller, being a finite average of these pair operators.
The uniformly convergent Poisson series proves the Feller property for $(P_t)_{t\ge0}$, giving item \emph{(iii)}.
The uniqueness assertion in item \emph{(iv)} follows from the convergence proved below.

\subsection{Ergodicity and $L^2$-mixing}

We prove item \emph{(v)} of Proposition~\ref{thm:MBC-properties-app} by comparison with the spectral gap for Kac's process \cite[Theorem~1.2]{CarlenCarvalhoLoss2020}.
Let $\Pi_{ij}F$ be the uniform average of $F$ over the relative direction, with all other pair coordinates fixed.
The first nonzero eigenvalue of $-\Delta_{\S^2}$ is $2$.
The spherical harmonic expansion on each conditional sphere therefore gives
\[
 \langle F,(I-\mathcal Q_{ij})F\rangle
 \ge\int_{\mathcal M}(1-\ee^{-2hr_{ij}^\gamma})
 |F-\Pi_{ij}F|^2\md\mu_{\mathrm{mc}}.
\]
On the fixed shell, $r_{ij}\le R:=\sqrt{2N(2E-|P|^2)}$.
Set $\alpha=\max\{\gamma,0\}\in[0,1]$.
If $\gamma\le0$, then $1-\ee^{-2hr^\gamma}\ge1-\ee^{-2hR^\gamma}>0$ for $r>0$.
If $\gamma>0$, concavity of $x\mapsto1-\ee^{-2hx}$ gives
\[
 1-\ee^{-2hr^\gamma}
 \ge\frac{1-\ee^{-2hR^\gamma}}{R^\gamma}r^\gamma
 \qquad(0\le r\le R).
\]
Thus there is $a>0$, depending only on the fixed parameters, such that
\[
 \langle F,(I-\mathcal Q)F\rangle
 \ge a\binom N2^{-1}\sum_{i<j}
 \int_{\mathcal M}r_{ij}^\alpha|F-\Pi_{ij}F|^2\md\mu_{\mathrm{mc}}.
\]
Up to the factor $N$ in its jump rate, the right-hand side is the Dirichlet form of the momentum- and energy-conserving Kac process with uniform angular scattering and collision rate $r^\alpha$.
The cited spectral-gap theorem, after translating and rescaling the shell, yields $\delta>0$ such that
\[
 \langle F,(I-\mathcal Q)F\rangle\ge\delta\|F\|_2^2
 \qquad\text{if }\mu_{\mathrm{mc}}(F)=0.
\]
Each pair heat operator is nonnegative and self-adjoint, hence $0\le\mathcal Q\le I$.
By the spectral theorem, the spectrum of $\mathcal Q$ on the orthogonal complement
\[
L^{2}_{0}(\mu_{\mathrm{mc}})
:=\big\{F\in L^{2}(\mathcal{M},\mu_{\mathrm{mc}}):
\mu_{\mathrm{mc}}(F)=0\big\},
\]
is contained in $[0,1-\delta]$.
Decreasing $\delta$ if necessary, we may assume $0<\delta<1$.
Hence there exists $\rho\in(0,1)$ such that
\[
 \|\mathcal Q^nF\|_{L^2(\mu_{\mathrm{mc}})}
 \le\rho^n\|F\|_{L^2(\mu_{\mathrm{mc}})}
 \qquad\forall\,n\in\N,\quad F\in L^2_0(\mu_{\mathrm{mc}}).
\]
For the continuous-time semigroup, recall
$\mathcal{L}=\lambda(\mathcal{Q}-I)$ and
$P_{t}=\exp(t\mathcal{L})$.
By the spectral theorem,
the spectrum of $\mathcal{L}$ on $L^{2}_{0}(\mu_{\mathrm{mc}})$
is contained in $(-\infty,-c]$ for some $c>0$
corresponding to the above spectral gap of $\mathcal{Q}$.
Thus
\[
\|P_{t}F\|_{L^{2}(\mu_{\mathrm{mc}})}
\le\ee^{-ct}\,\|F\|_{L^{2}(\mu_{\mathrm{mc}})}
\qquad\forall\,t\ge0,\ \forall\,F\in L^{2}_{0}(\mu_{\mathrm{mc}}),
\]
giving the claimed $L^{2}$-mixing rate.

Finally, let $\nu$ be any initial probability measure on $\mathcal M$ with density $g\in L^1(\mu_{\mathrm{mc}})$.
Approximate $g$ in $L^1$ by probability densities $g_k\in L^2(\mu_{\mathrm{mc}})$, using truncation and normalization.
The $L^1$ contraction and the preceding $L^2$ estimate give
\[
 \|\mathcal Q^ng-1\|_1
 \le\|g-g_k\|_1+\|\mathcal Q^ng_k-1\|_2\longrightarrow0,
\]
first letting $n\to\infty$ and then $k\to\infty$.
The same argument applies to $P_t$.
This proves convergence in total variation for every absolutely continuous initial law.
If $\nu$ is invariant, its law is unchanged by the dynamics, so the limit forces $\nu=\mu_{\mathrm{mc}}$.
This also proves item \emph{(iv)} of Proposition~\ref{thm:MBC-properties-app}.

\begin{remark}\label{rem:gap-N}
The constants $\rho\in(0,1)$ and $c>0$ in Proposition~\ref{thm:MBC-properties-app}\,\emph{(v)} may depend on $(N,E,h,\gamma)$.
The comparison above shows existence of a positive gap for each fixed $N$, but does not provide uniform control as $N\to\infty$.
We do not address uniform-in-$N$ mixing bounds here.
\end{remark}

\section{Auxiliary information estimates}
\label{sec:auxiliary-information}

This appendix contains the variational and affinity arguments for directional information and the angular regularization used in the fourth-power estimate.

\begin{proof}[Proof of Lemma~\ref{disc:prop:active-calculus}]
Since $X^m\psi$ and $\psi^2$ are continuous and compactly supported, \eqref{disc:eq:FXm} is a supremum of continuous affine functions of $P$, hence convex and weakly lower semicontinuous.

For marginal contraction, let $\psi$ be a compactly supported smooth test on the coordinates retained in $Q$.
Choose smooth compactly supported cutoffs $0\le\zeta_R\le1$ in the coordinates being integrated out, with $\zeta_R\to1$ pointwise, and test \eqref{disc:eq:FXm} with $\psi\zeta_R$.
Since $X$ acts only on retained coordinates,
\[
 X\zeta_R=0,\qquad X^m(\psi\zeta_R)=\zeta_RX^m\psi.
\]
By the dominated convergence theorem, the linear and quadratic terms converge to $2(-1)^m\int X^m\psi\md Q$ and $\int\psi^2\md Q$, respectively.
Each test value for $Q$ is therefore bounded by $\mathcal{K}_{X,m}(P)$, and taking the supremum proves the contraction.

If $\mathcal{K}_{X,m}(P)<\infty$, substituting $\lambda\psi$ in \eqref{disc:eq:FXm} and optimizing over $\lambda\in\R$ gives
\[
 \Big|\int X^m\psi\md P\Big|
 \le \mathcal{K}_{X,m}(P)^{1/2}\|\psi\|_{L^2(P)}.
\]
The linear functional $\psi\mapsto(-1)^m\int X^m\psi\md P$ vanishes on tests that are zero $P$-almost everywhere and is bounded on the dense subspace $C_c^\infty\subset L^2(P)$.
Its continuous extension and the Riesz representation theorem give the unique $S\in L^2(P)$ satisfying $X^mP=SP$.

Conversely, if such an $S$ exists, completing the square in the variational formula gives
\[
 2(-1)^m\int X^m\psi\md P-\int\psi^2\md P
 =\int|S|^2\md P-\int|\psi-S|^2\md P.
\]
Taking the supremum over the dense class $C_c^\infty\subset L^2(P)$ proves both finiteness and the norm identity \eqref{disc:eq:Riesz-score}.
For a smooth positive density $p$, integration by parts identifies the distributional derivative as $X^mp$, and the equivalence above gives the stated formula, including infinite values.
\end{proof}

\begin{proof}[Proof of Lemma~\ref{disc:thm:score-affinity}]
Let $\nu$ be a random probability measure with law $\Pi$.
Given $\nu$, sample $V_1,V_2,\ldots$ independently with common law $\nu$.
Then $\Pi^j=\Law(V_1,\ldots,V_j)=\E[\nu^{\otimes j}]$.
We take expectations over both $\nu$ and the samples, and set $\mathcal F_j:=\sigma(V_1,\ldots,V_j)$.

By Lemma~\ref{disc:prop:active-calculus}, $\mathcal{K}_{X,m}(\Pi^j)$ increases with $j$.
Convexity and the definition of $\mathcal{K}_{X,m}$ give
\begin{equation}
 \mathcal{K}_{X,m}(\Pi^j)
 \le\E\mathcal{K}_{X,m}(\nu^{\otimes j})
 =\E\mathcal{K}_{X,m}(\nu^{\otimes k_0}).
 \label{disc:eq:affinity-easy}
\end{equation}
If the information is infinite, marginal contraction gives the same equality.
If the limit on the left is infinite, \eqref{disc:eq:affinity-easy} proves the result; we therefore assume this limit is finite.

For the reverse inequality, we construct an $L^2$-bounded martingale.
By Lemma~\ref{disc:prop:active-calculus}, write $X^m\Pi^j=S_j\Pi^j$ with $S_j\in L^2(\Pi^j)$.
Below, $S_j$ also denotes the random variable $S_j(V_1,\ldots,V_j)$.
Then $\E|S_j|^2=\mathcal{K}_{X,m}(\Pi^j)$ is uniformly bounded.
Recall that $\Pi^j$ is the marginal of $\Pi^{j+1}$ on the first $j$ coordinates.
For $\phi\in C_c^\infty(\R^{3j})$, use a cutoff in the added coordinate to obtain
\[
 \E[S_{j+1}\phi(V_{1:j})]
 =(-1)^m\E[(X^m\phi)(V_{1:j})]
 =\E[S_j\phi(V_{1:j})].
\]
By density in $L^2(\Pi^j)$, it follows that $\E[S_{j+1}\mid\mathcal F_j]=S_j$.
Thus $(S_j)$ is an $L^2$-bounded martingale, and the martingale convergence theorem gives
\begin{equation}
 S_j\longrightarrow S_\infty\quad\text{in }L^2,
 \qquad
 \E|S_\infty|^2=\lim_{j\to\infty}\mathcal{K}_{X,m}(\Pi^j).
 \label{disc:eq:score-limit}
\end{equation}

For $j>k_0$, the conditional strong law of large numbers gives
\[
 \eta_j:=\frac1{j-k_0}\sum_{\ell=k_0+1}^j\delta_{V_\ell}
 \rightharpoonup\nu\quad\text{almost surely}.
\]
For test functions $\psi\in C_c^\infty(\R^{3k_0})$ and $g\in C_b(\Pp(\R^3))$, we have
\[
 \E[S_j\psi(V_{1:k_0})g(\eta_j)]
 =(-1)^m\E[(X^m\psi)(V_{1:k_0})g(\eta_j)].
\]
Let $j\to\infty$.
By the $L^2$ convergence of $S_j$ and the bounded convergence of $g(\eta_j)$, we obtain
\begin{equation}
 \E[S_\infty\psi(V_{1:k_0})g(\nu)]
 =(-1)^m\E[(X^m\psi)(V_{1:k_0})g(\nu)].
 \label{disc:eq:affinity-limit-test}
\end{equation}
We choose a jointly measurable function $\bar S$ such that
\[
 \bar S(\nu,V_{1:k_0})
 =\E[S_\infty\mid\nu,V_1,\ldots,V_{k_0}].
\]
Since $\psi$ and $g$ are arbitrary, conditioning \eqref{disc:eq:affinity-limit-test} on $\nu$ gives, for $\Pi$-almost every $\rho$,
\[
 \bar S(\rho,\cdot)\,\rho^{\otimes k_0} = X^m\rho^{\otimes k_0}
 \quad\text{in }\mathcal D'.
\]
Jensen's inequality ensures that $\bar S(\rho,\cdot)\in L^2(\rho^{\otimes k_0})$ for $\Pi$-almost every $\rho$.
Lemma~\ref{disc:prop:active-calculus} and \eqref{disc:eq:score-limit} therefore give
\[
 \E\mathcal{K}_{X,m}(\nu^{\otimes k_0})
 =\E|\bar S(\nu,V_{1:k_0})|^2
 \le\E|S_\infty|^2
 =\lim_{j\to\infty}\mathcal{K}_{X,m}(\Pi^j).
\]
Together with \eqref{disc:eq:affinity-easy}, this proves \eqref{disc:eq:score-affinity}.
\end{proof}

\begin{proof}[Proof of Lemma~\ref{stat:lem:active-difference}]
Let $P=\rho^{\otimes2}=F\md v\md w$, and fix one annular field $Y=\chi_n(r)X_k$ from Remark~\ref{disc:rem:singular-affinity}.
Its first distributional derivative is represented by
\[
 U(v,w)=\chi_n(r)r^{\gamma/4-1/2}
 e_k\cdot\big[(v-w)\times\big(\nabla\log \rho(v)-\nabla\log \rho(w)\big)\big].
\]
The annular cutoff and $I(\rho)<\infty$ give $U\in L^2(P)$.
Lemma~\ref{disc:prop:active-calculus} also gives $Y^2P=VP$ with $V\in L^2(P)$ and $\int|V|^2\md P=\mathcal{K}_{Y,2}(P)$.
We claim that
\begin{equation}
 \int|U|^4\md P\le\frac94\mathcal{K}_{Y,2}(P).
 \label{stat:eq:weak-score-fourth}
\end{equation}
For smooth positive densities, this bound follows from \cite[Proposition~2.1, with $\beta=1$]{Tabary2025_PoC_Landau_Coulomb}.

To justify the integration by parts under the stated finite-information assumptions, set
\[
 B_k:=b_k(v-w)\cdot(\nabla_v-\nabla_w).
\]
Let $\Phi_\theta$ be its periodic, volume-preserving flow, and let $h_\varepsilon$ be the strictly positive heat kernel on $\R/(2\pi\mathbb Z)$.
The flow fixes $(v+w)/2$ and rotates $v-w$ by angle $2\theta$ about $e_k$, so $2\pi$ is a period.
Define the angular convolution
\[
 F_\varepsilon(x)=\int_0^{2\pi}h_\varepsilon(\theta)
 F(\Phi_{-\theta}x)\md\theta,
 \qquad x=(v,w).
\]
Then $F_\varepsilon\to F$ in $L^1$.
Differentiating along the flow in distributions and integrating by parts in the periodic parameter gives
\[
 B_k^jF_\varepsilon(x)
 =\int_0^{2\pi}h_\varepsilon^{(j)}(\theta)
 F(\Phi_{-\theta}x)\md\theta,
 \qquad j=1,2.
\]
The bounded coefficient $\chi_n(r)r^{\gamma/4-1/2}$ is invariant under this flow, so angular convolution commutes with $Y$ and $Y^2$.
In particular,
\[
 Y^2F_\varepsilon(x)=\int_0^{2\pi}h_\varepsilon(\theta)
 (VF)(\Phi_{-\theta}x)\md\theta,
\]
and Jensen's inequality followed by volume preservation gives
\[
 \int\frac{|Y^2F_\varepsilon|^2}{F_\varepsilon}
 \le\int|V|^2\md P.
\]
The ratios $|h_\varepsilon'|/h_\varepsilon$ and $|h_\varepsilon''|/h_\varepsilon$ are bounded on the circle.
Consequently, both $Y\log F_\varepsilon=YF_\varepsilon/F_\varepsilon$ and $Y^2F_\varepsilon/F_\varepsilon$ are bounded, with both quantities set to zero on orbits of zero mass.
Along the other orbits, angular convolution is positive and smooth and gives
\[
 Y^2\log F_\varepsilon
 =\frac{Y^2F_\varepsilon}{F_\varepsilon}-(Y\log F_\varepsilon)^2.
\]
It follows that
\[
 Y\big(F_\varepsilon(Y\log F_\varepsilon)^3\big)
 =3(Y\log F_\varepsilon)^2Y^2F_\varepsilon
 -2F_\varepsilon(Y\log F_\varepsilon)^4.
\]
Thus $F_\varepsilon(Y\log F_\varepsilon)^3$ and its directional derivative are integrable.
Integrate $Y\big(F_\varepsilon(Y\log F_\varepsilon)^3\big)$ against a compactly supported cutoff of $|v|^2+|w|^2$.
This cutoff is invariant under $Y$, so there is no boundary term.
Letting the cutoff increase to one and applying the Cauchy--Schwarz inequality gives
\[
\begin{aligned}
 2\int F_\varepsilon(Y\log F_\varepsilon)^4
 &=3\int(Y\log F_\varepsilon)^2Y^2F_\varepsilon,\\
 \int F_\varepsilon(Y\log F_\varepsilon)^4
 &\le\frac94\int\frac{|Y^2F_\varepsilon|^2}{F_\varepsilon}.
\end{aligned}
\]
Finally, use the variational formula
\[
 \int|U|^4\md P
 =\sup_{\zeta\in C_c^\infty}
 \Big\{-4\int Y\zeta\md P
       -3\int|\zeta|^{4/3}\md P\Big\}
\]
to obtain lower semicontinuity under $L^1$ convergence.
Letting $\varepsilon\downarrow0$ then proves \eqref{stat:eq:weak-score-fourth}.

Sum \eqref{stat:eq:weak-score-fourth} over $k$, use $|a|^4\le3\sum_k|a_k|^4$ and $\langle v-w\rangle^{\gamma-2}\le|v-w|^{\gamma-2}$, and increase the nested cutoffs $\chi_n$ to one.
Monotone convergence proves \eqref{stat:eq:active-difference}; the diagonal is $P$-null because $\rho$ is a density.
\end{proof}

\bibliographystyle{amsalpha}
\bibliography{landau-ref}

\end{document}